\documentclass{article}
\usepackage[utf8]{inputenc}
\usepackage{amsmath, amsthm, epsfig,amssymb, multicol, tikz}
\usepackage{epsfig}
\usepackage{listings}
\usepackage{tikz-cd}
\usepackage{esvect}
\usetikzlibrary{decorations}
\usepackage{todonotes}

\usepackage{xcolor}
\usepackage{lipsum}
\usepackage{amsfonts,amsmath,amssymb,amsthm}
\usepackage[shortlabels]{enumitem}
\usepackage{graphicx}
\usepackage{float}
\usepackage{url}

\usepackage[pagewise]{lineno}
\usepackage{todonotes}
\usepackage{setspace}
\usepackage{lineno}

\usepackage{geometry}
\newtheorem{lemma}{Lemma}
\newtheorem{theorem}{Theorem}
\newtheorem{corollary}{Corollary}
\newtheorem{definition}{Definition}

\newtheorem{remark}{Remark}

\usepackage{datetime}

\date{}
\newcommand{\rood}[1]{}

\title{Ricci-flat graphs with maximum degree at most $4$}
\author{
Shuliang Bai \thanks{Southeast University ({\tt sbai@seu.edu.cn}). }
\and 
Linyuan Lu
\thanks{University of South Carolina, Columbia, SC 29208
({\tt lu@math.sc.edu}). Supported in part by NSF
grant DMS 1600811 and ONR grant N00014-17-1-2842.}
\and 
Shing-Tung Yau
\thanks{Harvard University, Cambridge, MA 02138 ({\tt yau@math.harvard.edu}).
Supported by NSF
grant DMS-1607871: Analysis, Geometry and Mathematical Physics
 and National Science Foundation DMS-1418252: Collaborative Research: Geometric Analysis for Computer and Social Networks. }}

\begin{document}
\maketitle

\begin{abstract}
A graph is called Ricci-flat if its Ricci curvatures vanish on all edges, here the  definition of Ricci curvature on graphs was given by
Lin-Lu-Yau \cite{LLY}. The authors in \cite{LLY2} and \cite{errorcrect}  obtained a complete characterization for all Ricci-flat graphs with girth at least five. In this paper, we completely determined all Ricci-flat graphs with maximum degree at most $4$. 
\end{abstract}

\textsl{keywords:} Ricci curvature, 
Ricci-flat graph, Maximum degree
\\

\section{Introduction}
There is an increasing interest in  applying tools and ideas from continuous geometry to
discrete setting such as graphs.  One of the principal developments in this area concerns curvature for graphs. 
It is known that Ricci curvature plays a very important role on geometric analysis on Riemannian manifolds. As for graphs, 
the first definition of Ricci curvature was introduced by Fan
Chung and Yau\cite{Fchung} in 1996, their definition provides a curvature at each vertex. 
In 2009, 
Ollivier\cite{Ollivier} gave a notation of coarse Ricci curvature of Markov chains valid on arbitrary metric spaces, including graphs. His definition on graphs provides a curvature on each edge and depends on a so-called idleness parameter. 
For a more general definition of Ricci curvature, Lin and Yau \cite{LYau}
gave a generalization of lower Ricci curvature bound in the framework of
graphs in term the notation of Bakry and Emery.
In 2011, Lin-Lu-Yau \cite{LLY} modified Ollivier's definition. The modified version is a more suitable definition for graphs. 

By Lin-Lu-Yau's definition, a  Ricci-flat graph is a graph where Ricci curvature
vanishes on every edge.
In Riemannian manifold, there are many works on constructing Calabi-Yau manifolds which is a class of manifolds with zero Ricci curvature. 
As an analog,  we want to know what do Ricci-flat graphs look like? There have been works on classifying Ricci-flat graphs. At first,  \cite{LLY2} and \cite{errorcrect} classified all Ricci-flat graphs with girth at least five. Then the authors in \cite{HLYY} characterized all Ricci-flat graphs of girth four with vertex-disjoint $4$-cycles, their results show that there are two such graphs. While the fact is there are infinitely many Ricci-flat graphs with girth three or four.  In this paper,  we will completely classify these Ricci-flat graphs with maximum vertex degree at most $4$. 

Throughout this paper, let $G=(V, E)$ represent an undirected connected graph with vertex set $V$ and edge set $E$ without multiple edges or self loops. For any vertices $x, y\in V$, let  $d(x)$ denote the degree of vertex $x$, $d(x, y)$ denote the distance from $x$ to $y$, i.e. the length of the shortest path from $x$ to $y$. 
Denote $\Gamma(x)$ as the set of vertices that are adjacent to $x$, and $N(x)=\Gamma(x)\cup \{x\}$. Notation $x\sim y$ represent that two vertices $x$ and $y$ are adjacent and $(x, y)$ represent the edge. Let $C_3, C_4, C_5$ represent any  cycle of length $3, 4, 5$ respectively. 
Let $\mathcal{G}$ be the set of simple graphs with maximum degree at most $4$ that contains at least one copy of $C_3$ or $C_4$. Since the Ricci-flat graphs with girth at least five have been completely determined \cite{LLY}, we will only need to find the Ricci-flat graphs from the class $\mathcal{G}$ for our purpose. 

\begin{definition}\label{probabilitydistribution}
 A probability distribution over the vertex set $V$ is a mapping $\mu: V\to [0,1]$ satisfying $\sum_{x \in V} \mu (x)=1$. Suppose that two probability distributions $\mu_1$ and $\mu_2$ have finite support. A coupling between $\mu_1$ and $\mu_2$ is a mapping $A: V\times V\to [0, 1]$ with finite support so that 
$$\sum\limits_{y \in V} A(x, y)=\mu_1(x) \ \text{and} \sum\limits_{x \in V} A(x, y)=\mu_2(y). $$

The transportation distance between two probability distributions $\mu_1$ and $\mu_2$ is defined as follows:
$$W(\mu_1, \mu_2)=\inf\limits_{A} \sum\limits_{x, y\in V} A(x, y)d(x, y),$$
where the infimum is taken over all coupling $A$ between $\mu_1$ and $\mu_2$. 
\end{definition}

A coupling function provides a lower bound for the transportation distance, the following definition can provide an upper bound for the transportation distance. 
\begin{definition}\label{Lipschitz}
Let $G=(V, E)$ be a locally finite graph. Let $f: V \to \mathbb{R}$. We say $f$ is 
1-Lipschitz if 
$$f(x)-f(y)\leq d(x, y),$$
for each $x, y\in V$. 
\end{definition}
By the duality theorem of a linear optimization problem, the transportation distance can also be written as follows:
$$W(\mu_1, \mu_2)=\sup\limits_{f} \sum\limits_{x\in V} f(x)[\mu_1(x)-\mu_2(x)],$$
where the supremum is taken over all 1-Lipschitz functions $f$.

For any vertex $x\in V$ and any value $\alpha \in [0, 1]$,  the  probability distribution $\mu_x^{\alpha}$ is defined as:
\[
\mu_x^{\alpha}(z)=
\begin{cases}
 \alpha,  & \text{if $z=x$}, \\
 \frac{1-\alpha}{d(x)},  & \text{if $z\sim x$},\\
 0,& \text{otherwise}. 
\end{cases}
\]

For any $x, y\in V$, the $\alpha$-Ricci curvature $k_{\alpha}$ is defined to be 
$$k_{\alpha}(x, y)=1-\frac{W(\mu_x^{\alpha}, \mu_y^{\alpha})}{d(x, y)}. $$

Then the Ollivier-Ricci curvature $k(x, y)$ is defined by Lin-Lu-Yau as 
$$k(x, y)=\lim\limits_{\alpha\to 1} \frac{k_{\alpha}(x, y)}{1-\alpha}. $$

In \cite{BCLMP}, the authors introduced the concept of idleness function and studied its several properties. 
In the $\alpha$-Ollivier-Ricci curvature,  for every edge $(x, y)$ in $G=(V, E)$, the value $\alpha$ is called the {\it idleness},  and function $\alpha \to k_{\alpha}(x, y) $ is called the {\it Ollivier-Ricci idleness function}. 

\begin{theorem}
Let $G=(V, E)$ be a locally finite graph. Let $x, y \in V$ such that $x\sim y$ and $d(x)\geq d(y)$. Then 
$\alpha \to k_{\alpha}(x, y) $ is a piece-wise linear function over $[0, 1]$ with at most 3 linear parts. 
Furthermore, $k_{\alpha}(x, y)$ is linear on $[0, \frac{1}{lcm(d(x), d(x))+1}]$ and is also linear on 
$[\frac{1}{\max(d(x), d(x))+1}, 1]$. Thus, if 
we have further condition $d(x)=d(y)$, then $k_{\alpha}(x, y)$ has at most two linear parts. 
\end{theorem}

By above theorem, to study the local structure of an edge $(x, y)$ such that $k(x, y)=0$, we only need to consider 
$$k_{\alpha}(x, y)=0,  \quad \text{for $\alpha=\frac{1}{\max(d(x), d(x))+1}$}, $$
equivalently, 
$$W(\mu_x^{\alpha}, \mu_y^{\alpha})=1, \quad \text{for $\alpha=\frac{1}{\max(d(x), d(y))+1}$}.$$

Here are some helpful lemmas. 
\begin{lemma}\label{lemma1}\cite{LLY}
Suppose that an edge $(x, y)$ in a graph $G$ is not in any $C_3, C_4$ or $C_5$. Then $k(x, y)=\frac{2}{d(x)}+\frac{2}{d(y)}-2$. 
\end{lemma}

\begin{corollary}\label{noleaf}\cite{LLY}
Suppose that $x$ is a leaf-vertex(i.e. $d(x)=1$). Let $y$ be the only neighbor of $x$. Then $k(x, y)>0$.
\end{corollary}

\begin{lemma}\label{lemma2}\cite{LLY}
Suppose that an edge $(x, y)$ in a graph $G$ is not in any $C_3$ or $C_4$. Then $k(x, y)\leq\frac{1}{d(x)}+\frac{2}{d(y)}-1$. 
\end{lemma}

For any edge $(x, y)$ in a Ricci-flat graph $G$, we require $k(x, y)=0$,   by Lemma \ref{lemma1}, 
 if $0=k(x, y)\neq \frac{2}{d(x)}+\frac{2}{d(y)}-2$, then $(x, y)$ must be either in $C_3$ or $C_4$ or $C_5$. Similarly, by Lemma \ref{lemma2},
if $0=k(x, y)>\frac{1}{d(x)}+\frac{2}{d(y)}-1$, then $(x, y)$ must be either in $C_3$ or $C_4$; by Corollary \ref{noleaf},  there is no leaf-vertex in a Ricci-flat graph. 

We first study the local structure in any Ricci-flat graph, that is, for every edge $(x, y)\in E(G)$, analysis the distance between each pair of vertices in $\Gamma(x)\times \Gamma(y)$, then take all``good" structures into consideration to construct Ricci-flat graphs in class $\mathcal{G}$. All Ricci-flat graphs with maximum degree at most $4$ are showing in Theorems \ref{thm:44C3}, \ref{vertex-disjoint},  \ref{thm:C424}, \ref{thm:C444edge},  and \ref{thm:44noedgesharec4Typec}. 

\section{Study of local structures in $\Gamma(x)\cup \Gamma(y)$}
In this section, let $(x, y)$ be an edge in any Ricci-flat graph, we will study the distance between any pair of vertices in $\Gamma(x)\cup \Gamma(y)$. Let us first recall the local structure when edges are not contained in any $C_3$ or $C_4$.

\begin{lemma}\cite{LLY}
Suppose that an edge $(x, y)$ in a graph $G$ is not in any $C_3$ or $C_4$.
Without loss of generality, we assume $d(x) \geq d(y)$. If $k(x, y) = 0$, then one of the
following statements holds:
\begin{enumerate}
\item $d(x)=d(y)=2$. In this case, $(x, y)$ is not in any $C_5$. 
\item $d(x)=d(y)=3$. In this case, $(x, y)$  is shared by two $C_5$s.
\item $d(x)=2, d(y)=2$. In this case, let $y_1$ be the other neighbor of $x$ other than
$y$. Let $x_1,  x_2$ be two neighbors of $x$ other than $y$. Then $\{d(x_1, y_1),  d(x_2, y_1)\} =\{2, 3\}.$

\item $d(x) = 4$ and $d(y) = 2$. In this case, let $y_1$ be the other neighbor of $y$ other
than $x$. Let $x_1, x_2, x_3$ be three neighbors of $x$ other than $y$. Then at least
two of $x_1, x_2, x_3$  have distance $2$ from $x$.

\end{enumerate}
\begin{figure}[H]
\centering
\includegraphics[scale=0.4]{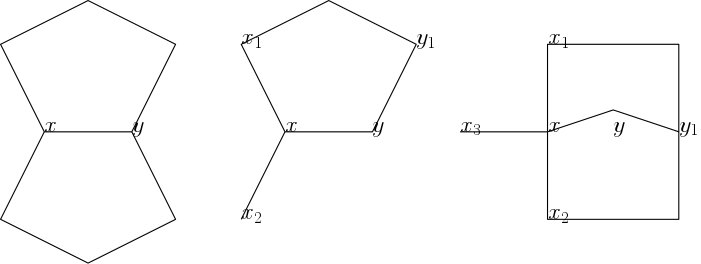}
\end{figure}
\end{lemma}

In the following, we consider the case when an edge is supported on $C_3$ or $C_4$ or both.
For convenience, we label all vertices by nonnegative integers, where $x$ and $y$ are labeled  by $0, 1$ respectively. Vertices in $\Gamma(x)$ are labeled by first $d(x)$ positive integers, vertices in $\Gamma(y)$ are labeled by the succeeding integers. 
Note in \cite{HLYY}, the authors also analyzed the local structures of edge with Ricci curvature $0$, the difference is that  their conclusions are based on that the graph has  girth $4$ and  the $4$-cycles in the graph  are mutually vertex-disjoint. Our results on these local structures can be applied to graphs without any restriction. 


\subsection{$d(x)=2, d(y)=2$}
%

For the case $k(x, y) = 0$ with $d(x)=d(y)=2$, by Definition \ref{probabilitydistribution},  shortening the distance for any pair of vertices would not increase the value of $W(\mu_x^{\frac{1}{3}}, \mu_y^{\frac{1}{3}})$, thus $(x, y)$ cannot appear in any $C_3$ or $C_4$ or $C_5$.

\subsection{$d(x)=3, d(y)=2$}
By Lemma \ref{lemma1}, $\frac{2}{d(x)}+\frac{2}{d(y)}-2=\frac{2}{3}+\frac{2}{2}-2\neq 0$. Thus $(x, y)$ must be either in $C_3$, or $C_4$ or $C_5$. 

\begin{itemize}

\item If $(x, y)$ in $C_4$, see the following graph, taking 
$A(3, 0)=A(3, 4)=\frac{1}{8}, A(2, 4)=\frac{1}{4}$ and other values $A(i, j)=0$, we have $W(\mu_x^{\frac{1}{4}}, \mu_y^{\frac{1}{4}})\leq A(3, 0)\times d(3, 0)+A(3, 4)\times d(3, 4)+ A(2, 4)\times d(2,4)=\frac{1}{8}\times 1 + \frac{1}{8}\times 3 +\frac{1}{4}\times 1 =0.75<1$. 

\begin{center}
\includegraphics[scale=0.4]{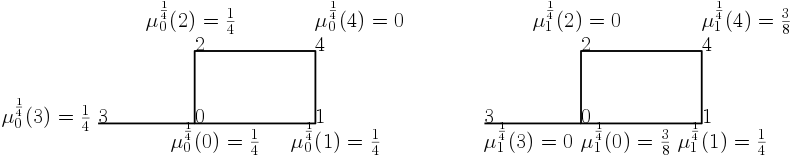} 
\end{center}
Thus $(x, y)$ cannot appear in any $C_3$ or $C_4$. 

%
%

\end{itemize}

\subsection{$d(x)=4, d(y)=2$}
By Lemma \ref{lemma1}, $\frac{2}{d(x)}+\frac{2}{d(y)}-2=\frac{2}{4}+\frac{2}{2}-2 \geq 0$. Thus $(x, y)$ must be either in $C_3$, or $C_4$ or $C_5$.

\begin{itemize}
\item If $(x, y)$ in $C_3$, see the following graph, $W(\mu_x^{\frac{1}{5}}, \mu_y^{\frac{1}{5}})<1$  by taking 
$A(3, 0)=A(4, 2)=\frac{1}{5}$  and other values $A(i, j)=0$.
\begin{center} 
 \includegraphics[scale=0.4]{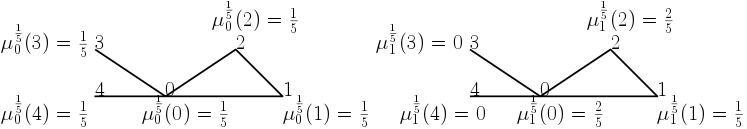} 
\end{center}

\item If $(x, y)$ is in $C_4$. See the following $W(\mu_x^{\frac{1}{5}}, \mu_y^{\frac{1}{5}})=1$. For $W(\mu_x^{\frac{1}{5}}, \mu_y^{\frac{1}{5}})< 1$, we take $A(0, 3)=A(5, 2)=A(5, 4)=\frac{1}{5}$  and other values $A(i, j)=0$. 
For $W(\mu_x^{\frac{1}{5}}, \mu_y^{\frac{1}{5}})\leq 1$, we take $f(0)=2, f(1)=2, f(2)=f(3)=f(4)=1$ and $f(5)=3$. 

\begin{center}
\includegraphics[scale=0.4]{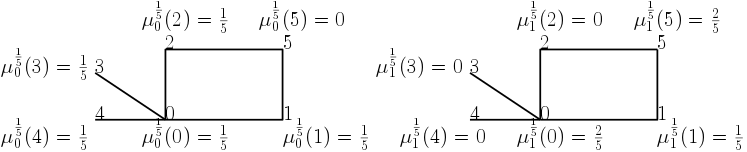} 
\end{center}

Note if $d(3, 5)=2$, $W(\mu_x^{\frac{1}{5}}, \mu_y^{\frac{1}{5}})< 1$ by taking
$A(4, 0)=A(2, 5)=A(3, 5)=\frac{1}{5}, $ and other values $A(i, j)=0$. Similarly for $d(4, 5)$. 
Thus if $(0, 1)$ is in $C_4$ with $d(2, 5)=1$, then $d(3, 5)$ and $d(4, 5)$ must be $3$.

%
%
%
%
%
\end{itemize}


\subsection{$d(x)=3, d(y)=3$}
By Lemma \ref{lemma1}, $\frac{2}{d(x)}+\frac{2}{d(y)}-2=\frac{2}{3}+\frac{2}{3}-2\neq 0$. Thus $(x, y)$ must be either in $C_3$, or $C_4$ or $C_5$.

\begin{itemize}
\item If $(x, y)$ is  in $C_3$, see the following graph, $W(\mu_x^{\frac{1}{4}}, \mu_y^{\frac{1}{4}})<1$  by taking
$A(3, 4)=\frac{1}{4}$ and other values $A(i, j)=0$.  Contradiction. 

\begin{center}
\includegraphics[scale=0.4]{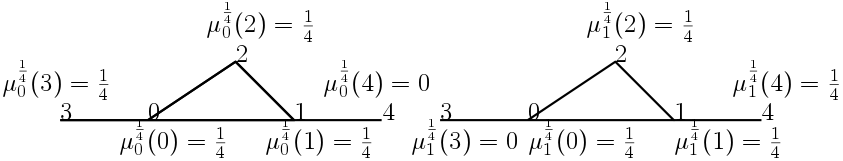} 
\end{center}

\item If $(x, y)$ is in $C_4$, see the following graph, $W(\mu_x^{\frac{1}{4}}, \mu_y^{\frac{1}{4}})=1$. 
For $W(\mu_x^{\frac{1}{4}}, \mu_y^{\frac{1}{4}})\leq 1$, we can take $A(2, 4)=A(3, 5)=\frac{1}{4}$.
For $W(\mu_x^{\frac{1}{4}}, \mu_y^{\frac{1}{4}})\geq 1$, we can take $f(0)=2, f(1)=1, f(2)=f(3)=3, f(4)=2, f(5)=0$. 

\begin{center}
\includegraphics[scale=0.4]{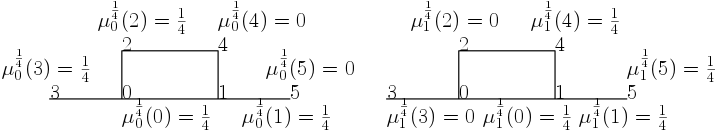} 
\end{center} 
Note if $d(2, 4)=1$, then $d(3, 5)$ must be $3$. 
If further $d(2, 5)=1$ or $2$ then  $d(3, 4)$ must be $3$, and we still have $W(\mu_x^{\frac{1}{4}}, \mu_y^{\frac{1}{4}})=1$.
Similarly, if further $d(3, 4)=1$ or $2$ then $d(2, 5)$ must be $3$, still $W(\mu_x^{\frac{1}{4}}, \mu_y^{\frac{1}{4}})=1$.

%
%
%
%

\end{itemize}

\subsection{$d(x)=3, d(y)=4$}

By Lemma \ref{lemma2}, $\frac{1}{d(x)}+\frac{2}{d(y)}-1=\frac{1}{3}+\frac{2}{4}-1<0$. Thus $(x, y)$ must be either in $C_3$ or $C_4$.

\begin{itemize}
\item If $(x, y)$ is in $C_3$, see the following graph. $W(\mu_x^{\frac{1}{5}}, \mu_y^{\frac{1}{5}})=1$.
For $W(\mu_x^{\frac{1}{5}}, \mu_y^{\frac{1}{5}})\leq 1$, we can take $A(1, 4)=A(1, 5)=\frac{1}{30}$, $A(3, 4)=A(3, 5)=\frac{2}{15}$, $A(2, 4)=A(2, 5)=\frac{1}{30}$  and other values $A(i, j)=0$.
For $W(\mu_x^{\frac{1}{5}}, \mu_y^{\frac{1}{5}})\geq 1$, we can take $f(0)=2, f(1)=1, f(2)=1, f(3)=3, f(4)=f(5)=0$. 

\begin{center}
\includegraphics[scale=0.4]{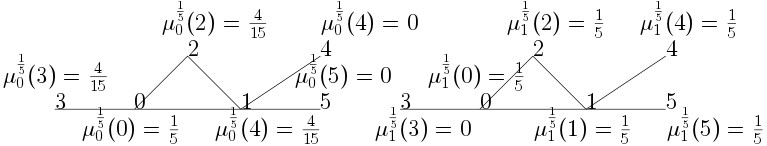} 
\end{center}

If $d(3, 5)=2$, taking the same function $A(i, j)$, we have $W(\mu_x^{\frac{1}{5}}, \mu_y^{\frac{1}{5}})<1$. Thus $d(3, 5), d(3, 4)$ must be $3$. 

\item If $(x, y)$ is in $C_4$, see the following graph, $W(\mu_x^{\frac{1}{5}}, \mu_y^{\frac{1}{5}})=1.2666$ by taking  $f(0)=2, f(1)=1, f(2)= f(3)=3, f(4)=f(5)=f(6)=0$. 

\begin{center}
\includegraphics[scale=0.4]{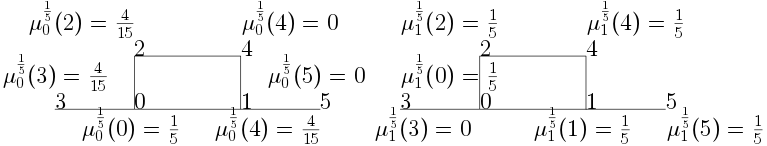} 
\end{center}

We calim $(x, y)$ cannot  share two $C_4$s, otherwise if both $d(2, 4), d(3, 5)$ are $1$, that is, $(0, 1)$ shares two $C_4$s, then $W(\mu_x^{\frac{1}{5}}, \mu_y^{\frac{1}{5}})\leq 0.8666$ by taking $A(2, 4)=A(3, 5)=\frac{1}{5}$  and other values $A(i, j)=0$.  

Thus there are two main cases left: 

\begin{itemize}

\item $(x, y)$ is contained in exactly two $C_4$s. 

\subitem Case a: $(0, 1)$ is in two $C_4$s with  $C_4=0-2-4-1-0$ and $C_4=0-2-5-1-0$.
 
We claim that $d(3, 6)$ cannot be $2$ under this assumption. Otherwise, take  $A(1, 5)=A(2, 5)=\frac{1}{15}, A(2, 4)=\frac{3}{15},  A(3, 5)=\frac{1}{15},  A(3, 6)=\frac{2}{15}$, then $W(\mu_0^{\frac{1}{5}}, \mu_1^{\frac{1}{5}})\leq \frac{1}{15}+\frac{1}{15}+\frac{3}{15}+ d(3, 5)\times \frac{1}{15}+d(3, 6)\times\frac{2}{15}=\frac{9}{15}+\frac{d(3, 5)}{15}\leq \frac{12}{15}<1$, a contradiction. 

We claim  whatever $d(2, 6)$ is ($2$ or $3$), 
$d(3, 4), d(3, 5)$ cannot be both $3$, since we can take 
$f(0)=2, f(1)=1, f(2)=1,  f(3)=3, f(4)=f(5)=f(6)=0$, then $W(\mu_0^{\frac{1}{5}}, \mu_1^{\frac{1}{5}})
=f(0)\times 0+f(1)\times \frac{1}{15}+(f(2)+f(3))\times \frac{4}{15}-(f(4)+f(5)+f(6))\times \frac{3}{15}
= \frac{1}{15}+\frac{16}{15}> 1$.  Thus one of $d(3, 4), d(3, 5)$ must be $2$.

By the symmetry of vertex $4$ and vertex $5$, wlog,  let $d(3, 5)=2$. See the following graph, we have  $W(\mu_1^{\frac{1}{5}}, \mu_1^{\frac{1}{5}})=1$ by taking $A(1, 6)=A(2, 5)=\frac{1}{15}, A(3, 6)=A(3, 5)=\frac{2}{15},  A(2, 4)=\frac{3}{15}, $ and 
$f(0)=2, f(1)=1, f(2)=2,  f(3)=3, f(4)=f(5)=1, f(6)=0$. 

\subitem Case b: $(0, 1)$ is in two $C_4$s with $C_4=0-1-4-2-0$,  $C_4=0-1-4-3-0$.

We analysis the values of $d(2, 5), d(2, 6), d(3, 5), d(3, 6)$. Currently we have  $W(\mu_0^{\frac{1}{5}}, \mu_1^{\frac{1}{5}})\geq 1.26666$. Thus one of  $d(2, 5), d(2, 6), d(3, 5), d(3, 6)$ must be $2$. Wlog,  let $d(2, 5)=2$. We calim that $d(3, 6)=3$, otherwise let $d(2, 6)=d(3, 5)=3$,  take $A(1, 5)=\frac{1}{15}, 
A(2, 4)=\frac{2}{15}, A(3, 4)=\frac{1}{15}, A(3, 6)=\frac{3}{15}$ and others $A(i, j)=0$,  then $W(\mu_0^{\frac{1}{4}}, \mu_1^{\frac{1}{4}})\leq 0.9333$, thus $d(3, 6)$ must be $3$. Similarly, if $d(2, 6)=2$, then $d(3, 5)$ must be $3$. 
Assume $d(2, 6)=3, d(3, 5)=2$, take the 1-Lipschitz function $f(0)=2, f(1)=1, f(2)=f(3)=3, f(4)=2, f(5)=1, f(6)=0$, then $W(\mu_0^{\frac{1}{4}}, \mu_1^{\frac{1}{4}})\geq \frac{16}{15}>1$. Thus $d(2, 6)=2$,  $d(3, 5)=3$.

\item $(x, y)$ is contained in exactly one $C_4$, then at least one of $d(3, 5), d(3, 6)$  is $2$. 

Wlog, let $d(3, 5)=2$, we have $W(\mu_x^{\frac{1}{5}}, \mu_y^{\frac{1}{5}})\geq 1.066666$ by taking  $f(0)=2, f(1)=1, f(2)= f(3)=3, f(4)=2, f(6)=0, f(5)=1$. 
Now there are two cases to consider. 

\begin{enumerate}

\item{If further $d(3, 6)=2$, then $d(2, 5)=d(2, 6)=3$. }
%
For $W(\mu_x^{\frac{1}{5}}, \mu_y^{\frac{1}{5}})\leq 1$, we can take
$A(1, 5)=A(1, 6)=A(2, 5)=A(2, 6)=\frac{1}{30}$, $A(3, 5)=A(3, 6)=\frac{2}{15}$, $A(2, 4)=\frac{3}{15}$  and other values $A(i, j)=0$.
For $W(\mu_x^{\frac{1}{5}}, \mu_y^{\frac{1}{5}})\geq 1$, we can take $f(0)=2, f(1)=1, f(2)=3, f(3)=2, f(4)=2, f(5)=f(6)=0$. From the values of $A(i, j)$, $d(2, 5), d(2, 6)$ must be $3$. 

\item{If further $d(3, 6)=3$, then $d(2, 6)=2$, $d(2, 5)$ could be $2$ or $3$}.

When $d(2, 6)=2$, taking $A(1, 6)=A(2, 6)=\frac{1}{15}$, $A(3, 5)=\frac{2}{15}$, $A(2, 4)=\frac{3}{15}$  and other values $A(i, j)=0$, we get $W(\mu_0^{\frac{1}{5}}, \mu_1^{\frac{1}{5}})\leq \frac{3}{15}+2\times \frac{1}{15} + 2\times \frac{3}{15}+\frac{1}{15}+3\times \frac{1}{15}=1$. Taking $f(0)=2, f(1)=1, f(2)=2, f(3)=3, f(4)=1, f(5)=1, f(6)=0$, we have $W(\mu_0^{\frac{1}{5}}, \mu_1^{\frac{1}{5}})\geq 1$.
From the values of $f(i)$, $d(2, 6)$ must be $2$. 
\end{enumerate}

\end{itemize}

\end{itemize}

\subsection{$d(x)=4, d(y)=4$}
%

By Lemma \ref{lemma2}, $\frac{1}{d(x)}+\frac{2}{d(y)}-1=\frac{1}{4}+\frac{2}{4}-1<0$. Thus $(x, y)$ must be either in $C_3$ or $C_4$.

\begin{itemize}

\item If $(x, y)$ is in $C_3$, see the following graph. $W(\mu_x^{\frac{1}{5}}, \mu_y^{\frac{1}{5}})\geq 1.2$, since  we can take $f(0)=2, f(1)=f(2)=1, f(3)=f(4)=3, f(5)=f(6)=0$. 

\begin{center}
\includegraphics[scale=0.5]{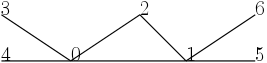}
\end{center}


Note $(x, y)$ cannot share both $C_3$ and $C_4$. For example, if $d(3, 5)=1$, ever $d(4, 6)=3$, we have $W(\mu_0^{\frac{1}{5}}, \mu_1^{\frac{1}{5}})\leq 0.8$ by taking $A(3, 5)=A(4, 6)=\frac{1}{5}. $ It implies that $(x, y)$ cannot share two $C_3$s.
\begin{center}
\includegraphics[scale=0.4]{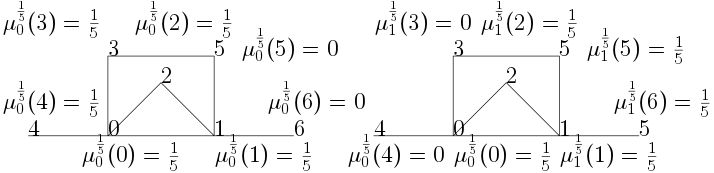}
\end{center}

Then we need to solve $d(3, 5)\times A(3, 5)+d(3, 6)\times A(3, 6)+d(4, 5)\times A(4, 5)+d(4, 6)\times A(4, 6)$, by symmetric of these four vertices, there is only one solution: if $d(3, 5)=2$, then $d(4, 6)$ must be $3$ and  $A(3, 5)=A(4, 6)=\frac{1}{5}$, $d(3, 6), d(4, 5)$ could be $2$ or $3$. 
To see this, let $d(3, 5)=2$, if $d(4, 6)=2$, then $W(\mu_0^{\frac{1}{5}}, \mu_1^{\frac{1}{5}})\leq 0.8$ by taking $A(3, 5)=A(4, 6)=\frac{1}{5}$  and other values $A(i, j)=0$. Thus we have when the edge $(0, 1)$ is in the $C_3:=0-1-2-0$, then if $d(3, 5)=2$, then $d(4, 6)=3$. 
%

If further $d(3, 6)=3$, then $d(4, 5)$ could  be $2$ or $3$, since we can take $f(0)=2, f(1)=f(2)=1=f(5), f(3)=f(4)=3$  and $f(6)=0$ such that  $W(\mu_0^{\frac{1}{5}}, \mu_1^{\frac{1}{5}})\geq 1 $ and take $A(3, 5)=A(4, 6)=\frac{1}{5}$  and other values $A(i, j)=0$ such that $W(\mu_0^{\frac{1}{5}}, \mu_1^{\frac{1}{5}})\leq 1 $. 

\begin{center}
\includegraphics[scale=0.4]{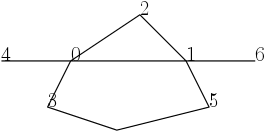}
\end{center}

%

\item When $(x, y)$ is in one $C_4$, see the following graph.  $W(\mu_x^{\frac{1}{5}}, \mu_y^{\frac{1}{5}})=1.4$.
\begin{center}
\includegraphics[scale=0.4]{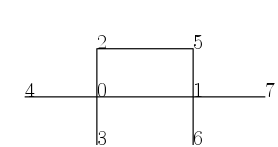}
\end{center} 

Then we analysis the other pair of neighbors $(3, 6), (3, 7), (4, 6), (4, 7)$.

\begin{itemize}
\item{Assume $d(3, 6)=1$, then $d(4, 7)$ must be $3$}. Otherwise $W(\mu_0^{\frac{1}{5}}, \mu_1^{\frac{1}{5}})\leq 0.8$ by taking
$A(2, 5)=A(3, 6)=A(4, 7)=\frac{1}{5}$. Let $d(4, 7)=3$, then $(x, y)$ shares two $C_4$s, we have 
$W(\mu_0^{\frac{1}{5}}, \mu_0^{\frac{1}{5}})\leq 1$ by taking
$A(2, 5)=A(3, 6)=A(4, 7)=\frac{1}{5}$;  
and $W(\mu_x^{\frac{1}{5}}, \mu_y^{\frac{1}{5}})\geq 1$ by taking take $f(0)=2, f(1)=1, f(2)=f(3)=f(4)=3, f(5)=f(6)=2, f(7)=0$.

When $d(2, 6)=1$ (and $d(3, 5)=1$), then $d(4, 5), d(4, 6)$ could be $1$.
When $d(2, 7)=1$ or $d(3, 7)=1$, then $d(4, 5), d(4, 6)$ must be $3$. 
When $d(2, 7)=2$ (and $d(3, 7)=2$ ), $d(4, 5), d(4, 6)$ could be $2$ or $3$.

\item{Assume $d(3, 6)=2$,  then $d(4, 7)$ must be $2$}.  Otherwise we can take $f(0)=2, f(1)=2, f(2)=f(3)=f(4)=3, f(5)=2, f(6)=1, f(7)=0$, then  $W(\mu_0^{\frac{1}{5}}, \mu_1^{\frac{1}{5}})\geq 1.2$. Let $d(3, 6)=2=d(4, 7)=2$, then
$W(\mu_0^{\frac{1}{5}}, \mu_1^{\frac{1}{5}})=1$ by taking 
$f(2)=f(3)=f(4)=3, f(5)=2, f(7)=1, f(6)=1$ and $A(2, 5)=A(3, 6)=A(4, 7)=\frac{1}{5}$.

\begin{center}
\includegraphics[scale=0.4]{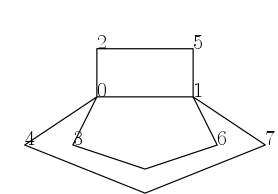} 
\end{center}
\end{itemize}

If further $d(2, 6)=d(2, 7)=1$, $d(3, 5)=d(4, 5)=d(4, 6)=d(3, 7)=2$, then still $W(\mu_0^{\frac{1}{5}}, \mu_1^{\frac{1}{5}})=1$. 

\end{itemize}

\subsection{Conclusion on local structures}
Based on above analysis, we get the following summary. 
\begin{lemma}\label{222324333444}
Suppose $(x, y)$ is an edge in a graph $G$ with  Ricci curvature $k(x, y) = 0$.  Then one of the following statements holds.
\begin{enumerate}

\item[Type 1:]  $d(x)=2, d(y)=2$, then $(x, y)$ is not in any $C_3, C_4$ or $C_5$. 

\item[Type 2:]  $d(x)=3, d(y)=2$, then $(x, y)$ is in exactly one $C_5$. In other words, let 
$x_1, x_2$ be the neighbors of vertex $x$, $y_1$ be the neighbor of $y$, then $\{d(x_1, y_1), d(x_2, y_2)\}=\{2,3\}$.

\item[Type 3:]  $d(x)=4, d(y)=2$, then $(x, y)$ is either in exactly one $C_4$ or is in at least two $C_5$s. In other words,
let 
$x_1, x_2, x_3$ be the neighbors of vertex $x$, $y_1$ be the neighbor of $y$, if $d(x_1, y_1)=1$, then $d(x_2, y_1)=d(x_3, y_1)=3$. 
If $d(x_1, y_1)=2$, then at least one of $d(x_2, y_1), d(x_3, y_1)$ is $2$.

\item[Type 4:]  $d(x)= 3, d(y)=3$,  then $(x, y)$ either is in at least one $C_4$ or shares  two $C_5$s. In other words,
let 
$x_1, x_2$ be the neighbors of vertex $x$, $y_1, y_2$ be the neighbors of $y$. Then there are two main cases:

\begin{itemize}
\item[Case 1:]
If $d(x_1, y_1)=1$, then $d(x_2, y_2)=3$. If further $d(x_1, y_2)=1, 2$, then $d(x_2, y_1)=3$. Similarly, if $d(x_2, y_1)=1, 2$, then $d(x_1, y_2)=3$.
\item[Case 2:] 
If $d(x_1, y_1)=2$, then $d(x_2, y_2)=2$. 

\end{itemize}

\item[Type 5:]  $d(x)=3, d(y)=4$, then $(x, y)$ either shares one $C_3$ or shares one $C_4$ plus one $C_5$. In other words,
 let 
$x_1, x_2$ be the neighbors of vertex $x$, $y_1, y_2, y_3$ be the neighbors of $y$. Then there are three main cases:

\begin{itemize}
\item[Case 1:]\label{Type5case1}
If $x_1= y_1$, then $d(x_2, y_2)=d(x_2, y_3)=3$; 
\item[Case 2:]\label{Type5cas2}
If $d(x_1, y_1)=d(x_1, y_2)=1$, then $d(x_2,  y_3)=3$ and at least one of $d(x_2, y_1), d(x_2, y_2)$ is $3$;

If $d(x_1, y_1)=d(x_2, y_1)=1$, then $\{d(x_1,  y_2), d(x_1, y_3), d(x_2, y_2), d(x_2, y_3)\}$ is in the set $\{2, 2, 3, 3\}$ or $\{3, 3, 2, 2\}$;
\item[Case 3:]\label{Type5cas31}
If $d(x_1, y_1)=1$, then at least one of  $d(x_2, y_2),  d(x_2, y_3)$ is $2$. 
Further, there are two types. 
\begin{itemize}
\item[Type 5a:]
If $d(x_2, y_2)=d(x_2, y_3)=2$, then  $d(x_1, y_2)=d(x_1, y_3)=3$. 
\item[Type 5b:] 
If $d(x_2, y_2)=2, d(x_2, y_3)=3$,  then $d(x_1, y_3)=2$, $d(x_1, y_2)=2$  or $3$.

\end{itemize}
\end{itemize}

\item[Type 6:]  $d(x)=d(y)=4$, then $(x, y)$ either shares one $C_3$ plus one $C_5$ or shares two $C_4$s or shares one $C_4$ plus two $C_5$s. In other words, let 
$x_1, x_2, x_3$ be the neighbors of vertex $x$, $y_1, y_2, y_3$ be the neighbors of $y$. Then there are three main cases:
\begin{itemize}
\item[Type 6a:]
$x_1=y_1$ and $d(x_2, y_2)=2$, then $d(x_3, y_3)$ must be $3$. 
\item[Type 6b:]
$d(x_1, y_1)=1$ and $d(x_2, y_2)=1$, then $d(x_3, y_3)$ must be $3$. 

\item[Type 6c:]
$d(x_1, y_1)=1$,  $d(x_2, y_2)=2$ and $d(x_3, y_3)=2$. 
\end{itemize}
\end{enumerate}
\end{lemma}


\subsection{Further analysis}
Above conclusions can be applied to Ricci-flat graphs without degree restriction. We now give several simple facts and further study these conclusions on Ricci-flat graphs with maximum degree at most $4$. 

\begin{lemma}\label{C3C4C5}
Let G be a Ricci-flat graph that does not contain edge with endpoints degree $\{2, 2\}$, then any edge in $G$ is either contained in a $C_3$ or  $C_4$ or $C_5$. 
\end{lemma}

\begin{lemma}\label{noc3c4}
Let G be a Ricci-flat graph  with maximum degree at most $4$, then no edge in $G$ shares $C_3$ and $C_4$. 
\end{lemma}

\begin{lemma}\label{34noc3}
Let G be any Ricci-flat graph   with maximum degree at most $4$ that contains an edge $(x, y)$ with  $d(x)=3$ and  $d(y)=4$, then $(x, y)$ is not contained in any $C_3$.  Thus Case 1 in ``Type 5" is excluded for our purpose. 
\end{lemma}

\begin{proof}
We suppose $G$ contains and edge $(0, 1)$ with $d(0)=3$, $d(1)=4$ such that $0\sim 2, 3, 1\sim 2, 3, 4$. See the following graph.  By previous analysis $d(3, 5), d(3, 4)$ must be $3$.  Observe that the edge $(0,2)$ is in a $C_3$ with $d(0)=3$,  then $d(2)=4$. By Lemma \ref{noc3c4}, $2\not\sim 4$ or $2\not\sim 5$.  Let $2\sim 6$ and $2\sim 7$.  Since the edge $(0, 2)$ is contained in a $C_3$ sharing vertex $1$, thus $d(3, 6), d(3, 7)$ must be $3$. Then the edge $(0, 3)$ cannot be contained in  any $C_3$ or $C_4$ or $C_5$, a contradiction to Lemma \ref{C3C4C5}. A contradiction. 

\begin{center}
\includegraphics[scale=0.4]{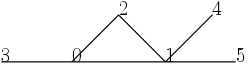} 
\end{center}

\end{proof}


\begin{lemma}\label{34notwoc4}
Let G be a Ricci-flat graph  with maximum degree at most $4$, if there exists an edge $(x, y)$ with  $d(x)=3$ and  $d(y)=4$,  then  $(x, y)$ must be contained in exactly one $C_4$. Thus Case 2 in ``Type 5" is excluded for our purpose. 
\end{lemma}

\begin{proof}
Let $x=0, y=1$ and $0\sim 2, 3, 1\sim 4, 5, 6$, Let $(0, 1)$ is in the $4$-cycle $C_4=0-2-4-1-0$.  There are two different situations for the edge $(0, 1)$ to be contained in two $C_4$s. 

\begin{itemize}
\item $(0, 1)$ is in two $C_4$s with  $C_4=0-2-4-1-0$ and $C_4=0-2-5-1-0$. By previous analysis, we have $d(3, 6)=3$ and one of $d(3, 4), d(3, 5)$ is $2$. Wlog, let $d(3, 5)=2$. 
 
\begin{center}
\includegraphics[scale=0.4]{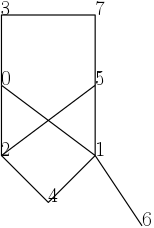} 
\hfil
\includegraphics[scale=0.4]{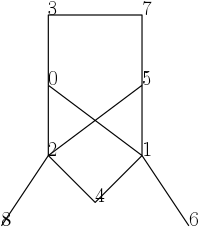} 
\end{center}

Let us consider the edge $(0, 2)$, note $d(2)\neq 3$, since the edge $(0, 2)$ does not satisfy ``Type 4".  Thus $d(2)$ must be $4$. B Lemma \ref{noc3c4}, vertex $2$ is not adjacent to vertex $3$, $6$ or $7$,  we need a new vertex $8\sim 2$.  We have $d(3, 8)=3$ with the same reason as $d(3, 6)=3$. 
From above local structure,  the edge $(0,3)$ is not ``Type 2", then $d(3)\neq 2$. Since edge $(0, 3)$ cannot be in any $C_4$, then $d(3)\neq 4$. Then it must be $d(3)=3$, and edge $(0, 3)$ must share two separated $C_5$s which would lead to $d(3, 6)=2$ or $d(3, 8)=2$, a contradiction. 

\item $(0, 1)$ is in two $C_4$s with $C_4=0-1-4-2-0$,  $C_4=0-1-4-3-0$.  By previous analysis, wlog, let $d(2, 5)=d(2, 6)=2, d(3, 5)=d(3, 6)=3$.
%
%
%
Let vertex $7$ connect $2$ and $5$.  
\begin{center}
\includegraphics[scale=0.4]{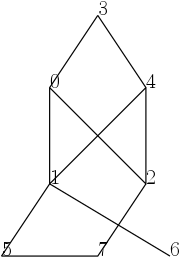}
\end{center}

Now we have $d(2), d(4)\geq 3$, there are several different cases to consider.

\begin{description}

\item[If $d(2)=3, d(4)=3$, ] then we have $W(\mu_2^{\frac{1}{5}}, \mu_4^{\frac{1}{5}})\leq 0.75$ as the edge $(2, 4)$ does not satisfies ``Type 4".  

\item[If $d(2)=3, d(4)=4$,] consider the edge $(2, 4)$, it is contained in two $C_4:= 2-4-1-0-2, 2-4-3-0-2$, which  can be compared with the first situation of edge $(0, 1)$,  this case is excluded. 

\item[If $d(2)=4, d(4)=3$, ] let $2\sim 8$.  For edge $(0, 2)$ , it is contained in two $C_4:=0-2-4-3-0, 0-2-4-1-0$ and $d(1, 7)=2$,  then $d(1, 8)$ must be $2$, $d(3, 7), d(3, 8)=3$.  Since $(0, 3)$ is contained in two $C_4$s, then $d(3)\neq 2$. If $d(3)=3$ with the third neighbor $9$, then $d(9, 1)=d(9, 2)=3$ which would result that  the edge $(3, 9)$ cannot be in any $C_3, C_4$ or $C_5$, a contradiction. Let $d(3)=4$ with $3\sim 9$ and $3\sim 10$, compare the edge $(0, 3)$ with the second situation of edge $(0, 1)$, then either $d(1, 9)=d(2, 9)=2, d(1, 10)=d(2, 10)=3$ or $d(1, 9)=d(2, 9)=3, d(1, 10)=d(2, 10)=2$, wlog, we use the latter case, then the edge $(3, 9)$ cannot be in any $C_3, C_4$ or $C_5$, a contradiction. 

\item[If $d(2)=4, d(4)=4$, ] let $2\sim 8$, $4\sim 9$.  Consider the edge $(0, 2)$, since $d(1, 7)=2$, then $d(1, 8)=2$, $d(3, 7)=d(3, 8)=3$. Note $d(3, 5)=d(3, 6)=3$, then  any neighbor of $3$ must have distance $3$ to both vertices $1$ and $2$. Then $d(3)\neq 3$ by considering edge $(3, 0)$. Let $d(3)=4$ with $3\sim 10, 3\sim 11$. Then for the edge $(0, 3)$, we have both $d(1, 10)=d(1, 11)=3$ and  $d(2, 10)=d(2, 11)=3$, 
 a contradiction to ``Type 5". 
\end{description}

\end{itemize}
We exclude the two situations, thus $(0, 1)$ must be contained exactly one $C_4$ in a Ricci-flat graph. The result follows. 
\end{proof}

\begin{lemma}\label{33noc4}
Let G be a Ricci-flat graph  with maximum degree at most $4$, if there exists an edge $(x, y)$ with endpoint degree $(d(x), d(y))=(3, 3)$, then $(x, y)$ is not contained in any $C_4$.  Thus Case 1 in ``Type 4" is excluded for our purpose. 
\end{lemma}

\begin{proof}
Let $x=0, y=1$, $0\sim 2, 3$, $1\sim 4, 5$. By contradiction, assume $(0, 1)$ is contained in a $C_4:=0-1-4-2-0$. 
\begin{center}
\includegraphics[scale=0.4]{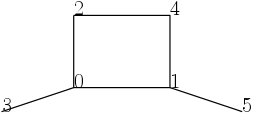} 
\end{center}
There are three different cases. 
\begin{itemize}
\item{Case 1: $d(2, 4)=d(2, 5)=1$}, then $d(3, 4), d(3, 5)$ must be $3$ by ``Type 4".  If $d(2)=3$, then the edge $(0, 3)$ cannot be in any $C_3$ or $C_4$ or $C_5$ through vertex $1$ or $2$, contradicting to Lemma \ref{C3C4C5}. Then $d(2)=4$. For edge $(0, 2)$, we have $(d(0), d(2))=(3, 4)$ and it is the first situation of Lemma \ref{34notwoc4}, this case should be excluded. 

\item{Case 2: $d(2, 4)=1, d(2, 5)=2$},   then $d(3, 4), d(3, 5)$ must be $3$.  Let $2\sim 6\sim 5$. If $d(2)=3$, since $d(1, 4)=1$, then $d(3, 6)$ must be $3$.  Then the edge $(0, 3)$ cannot be in any $C_3$ or $C_4$ or $C_5$ through vertex $1$ or $2$, contradicting to Lemma \ref{C3C4C5}.  Let $d(2)=4$ with $2\sim 7$. Consider the edge $(1, 5)$, note since $d(0, 6)=d(4, 6)=2$, then $d(5)\neq 2$. Since $(1, 5)$ cannot be in any $C_4$, then $d(5)\neq 4$. Thus $d(5)=3$ and $(1, 5)$ shares two $C_5$s. However, since $d(3, 5)=3$, this case cannot happen either. 

\item{Case 3: $d(2, 4)=1, d(2, 5)=d(3, 4)= d(3, 5)=3$.} By the same reason as case 1, $d(2)$ must be $4$. 
Similarly, $d(4)=4$. 
Still the edge $(0, 3)$ cannot be in any $C_3$ or $C_4$, then $d(3)\neq 4$.  Assume $d(3)=3$ with new neighbors $6, 7$, note $(0, 3)$ must share two $C_5$s,  then we need $d(1, 6)=2$ or $d(1, 7)=2$, both would lead to $d(3, 5)=2$, a contradiction. Thus $d(3)=2$ with $3\sim 6$. Similarly, $d(5)=2$ with $5\sim 7$.  Let the edge $(0, 3)$ in a $C_5=0-3-6-8-2-0$, where $8$ is new neighbor of vertex $2$, the edge $(1, 5)$ in a $C_5=1-5-7-9-4-1$, where $9$ is new neighbor of vertex $4$, 
let the fourth neighbor of $2$ be $10$,  the fourth neighbor of $4$ be $11$. Since $d(2, 5)=d(3, 4)=3$, then $d(0, 11)=d(1, 10)=3$, for the edge $(0, 2)$, we need $d(3, 10)=2$, thus $6\sim 10$, similarly, for the edge $(1, 4)$, let $7\sim 11$. 

\begin{center}
\includegraphics[scale=0.3]{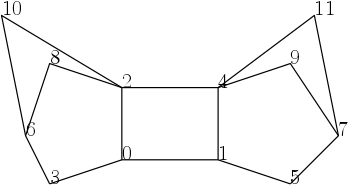} 
\end{center}

Consider the edge $(3, 6)$, since $d(0, 8)=d(0, 10)=2$, then $d(6)\neq 3$. Thus $d(6)=4$, similarly,  $d(7)=4$. 

Consider the edge $(2, 8)$, since it is in a $C_4$ and $C_5$, then $d(8)\neq 2$. By symmetric, $d(9), d(10), d(11)\neq 2$. 
For the edge $(2, 4)$, wlog, either we need $d(8, 9)=d(10, 11)=2$, either $d(8, 9)=1, d(10, 11)=3$. 
For the latter case, since edge $(2, 8)$ is in two separate $C_4$s, then $d(8)$ must be $4$.  Similarly,  $d(9)=4$. Let $12, 13$ be fourth neighbors of $8$ and $9$ respectively. Now let us consider the edge $6, 8$, note $d(6)=d(8)=4$ and $d(10, 2)=1$, both $d(3, 9), d(3, 12)\neq 1$, then one of them must be $2$.  If $d(3, 12)=2$, then  $0\sim 12$ which cannot happen or $6\sim 12$ which will contradict to Lemma \ref{noc3c4};  similarly it cannot be $d(3, 9)=2$, a contradiction.

Thus for the edge $(2, 4)$, we must have that $d(8, 9)=d(10, 11)=2$.  For $d(8, 9)=2$,  assume $8\sim 7$, then the edge $(7, 8)$ must be in a $C_4$, then either $9\sim 6$  or $11\sim 6$, wlog, let $9\sim 6$. Then $d(8)\neq 3$ since both $d(2, 5), d(2, 11)$ are $3$.  Let $d(8)=4$ with the fourth neighbor $12$, then we need $d(12, 11)=1$ for the edge $(8, 7)$.  Consider the edge $(2, 8)$, since both $d(0, 7), d(0, 12)$ are $3$, we need $d(4, 12)=1$ which cannot happen. A contradiction. 

For $d(8, 9)=2$, let $12$ be the common neighbor of $8$ and $9$.   We claim $d(8)\neq 4$.  Otherwise, let $13$ be the fourth neighbor of $8$,  for the edge $(2, 8)$, since $d(6, 10)=1$ and $d(12, 4)=2$, then $d(13, 0)$ must be $2$, however, this cannot happen. Thus $d(8)=3$ with neighbors $2, 6, 12$. Similarly, $d(9)=d(10)=d(11)=3$. 
\begin{center}
\includegraphics[scale=0.3]{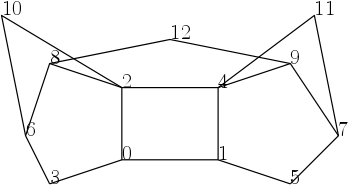} 
\end{center}

Note we need another new vertex to connect $10, 11$ such that $d(10, 11)=2$. Note it cannot be $12$ since then  $(8, 2)$ is contained in two $C_4$s which contradicts to Lemma \ref{34notwoc4}. Let this vertex be $13$. Observe that the edge $(8, 12)$ cannot be in any $C_4$, thus $d(12)\neq 4, d(13)\neq 4$.  When $d(12)=d(13)=2$, we only need to consider edge $(3, 6)$ and $(3, 7)$. 
Let $6\sim 14$, assume $14\sim 7$.  Consider the edge $(6, 8)$, we need $d(12, 14)=2$, let $12\sim 15\sim 14$, then $d(4)\geq 3$. Similarly, we need $d(13, 14)=2$ for the edge $(6, 10)$, thus the edge $(6, 14)$ cannot be in any $C_4$ which implies $d(14)\neq 3, 4$. A contradiction. 
Let vertex $14, 15$  be neighbors of $6, 7$ respectively. And $d(12)=d(13)=3$. 
\begin{center}
\includegraphics[scale=0.3]{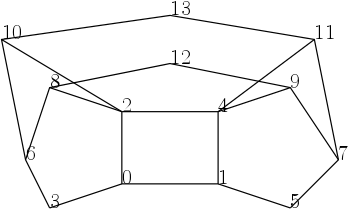} 
\end{center}

Now consider the edge $(6, 10)$  since $d(14, 13)$ must be $2$, we need new vertex $16$ to connect them. Similarly, consider edge $(7, 11)$, we need $d(13, 15)=2$. 
Then $15$ must be adjacent to vertex $16$ in order to keep $d(13)=3$.  Since $(13, 16)$ cannot be in a $C_4$, then $d(16)=3$.

For the edge $(6, 8)$, we need  $d(12, 14)=2$, thus need new vertex $17$ to connect them. Similarly, the edge $(7, 9)$, we need  $d(12, 15)=2$, since $(8, 12)$ cannot be in a $C_4$, then $d(12)\neq 4$, thus $15$ must be adjacent to $17$ in order to have $d(12, 15)=2$.
 \begin{center}
\includegraphics[scale=0.3]{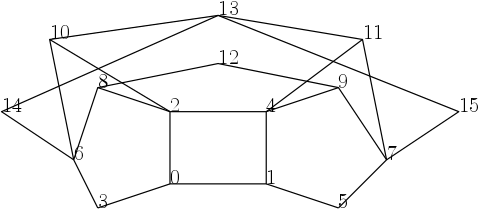} 
\end{center}
 
Observe that the edge $(6, 14), (7, 15)$ cannot be in any $C_3$ or $C_4$, thus $d(14), d(15)$ cannot be $3$ or $4$, a contradiction. 

\end{itemize}

\end{proof} 

\begin{lemma}\label{twoc5property}
Let G be a Ricci-flat graph  with maximum degree at most $4$, if there exists an edge $(x, y)$ with $d(x)= d(y)=3$, and $(x, y)$ shares two $C_5$s. Then $d(z)=d(w)$ if both $z, w$ are adjacent to $x$ or $y$. 
\end{lemma}
\begin{proof}
Look at the following graph with $x=0, y=1$. We will show that $d(2)=d(3)$, then $d(4)=d(5)$ will follow immediately. 
\begin{center}
\includegraphics[scale=0.4]{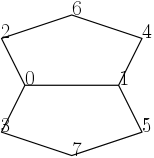} 
\end{center}

Suppose $d(2)=2$, then $d(3, 6)=3$ for the edge $(0, 2)$, thus the edge $(0, 3)$ cannot be in any $C_4$ or $C_5$ which lead to $d(3)=2$. 
Suppose $d(2)=3$ with the third neighbor $8$. Since the edge $(0, 2)$ must share two $C_5$s, then  the edge $(0, 3)$ cannot be in any $C_4$, thus $d(3)=3$. 
Suppose $d(2)=4$, then the edge $(0, 2)$ must be in a $C_4$ which can only pass through vertex $3$, thus the edge $(0, 3)$ is also in a $C_4$, then $d(3)$ must be $4$. 
\end{proof}

\begin{theorem}\label{23determined}\label{Them:halfdode}
Let $G$ be a Ricci-flat graph   with maximum degree at most $4$. If $G$ contains an edge $(x, y)$ with $d(x)=3, d(y)=2$. Then $(x, y)$ shares two $C_5$s, and $G$ is the Half-dodecaheral graph. 
\end{theorem}

\remark{Note Half-dodecaheral graph is proved to be Ricci-flat in \cite{LLY}, however we cannot use their results since we make no assumption on the girth of $G$. }
\begin{proof}
Look at the following graph with $x=0, y=1$.
\begin{center}
\includegraphics[scale=0.4]{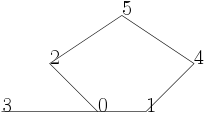}
 \end{center}
 
First note that $d(3, 4)$ must be $3$ by Lemma \ref{222324333444}.  If $d(2)=2$, we must have $d(3, 5)=3$ for the edge $(0, 2)$,  then $(0, 3)$ cannot be in any $C_3, C_4$ or $C_5$, a contradiction.
If $d(2)=4$, with new neighbors $6, 7$, then the edge $(0, 2)$ must be in a $C_4$.   By Lemma \ref{33noc4}, we get $d(3)=4$.  Let $C_4 :=2-0-3-6-2$, and  $3\sim 8, 9$. 
Consider the edge $(0, 3)$, we need one of $d(1, 8), d(1, 9)$ to be $2$, which would contradict to  $d(3, 4)=3$.

Thus $d(2)=3$, let $2\sim 6$. For the edge $(0, 2)$, it must be in two $C_5$s, and $d(3)=d(1)=2$, $d(5)=d(6)$. There are two cases:
1) $d(1, 6)=d(3, 5)=2$ and $d(3, 6)=3$, then $6\sim 4$ and need new vertex $7$ such that $3\sim 7\sim 5$. Under this situation, the edge $(2, 5)$ is contained in the $C_4:=2-5-4-6-2$, then $d(5)=4$, let $5\sim 8$. By Item 5 of Lemma \ref{222324333444}, since $d(0, 8)=3$,  we need $d(6, 8)=2$, let $6\sim 9\sim 8$. Note $d(6)=d(5)=4$ and the edge $(2, 6)$ cannot be in two $C_4$s by Lemma\ref{34notwoc4},   thus $6\not\sim 7$. Let $6\sim 10$, we need $d(0, 9)=2$ or $d(0, 10)=2$, both cannot happen, a contradiction.  
 Thus for the edge $(0, 2)$ sharing two $C_5$s.,  it must be the second case 2)$d(1, 5)=d(3, 6)=2$ with new vertex $7$ achieving this and $d(3)=2$. 
\begin{center}
\includegraphics[scale=0.4]{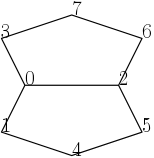}
 \end{center}

Let us focus on vertex $4$. Note $d(4)\neq 2$ and $d(1, 7)=d(3, 4)=3$, then $4\not\sim 7$.
Assume $d(4)=4$, and  assume $z, w$ are third and fourth neighbors of $4$, then one of $d(0, z), d(0, w)$ must be $2$, which which implies $4$ is adjacent to vertex $6$. 
Then the edge $(2, 5)$ is in a $C_4:=4-6-2-5-4$, implies that $d(5)=4$. Note we need new vertices as third and fourth neighbors of $5$, let them be $9, 10$. However both $d(0, 9), d(0, 10)$ cannot be $2$, a contradiction. 

Thus $d(4)=3$. Similarly, $d(7)=3$. Let the third neighbor of $4$ be $z$, by comparing edge $(0, 1)$ and edge $(4, 1)$, we obtained that $d(z)=2$ just like $d(3)=2$, thus $z\neq 6$. Let $4\sim 8$. Actually there is an isomorphism $f$ between the neighborhood $N(0)\cup N(1)$ and the neighborhood $N(1)\cup N(4)$, then $f(0)=4, f(1)=1, f(0)=4, f(2)=5, f(5)=2, f(8)=3$. Thus $d(5)=d(2)=3, d(8)=d(3)=2$. 
Assume $8\sim 7$. Then $f(7)=7$, then $d(7, 5)=d(f(7), f(5))=d(7, 2)=2$, note $5\not\sim 3$ and $5\not\sim 6$, thus we need new vertex to connect vertices $5, 7$, then $d(7)=4$, a contradiction. 
\begin{center}
\includegraphics[scale=0.4]{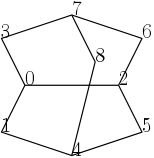}
 \end{center}
 
For $d(7)=3$, let $7\sim 9$. Since the distance from $3$ to any neighbor of $9$ cannot be $2$, then $d(9)$ cannot be $3$.  If $d(9)=4$, then the edge $(7, 9)$ must be in $C_4$ which can only through vertex $6$,  then by considering the edge $(7, 9)$, one of neighbors of $9$ must has distance  $2$ to vertex $3$, which cannot happen. 
 \begin{center}
\includegraphics[scale=0.4]{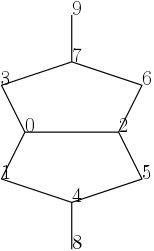}
 \end{center}

Thus $d(9)=2$, similarly, $d((8)=2$.  The edge $(7, 9)$ must be in a $C_5$ through $6$.  If the $C_5:=9-7-6-2-5-9$, then consider the edge $(1, 4)$, we need $d(8, 9)=2$, then either $8\sim 7$ or $9\sim 4$, both cannot happen. 
Thus the $C_5$ cannot pass through vertex $2$, we need new neighbor for $9$ and $6$, let them be $10, 11$ and $10\sim 11$ respectively. Note $d(6)=d(5)=3$.  Then for the edge $(2, 6)$, it must be in two separate $C_5$s, let the other one be $C_5:6-2-5-12-11-6$ where $12$ is a neighbor of $5$.  For the edge $(4, 5)$, we need the other $C_5:5-4-8-14-12-5$ where $14$ is a neighbor of $8$, then $14\sim 12$. 
Now for the edge $(8, 14)$,  any other new neighbor of $14$ cannot have distance $1$ or $2$ to $4$, thus $d(14)$ must be $3$. Similarly $d(10)=3$. Since there is no way for the edge $(6, 11)$ to be in any $C_4$, thus $d(11)=d(12)=3$. Then for the edge $(11, 12)$, we need $d(10, 14)=2$, let $10\sim 13\sim 14$. And $d(13)\neq 3, 4$.  
\begin{center}
\includegraphics[scale=0.4]{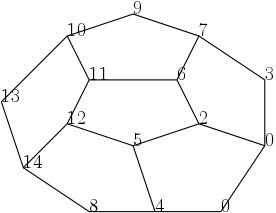}
 \end{center}
\end{proof}

\begin{lemma}\label{34withonec4property}
Let G be a Ricci-flat graph  with maximum degree at most $4$, if there exists an edge $(x, y)$ with $d(x)=3,  d(y)=4$, and $(x, y)$ is contained in exactly one $C_4:=x-x_1-y_1-y-x$ and let $x_2, y_2, y_3$ be the neighbors of $x$ and $y$ that are not on the $C_4$,   then $d(x_2)=3$, $d(x_1)=4$. 
%
\end{lemma}

\begin{proof}
Let $x=0, y=1$ and $0\sim 1, 2, 3$, $1\sim 4, 5, 6$ and the edge $(0, 1)$ be contained in exactly one $C_4=0-1-4-2-0$.  See the following graph, where $x_1=2, x_2=3$. 
\begin{center}
\includegraphics[scale=0.4]{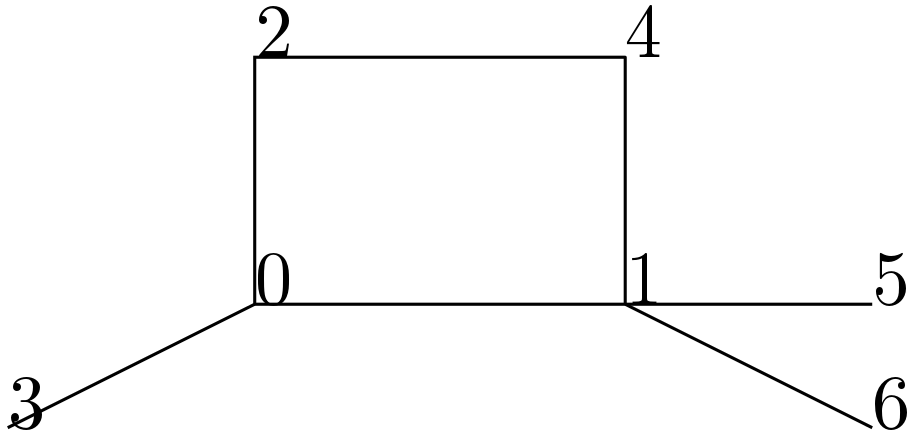}
\end{center}
By previous analysis, we know $d(3, 5)=2$, let $3\sim 7\sim 5$. Since the edge $(0, 2)$ is in a $C_4$, then $d(2)$ must be $4$ by Lemma \ref{33noc4}, thus $d(x_1)=4$. 
Consider the vertex $3$, since $d(3, 5), d(3, 4), d(3, 6)\neq 1$, and the edge $(0, 2)$ is not in two $C_4$s, then the edge $(0, 3)$ cannot be in any $C_4$, thus $d(3)\neq 4$.  Assume $d(3)=2$, consider the edge $(0, 3)$, we must have $d(2, 7)=3$. Let $8, 9$ the third and fourth neighbors of $2$. Thus for the edge $(0, 2)$, one of $d(3, 8), d(3, 9)$ must be $2$, this contradicts to the fact that $d(2, 7)=3$. Thus $d(3)=3$, then $d(x_2)=3$.

%
%
%
%
%
\end{proof}

\section{Ricci-flat graphs containing vertex with degree $3$}
We have determined the Ricci-flat graphs in class $\mathcal{G}$ that contains an edge with endpoint degree $\{2, 3\}$ in Theorem \ref{Them:halfdode}. In this section, we continue with endpoint degree $(3, 3)$ and $(3, 4)$. To determine the former case (see Theorem \ref{33twoC5}), we will need the Theorem \ref{34withonec41}. 
So we start with the graphs that contain an edge with endpoint degrees  $\{3, 4\}$. There are two cases ``Type 5a" and ``Type 5b". Case ``Type 5b" will be excluded for contracting a Ricci flat graph in class $\mathcal{G}$. See Lemma \ref{34withonec42}. 
For ``Type 5a", we have the following four Ricci-flat graphs. 
\begin{theorem}\label{34withonec41}
Let G be a Ricci-flat graph  with maximum degree at most $4$, let $(x, y)$ be an edge in G with $d(x)=3,  d(y)=4$ that satisfies ``Type 5a".  Then $G$ is isomorphic to one of the following  graphs.
 \begin{figure}[H]
\centering
\includegraphics[scale=0.3]{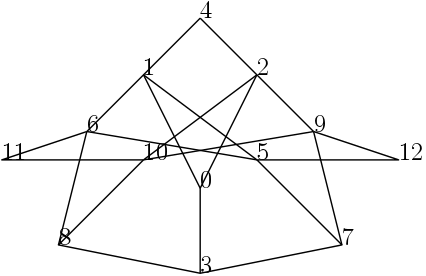} 
\hfil
\includegraphics[scale=0.35]{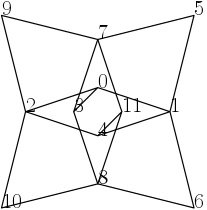} 
\caption{Ricci-flat graphs $G_1, G_2$}
\label{34Ricci-flat1}
\end{figure}

 \begin{figure}[H]
\centering
\includegraphics[scale=0.23]{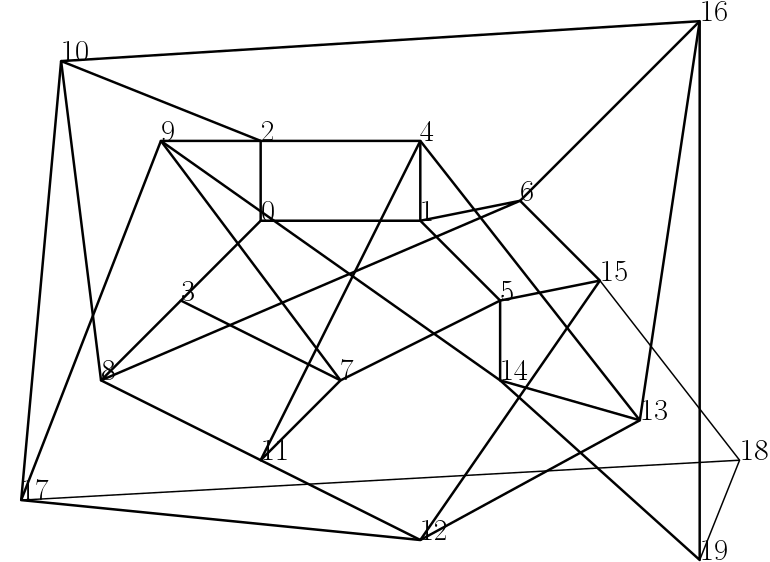}
\hfil
\includegraphics[scale=0.23]{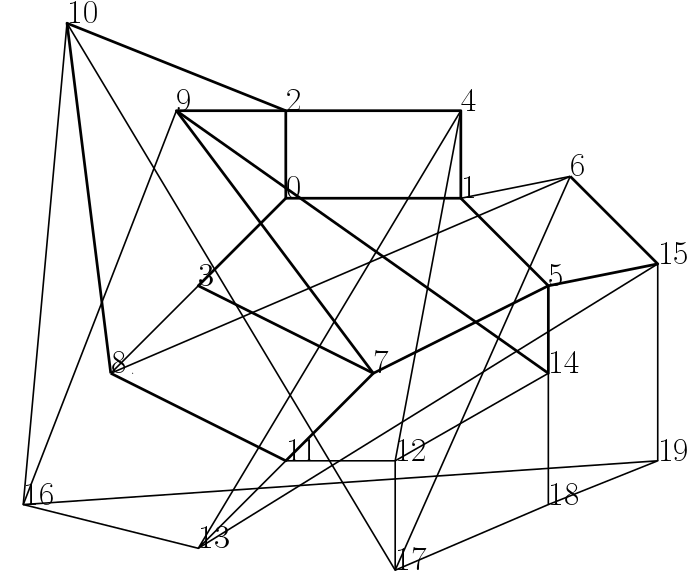} 
\caption{Ricci-flat graphs $G_3, G_4$}
\label{34Ricci-flat2}
\end{figure}

\end{theorem}

\begin{proof}
Let $x=0, y=1$ and $(0, 1)$ be contained in exactly one $C_4:0-2-4-1-10$, with $d(2, 5)=d(2, 6)=3$. By Lemma \ref{34withonec4property},  $d(3)=3$ and $d(2)=4$.  
Let $3\sim 7, 8$,  $7\sim 5$, and $2\sim 9,10$. Then we have the following subgraph of $G$:  
\begin{center}
\includegraphics[scale=0.4]{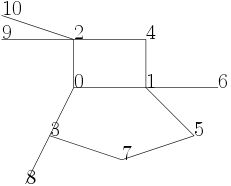} 
\end{center}

We focus on the degree of vertex $7$.  There are two cases:
\begin{enumerate}
\item Case 1a: $d(7)=3$. Since $d(2, 5), d(2, 6)=3$, then $d(1,9)=d(1, 10)=3$. Then for the edge $(0, 2)$,  we need $d(3, 9)=d(3, 10)=2$. Then vertex $7$ must be adjacent to one of $6, 9, 10$, otherwise we need $8$ to be adjacent to all three then $d(8)=4\neq d(7)$, a contradiction to Lemma \ref{twoc5property}.  If $7\sim 6$, then $8\sim 9$ and $8\sim 10$. Consider the edge $(3, 7)$,  it should share two $C_5$s, thus  one of $d(8, 5), d(8, 6)$ must  be $2$, wlog let $d(8, 5)=2$. Note $d(8)=d(7)=3$, then it must be either $5\sim 9$ or $5\sim 10$, which would contradict to the fact that $d(2, 5)= d(2, 6)= 3$. A contradiction. 

Now let $7\sim 9$. Since we require $d(3, 6)=2$, then $8\sim 6$. Still for the edge $(0, 2)$, we need $d(3, 10)=2$ which implies that $8\sim 10$. For the edge $(3, 7)$ to share two $C_5$s, there are two cases, either  $d(8, 9)=2$ or $d(8, 5)=2$, by symmetry of the current graph, wlog, let  $d(8, 9)=2$ which implies that $9\sim 10$, a $C_3$ appears,  thus $d(9)=d(10)=d(2)=4$. For the edge $(8, 10)$, it must be in a $C_4$ which can only  passes through vertex $6$, thus we need a new vertex $11$ as common of $6, 10$.  Similar analysis for the edge $(7, 9)$,  we need a new vertex $12$ as common of  $5, 9$.

Consider $d(4)$, if $d(4)=3$, then we need new vertex as its neighbor say it is $13$. For the edge $(2, 4)$, since$d(1,9)=d(1, 10)=3$,  we need  $d(13, 9)=d(13, 10)=2$, then it must be $13\sim 11$ and $13\sim 12$, however this is not good for the edge $(2, 9)$ and  $(2, 10)$ as for them we must have $d(4, 11)=3, d(4, 12)=3$. If $d(4)=4$, we need two new vertices $13, 14$ as neighbors of $4$. However, $d(10, 13), d(10, 14), d(9, 13), d(9, 14)$ must be $3$ which is not good for the edge $(2, 4)$. 
Thus $d(4)=2$. Similarly, $d(11)=d(12)=2$. 

 Note $d(6), d(5)\neq 3$ by Lemma \ref{33noc4}.  Observe that the edge $(1, 6)$ cannot be in any $C_4$, then it must be in the $C_5:=1-6-5-1$. Then $6\sim 5$. The resulting graph is $G_1$.  

\item Case 1b: $d(7)=4$. There are several cases with respect to the neighbors of vertex $7$.
\begin{itemize}
\item{Case 1b1: $7\sim 6$, $7\sim 9$}. Then $8\sim 10$ for the edge $(0, 2)$. 
 Since the edge $(3, 7)$ must be in a $C_4$, there are two distinct possible cases: $C_4:=7-3-8-5-7$ or $C_4:=7-3-8-9-7$. 
 
We will reject the former case. Otherwise,  we need $d(8, 9)=d(8, 6)=3$ for the edge $(3, 7)$ as $d(0, 5)=d(0, 6)=2$. 
Assume $8\sim 4$. Then both the edges $(2, 4), (5, 8)$ are in two $C_4$s, $d(4)=d(5)=4$. Let $4\sim 11$, a new vertex, then $5\not\sim 11$ for the edge $(1, 5)$. Let $5\sim 12$. However, we need $d(6, 12)=3$ for the edge $(1, 5)$,  $d(10, 12)=3$ for the edge $(8, 5)$,  $d(9, 12)=3$ for the edge $(7, 5)$, then the edge $(5, 12)$ cannot be in any $C_3, C_4$, then it must be in a  $C_5$ which needs $11\sim 12$, and $d(12)=2$.  Now we need $d(6, 11)=3$ for the edge $(1, 4)$, $d(9, 11)=3$ for the edge $(2, 4)$,  $d(10, 11)=3$ for the edge $(4, 8)$, then the edge $(4, 11)$ cannot be in any  $C_3$ and $C_4$. Thus $d(11)=2$, a contradiction for the edge $(11, 12)$. 
Thus vertex $8$ is adjacent to a new vertex. Let  $8\sim 11$.  For the edge $(3, 8)$, since $d(0, 11)$ cannot be $2$, then we need $d(7, 11)=2$ which would be true if $11\sim 6$ or $11\sim 9$, while both cannot happen since otherwise  $d(8, 9)=d(8, 6)=2$. Thus $8\not\sim 5$, similarly, $8\not\sim 6$. 
 
Thus $(3, 7)$ is in the $C_4:=7-3-8-9-7$ which implies $8\sim 9$. Note when $8\sim 9, 10$, by similar arguments as above, $7$ should not be adjacent to vertex $9$. A contradiction. 
As a result,  the edge $(3, 7)$ cannot be in $C_4$, a contradiction.

\item{Case 1b2: Assume $7\sim 6$, $7\sim 11$}.  Note for the edge $(0, 2)$,  we need $8\sim 9, 10$. Since  $(3, 7)$ must be in a $C_4$, and this $C_4$  can only be $7-3-8-11-7$. Thus $8\sim 11$. For edges $(3, 7), (3, 8)$, we need $d(8, 5)=d(8, 6)=3$ since $d(0, 5)=d(0, 6)=2$ and $d(7, 9)=d(7, 10)=3$ since $d(0, 5)=d(0, 10)=2$.  
 \begin{center}
\includegraphics[scale=0.3]{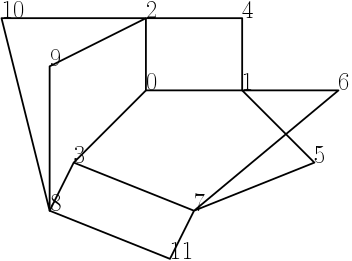} 
\end{center}
Now we consider $d(4)$. Observe that $d(5)\neq 2$, and any new neighbor of $5$ has distance $3$ to vertex $0$, thus at least one of these new neighbors of $5$ should have distance $1, 2$ to vertex $4$ by considering the edge $(1, 5)$. Thus $d(4)\neq 2$.  Assume $d(4)=3$, and assume $4\sim 11$. Then $d(11)$ cannot be just $3$. Let $11\sim 12$. Then $5\sim 12$ for the edge $(1, 5)$. However, this is not good for the edge  $(7, 11)$.  Thus $4\not\sim 11$. Let $12$ be the third neighbor of $4$, then $d(12, 5)=d(12, 6)=2$ for the edge $(1, 5)$ since $d(2, 5)=d(2, 6)=3$. Let $5\sim 13\sim 12$. Note $d(5)=3$, since  both third and fourth neighbors of $5$ have distance $3$ to $0$ and distance at least $2$ to vertex $4$ which is not good for the edge $(1, 5)$ when $d(5)=4$.  
 Since $(1, 5)$ must be in exactly one $C_4$, then $6\not\sim 13$. For $d(12, 6)=2$, let $6\sim 14\sim 12$. Similarly, $d(6)=3$. 
  Note when $d(4)=3$, $d(12)=3$ by similar analysis in Lemma \ref{34withonec4property}. Similarly $d(13)=d(14)=3$. 
  Now we need $d(12, 9)=d(12, 10)=2$ for the edge $(2, 4)$ since $d(1, 9)=d(1, 10)=3$. Then wlog, let $9\sim 13, 10\sim 14$. Now we cannot add new neighbors to vertices $12,13,14$, thus cannot guarantee the edge $(5, 7)$. A contradiction. Thus $d(4)\neq 3$. 

Let $d(4)=4$, similarly $d(11)=4$.  Still assume $4\sim 11$, then the edge $(4, 11)$ have to be in a $C_4$  that passes through new neighbors. Let $4\sim 12\sim 13\sim 11$. Consider the edge $(1, 4)$, since $d(5, 11)=d(6, 11)=2$, then at least one of $d(5, 12), d(6, 12)$ is $2$. Similarly we need at least one of $d(5, 13), d(6, 13)$ is $2$ for the edge $(7, 11)$. Wlog, let $d(5, 12)=2$ with $5\sim 14\sim 12$. Note $d(5)=3$, since  both third and fourth neighbors of $5$ have distance $3$ to $0$ and distance at least $2$ to vertex $4$ which is not good for the edge $(1, 5)$ when $d(5)=4$.  Thus we need $d(14, 11)=2$ for the edge $(5, 7)$, let $14\sim 13$. And $d(14)=3$. However, the edge $(12, 13)$ is sharing a $c_4:=12-13-11-4-12$ and a $C_3:=12-13-14-12$, a contradiction to Lemma \ref{noc3c4}. Thus $4\not\sim 11$. Let $4\sim 12, 13$. 
The Ricci-flat graph which contains an edge with endpoints degree $(3, 4)$ and $d(7)=4, 7\sim 5, 7\sim 6$ must have the following structure. 
 \begin{center}
\includegraphics[scale=0.5]{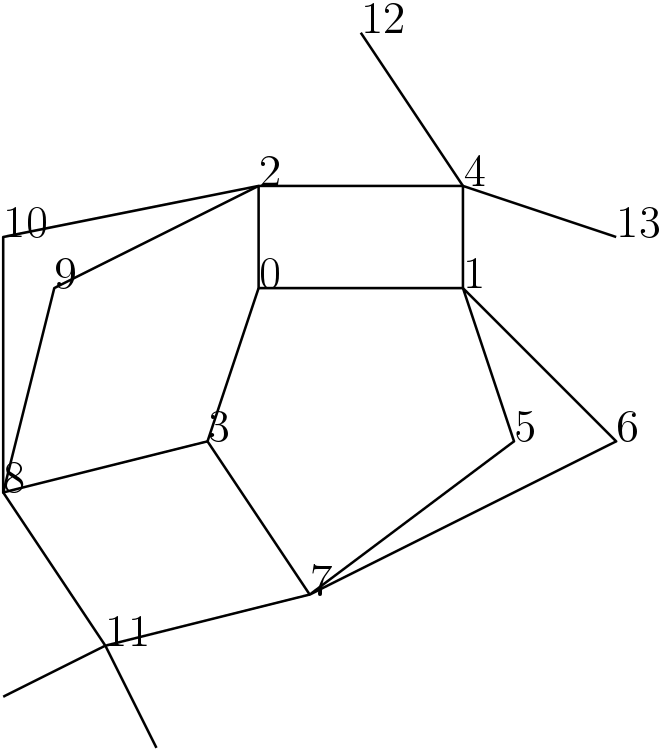} 
\end{center}

If $d(5)=3$ and let $z$ represent the third neighbor of vertex $5$,the $d(z)=4$. However, the edge $(5, 7)$ must be in exactly one $C_4$ that passes through edge $(5, z)$, similarly,  the edge $(1, 5)$ should also be in exactly one $C_4$  that passes through edge $(5, z)$, the edge $(5, z)$ share two $C_4$s , a contradiction to ``Type 5". Thus $d(5)=4$. Similarly, $d(6)=d(9)=d(10)=4$. 

Let $d(5)=4$, observe that we need $d(7, 9)=d(7, 10)=3$ and $d(2, 5)=d(2, 6)=3$, then $5, 6\not\sim 9$ and $5, 6\not\sim 10$.  Since any new neighbor of $5$ has distance $3$ to $0$, then we need that one of new neighbors of $5$ has distance $1$ to vertex $4$ for the edge $(1, 5)$, thus $5\sim 12$ or $5\sim 13$, wlog, assume the former.  Consider the edge $(1, 4)$, we need $d(6, 13)=3$. Then for the $(1, 6)$ we need $6\sim 12$ for the edge $(1, 5)$ and we need $d(5, 13)=3$ for the edge $(1, 4)$.  Let $5\sim 14$. Note $6\not\sim 14$ for the edge $(1, 5)$. Let $6\sim 15$. Consider the edge $(5, 7)$, we have $d(6, 1)=d(6, 12)=1$, $d(3, 1)=2$, 
$d(3, 12)=d(3, 14)=3$. Note $14\not\sim 11$, otherwise $d(6, 12)=1, d(1, 3)=2$ and $d(14, 11)=1$ would result $k(5, 7)\neq 0.8$. So we need $d(14, 11)=2$. Similarly, $d(11, 15)=2$ for the edge $(6, 7)$.  Let $11\sim 16\sim 14$. 

Now consider the edge $(4, 12)$ which is sharing two $C_4$s, thus $d(12)=4$. 
Since $d(2, 6)=d(13, 6)=3$, we need the fourth neighbor of $12$ to have distance $1$ to either $2$ or $13$. Assume $12\sim 9$, then consider the edge $(2, 4)$, then we need $d(10, 13)=3$. Observe that  the edge $(4, 13)$ cannot be in any $C_3$, if it is in a $C_4$, then we need $9\sim 13$. Consider the edge $(8, 9)$, we need $11\sim 13$. Then for $d(11, 15)=2$, we need $15\sim 16$. However, this situation is not good for the edge $(7, 11)$. Similarly, $12\not\ sim 10$. 
Thus, we need the fourth neighbor of $12$ to have distance $1$ to vertex $13$, which need new vertex $18\sim 12, 13$. 
Now consider the edge $(2, 4)$, we need at least one of $d(12, 9),  d(12, 10)$ is $2$. Wlog, let $d(12, 9)=2$, then $9\sim 18$. Then consider the edge $(2, 9)$, let $z$ be the fourth neighbor of $9$, since $d(0, z)=d(0, 18)=3$, $d(4, 18)=2$, then we need $d(4, z)=1$. Thus $9\sim 13$. However, this would result the edge $(13, 18)$ shares a $C_3:=13-18-9$ and a $C_4:=13-18-12-4$. A contradiction.  

Therefore, $d(5)\geq 5$. We should exclude ``Case 1b2".

\item{Case 1b3: $7\sim 9$, then $8\sim 6$ and $8\sim 10$.} For edge $(3, 7)$ in a $C_4$, we need new vertex $7\sim 11\sim 8$.  Note if $d(11)=2$, then $d(5)$ cannot be $2, 3, 4$. Similarly, $d(4)\neq 2$. If $d(11)=3$, note $11\not \sim 5, 6, 9, 10$.  Assume $4\sim 11$, then we get the  Ricci-flat graph $G_2$. Note $G_2$ is isomorphic to one of two graphs found in \cite{HLYY}. 

Next we claim there is no  Ricci-flat graphs with  $d(11)=3$ and $11\not\sim 4$:
\begin{proof}
Now we consider the case $d(11)=3$ and $11\sim 12$ a new vertex, then for the edge $(8, 11)$, since $d(7, 6)=d(7, 10)=3$, then $d(12, 6)$ and $d(12, 10)$ must be $2$. For the edge $(7, 11)$, since $d(8, 5)=d(8, 9)=3$, then $d(12, 5)$ and $d(12, 9)$ must be $2$.  Observe that $(11, 12)$ cannot be in any $C_4$, thus $d(12)\neq 4$. $d(12)$ must be $3$. 
Let us assume $d(5)=3$. Then the edge $(5, 7)$ must be in a $C_4$ which also must pass through vertex $9$. Let this  be $C_4:=5-7-9-13-5$, and since $d(1, 11)=3$, we need $d(13, 11)=2$, thus $13\sim 12$. Then consider the edge $(1, 5)$ which must also be contained in a $C_4$ and this $C_4$ can only pass through vertex $13$ and $6$ which implies that $6\sim 13$. Then the edge $(5, 13)$ are contained in two $C_4$s with endpoints degree $(3, 4)$, a contradiction. Thus For above graph, $d(5)=4$. By symmetric of vertices $5$ and $9$, $d(9)=4$, symmetric of vertices $7$ and $8$, $d(6)=d(10)=4$. 
Note when $d(5)=4$, Note $(5, 7)$ cannot be in any $C_3$,  we still need $5\sim 13\sim 9$ and $13\sim 12$. Similarly, the edge $(6, 8)$ must be in a $C_4$ which also must pass through vertex $10$. Let this  be $C_4:=6-8-10-z-6$, note the vertex $z$ cannot be $13$ since then $d(13)=5$. Thus we need a new vertex say $14$ as $z$. And we need $d(11, 14)=2$, thus  $14\sim 12$. Consider the edge $(12, 13), (12, 14)$, if $d(13$ must be $4$, by symmetry, $d(14)$ must be also $4$, then we need $d(9, 14)=2$ which implies $9\sim 10$. However, $d(7, 10)=2$, a contradiction. Thus $d(13)=4$, the edge $(12, 13)$ must be in a $C_4$ that passes through vertex $14$, then $d(14)=4$, let $13\sim 15\sim 14$. 
Now the graph has the following subgraph which is isomorphic to subgraph in the ``Case 1b2", thus we should exclude this case. 
 \begin{center}
\includegraphics[scale=0.3]{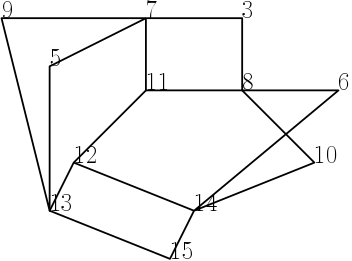} 
\end{center}
\end{proof}

In the following part, we consider the last case $d(11)=4$. By symmetry of vertex $4$ and vertex $11$, we have $d(4)=4$.  We will prove the Ricci-flat graphs exist and they are $G_3$ or $G_4$. 

\begin{itemize} 
\item For case $4\sim 11\sim 12$. Note $4\not\sim 12$ when consider edge $(4, 11)$. Let $13$ be the fourth neighbor of $4$. Then the edge $(4, 11)$ must  be in a $C_4$ which can only be $4-11-12-13-4$. Consider the edge $(1, 4)$, we have $d(5, 13), d(6, 13)$ must be $2$. Consider the edge $(7, 11)$, we have $d(5, 12), d(9, 12)$ must be $2$. Similarly, $d(9, 13), d(10, 13), d(10, 12), d(6, 12)$ are all $2$.

Consider the edge $(5, 7)$, since $d(5)\neq 2$, then $(5, 7)$ must be in a $C_4$ which must  pass through vertex $9$. Then we need a new vertex $14$  as a common for $5$ and $9$. Note $d(5)\neq 3$ since $d(1)\neq 3$ by Lemma \ref{34withonec4property}.  Let $15\sim 5$ be the fourth neighbor of $5$, for the edge $(5, 7)$, since $d(1, 3)=d(1, 11)=2$, then we need $d(15, 3)=2$ or $d(15, 11)=2$, obviously, it must be the latter case. Then $15\sim 12$.

For $d(5, 13)=2$,  note $13\sim 15$ would result edge $(12, 13)$ sharing $C_4$ and $C_3:=12-13-15$, which cannot happen. Thus $13\sim 14$. 
Then consider the edge $(1, 5)$ which must be in a $C_4$ and can only pass through edge $(1, 5)$. If $6\sim 14$, then the edge $(5, 14)$ share two $C_4$s but $d(15, 13)=2$, a contradiction. Thus $6\sim 15$. 

Now consider the edge $(6, 8)$ which must be in a $C_4$ that can only pass through edge $(8, 10)$. Thus we need a common vertex for $6$ and $10$.  Let $6\sim 16\sim 10$. Now consider the edge $(1, 6)$, since $d(8, 0)=d(8, 4)=2$ and $d(16, 0)=3$, we need $d(16, 4)=2$ which lead to $16\sim 13$.

Now consider the edge $(2, 9)$ which must be in a $C_4$ and can only pass through vertices $9-2-10-*-9$.  Since all vertices in the current graph have degrees at least $3$, we need new vertex $17$ as the common neighbor of $9, 10$. Now consider the edge $(7, 9)$, we need $d(11, 17)=2$, then it must be $17\sim 12$. Then $d(12)=4$. Consider the edge $(5, 15)$, by Lemma \ref{34withonec4property} we have $d(15)\neq 3$ as $d(12)\neq 3$. Then $d(15)=4$, similarly, $d(14)=4, d(16)=4, d(17)=4$. Observe that vertices $14, 15, 16, 17$ are not adjacent to each other, otherwise there would be  edges sharing $C_3$ and $C_4$. Let $15\sim 18$ a new vertex. 
Observe that the edge $(5, 15)$ can only be in one $C_4:=5-15-6-1-5$ and must be ``Type 6c", then we need $d(18, 14)=2$. Let $14\sim 19\sim 18$. Then consider the edge $(13, 14)$ which must be in a $C_4$ that passes through vertices $16$ and $19$, let $16\sim 19$. Consider the edge $(9, 14)$, we need $d(17, 19)=2$. There are two case, either we need a common vertex for $17, 19$ or $17\sim 18$. Then consider the edge $(12, 15)$ which must be in a $C_4$ that passes through vertices $16$ and $19$,we need $17\sim 18$. A new Ricci-flat graph is obtained: 
\begin{figure}[H]
\centering 
\includegraphics[scale=0.3]{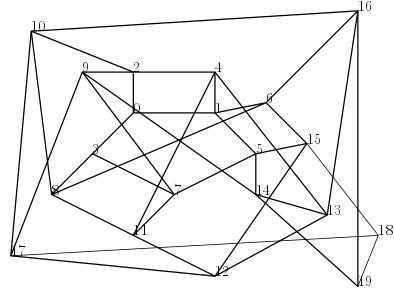}
\caption{Ricci-flat graph $G_3$}
\end{figure}

\item For case $4\not\sim 11$. 
Let $12, 13$ be the third and fourth neighbor of $11$.  If $d(5)=3$ and let $z$ represent the third neighbor of vertex $5$. Then the edge $(5, 7)$ must be in exactly one $C_4$ that passes through edge $(5, z)$, similarly,  the  edge $(1, 5)$ should also be in exactly one $C_4$  that passes through edge $(5, z)$, the edge $(5, z)$ share two $C_4$s , a contradiction to ``Type 5". Thus $d(5)=4$. Similarly $d(6)=d(9)=d(10)=4$. 

Consider the edge $(7, 11)$.  There are two cases: $d(5, 12)=1$ and  $d(9, 13)=3$; (2) $d(5, 12)=d(9, 13)=2$. 
\begin{itemize}
\item{Case (1):} When $5\sim 12$,  we claim $5\not\sim 13$. Otherwise we need $d(9, 12)=2$ for the edge $(5, 7)$. Then the $C_4$ that passes through edge $(7, 9)$ must pass through vertex $5$ then $12$ or $13$, which implies $9\sim 12$ or $9\sim 13$, a contradiction. 
Let $5\sim 14$, a new vertex.  We claim $9\sim 12$.
Also $d(1, 3)=2$ then $d(9, 14)\neq 1$ for the edge $(5, 7)$.  Since $d(9)=4$, then the edge $(7, 9)$ must be in $C_4$ that can only pass through vertices  $12$.   Thus $9\sim 12$.   Let $9\sim 15$. Now for the edge $(5, 7)$, consider the $C_4:-5-7-9-12-5$, we need $d(14, 11)=2$, but if we use $C_4:=5-7-11-12-5$, we need $d(9, 14)=2$. Similarly, for the edge $(7, 9)$, either $d(15, 11)=2$ or $d(15, 5)=2$. For $d(15, 11)=2$, either $15\sim 12$ then the edge $(9, 12)$ shares a $C_3$ and $C_4$, or $9\sim 13$, however, both cases cannot happen. Thus we need $d(15, 5)=2$, which implies $15\sim 4$. 
Since the edge $(5, 12)$ share two $C_4:=5-12-11-7-5, 5-12-9-7-5$, then $d(12)\neq 3$. Observe vertex $12$ cannot be adjacent to any existing vertices in the current subgraph. Let $12\sim 16$,a  new vertex.
 However, for the edge $(5, 12)$, $d(1, 9)=d(1, 11)=d(1, 16)=3$, we need $d(14, 16)=1$. Then for the edge $(9, 12)$, since $d(2, 11)=d(2, 16)=d(2, 5)=3$, we need $d(15, 16)=1$, this would result that the edge $(14, 16)$ shares a $C_3:=14-15-16-14$ and a $C_4:=14-16-12-5-14$.  A contradiction. 
Thus the edge $(7, 11)$ only shares exactly one $C_4$, similar analysis indicates each edges of $(8, 11), (1, 4)$, $(2, 4)$ shares exactly one $C_4$.

\item{Case (2)} For $d(5, 12)=2$, let $14$ be the common neighbor of $5, 12$. 
Let $15$ be the fourth neighbor of $5$.   The edge $(5, 7)$ must be in a $C_4$ which can only pass through $5-7-9-t-5$ where $t$ is either $14$ or $15$. 

We assume $15\sim 12$. 
Then by symmetry of vertices $14, 15$, wlog, let $9\sim 14$.  Note $9\not\sim 15$, otherwise for $d(9, 13)=2$, we need $14\sim 13$( similar for $15\sim 13$). Then we must have $6\sim 15$ for the $C_4$ that passes  through edge $(1, 5)$. However, this would result the edge $(5, 15)$ to share three $C_4$s, a contradiction. Thus $9\sim 15$.  Still consider $C_4$ that passes through edge $(1, 5)$.  Note $6\not\sim 14$ for the edge $(5, 14)$, thus $6\sim 15$. 
Then consider the fourth neighbor of vertex $9$,  we need new vertex as the fourth neighbor of $9$, let it be $16$. For the edge $(7, 9)$, we need $d(11, 16)=2$, then either $16\sim 12$ or $16\sim 13$. 
For $d(9, 13)=2$, either $14\sim 13$ or $16\sim 13$. Assume $14\sim 13, 16\not\sim 13$, then $16\sim 12$, then the edge $(12, 14)$ shares three $C_4$s, a contradiction. Assume $14\sim 13$ and $16\sim 13$. Now we have $d(14)=4$, then we need $10\sim 16$ for the $C_4$ that passes through edge $(2, 9)$. Then we need $d(4, 14)=2$ for the edge $(2, 9)$ which implies $4\sim 12$ or $4\sim 13$. 
Note $12\not\sim 4$ for the edge $(12, 14)$ and  $4\not \sim 13$ for the edge $(13, 14)$, a contradiction. 
Assume $14\not\sim 13$, then $16\sim 13$. Now consider the edge $(5, 14)$ which shares two $C_4$s, thus we need the fourth neighbor of $14$ to have distance $3$ to vertex $1$. For the edge $(2, 9)$,if $10\sim 16$, then we need $d(4, 14)=2$. For these two requirements, we need $4\sim 12$. Then we also need $4\sim 13$ for the $C_4$ passing through edge $(11, 12)$. Then consider the edge $(4, 13)$, we need the fourth neighbor of $13$ to have distance $2$ to vertex $1$. Then $13\sim 15$ which is not good for the edge $(5, 15)$ as it does not satisfy ``Type 6b".  The other case for the edge $(2, 9)$ is  $10\sim 14$. Then consider the edge $(1, 5)$, we need $d(4, 14)=2$ which implies $4\sim 12$ then $4\sim 13$. Then still consider the edge $(4, 13)$, we need $13\sim 15$, a contradiction.  
Thus $15\not\sim 12$. 

Then for the edge $(5, 7)$. There are two cases: either $9\sim 14$ or $9\sim 15$. 
\begin{itemize}
\item 
Assume $9\sim 14$. We need $d(15, 11)=2$ for the edge $(5, 7)$. Note $15\not\sim 12$, we must have $15\sim 13$.  
Consider the $C_4$ that passes through edge $(1, 5)$, assume $6\sim 14, 6\not\sim 15$. Then consider the $C_4$ passing through edge $(5, 15)$, we must have $9\sim 15$, a contradiction for the edge $(5, 14)$. Thus 
 Thus under the assumption  $9\sim 14$, we must have $13\sim 15$ and $6\sim 15$. Note $9\not\sim 15$ for the edge $(5, 15)$. Let $9\sim 16$. Consider the $C_4$ passing through vertices $2, 9, 10$, assume $10\not\sim 16$, then we need $10\sim 14$, however,  the edge $(9, 16)$ cannot be in any $C_3, C_4$, a contradiction. Thus $10\sim 16$. If $10\sim 14$, then consider the edge $(10, 14)$, we need $d(16, 5)=2$, which implies $16\sim 15$. Then for the $C_4$ passing through vertices $8, 6, 10$, we need $16\sim 6$ which is not good for the edges $(10, 16), (6, 16)$. Note $10\not\sim 15$ for the edge $(5, 15)$
  Thus we need a new vertex as the fourth neighbor of vertex $10$, let $10\sim 17$. 
 Still consider the $C_4$ passing through vertices $8, 6, 10$, either $6\sim 16$ or $6\sim 17$. 
If  $6\not\sim 17$, then $6\sim 16$, consider the edge $(6, 16)$, since both $d(1, 9), d(15, 9)$ are $3$, we need the fourth neighbor of vertex $16$ to have distance $1$ from vertex $15$. Thus $16\sim 13$. Then for the $C_4$ that passes through edge $(10, 17)$ we must have $17\sim 13$, however the edge $(13, 16)$ does not satisfy ``Type 6b",  a Contradiction. Thus, we need $6\sim 17$. 

Still for the edge $(8, 10)$, we need $d(16, 11)=2$, for the edge $(8, 11)$ we need either $d(6, 12)=2$ or $d(10, 12)=2$. Thus the we needs either $16\sim 12$ or $17\sim 12$. Since when  $5\sim 14\sim 12$, we have $6\sim 15\sim 13$, thus for $9\sim 14\sim 12$, and $10\sim 16$, we must have $16\sim 13$. 

For the edge $(1, 5)$, we need $d(4, 14)=2$. Note vertices $4, 14$ cannot be both adjacent to $16$ or $17$ by considering the edges $(2, 4), (1, 6)$ respectively. Then two cases are left:

Case a: If $4\sim 12$, then $4\sim 13$ in order to make $(11, 12)$ in a $C_4$. 
\begin{center}
\includegraphics[scale=0.3]{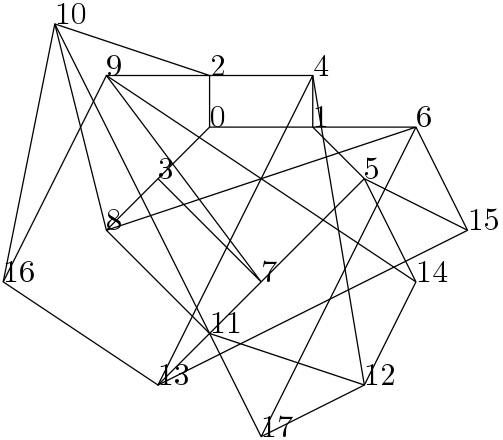}
\end{center}
Observe that the edge $(12, 14)$ must be in a $C_4$ that passes through edge $(12, 17)$, then we need a new vertex as the common of $14, 17$, let it be $18\sim 14, 17$. Then we need $d(15, 18)=2$ for the edge $(5, 14)$, let $15\sim 19\sim 18$. Observe that $d(18)$ cannot be $4$ since the new edge cannot be in any $C_3, C_4, C_5$. Similarly, $d(19)=3$. This is graph $G_4$. 

If we need a new vertex $18$ for $d(4, 14)=2$. Let $4\sim 18\sim 14$. Then $4\not\sim 12$, otherwise $4\sim 13$ which would make $d(4)=5$, a contradiction. Thus we need a new vertex as the fourth neighbor of vertex $14$. Let $4\sim 19$. 
For the edge $(5, 14)$, we need either $d(1, 12)=3, d(15, 18)=1$ or $d(1, 18)=2, d(15, 12)=2$. Note for the former case, we need $15\sim 18$, then the edge $(5, 15)$ does not satisfy ``Type 6". Thus we need the latter case.

 For $d(15, 12)=2$, note $12\not\sim 13$ as the edge $(11, 12)$ cannot satisfy ``Type 6".
We should also reject the case $12\sim 19\sim 15$. As the $C_4$ passing through edge $(4, 19)$ must pass through edge $4, 18$, then either $18\sim 15$ or $18\sim 12$ or $18$ is adjacent to a new vertex. The first case was rejected, the second case cannot guarantee the edge $(12, 14)$. For the third case, we have $d(2,15)=d(2, 12)=3, d(1, 15)=2, d(1, 12)=3$, then the edge $(4, 19)$ does not satisfy ``Type 6".

Thus we need a new vertex for $d(15, 12)=2$. Let $12\sim 20\sim 15$. Note vertex $16$ is not adjacent to vertex $12$ for the edge $(9, 16)$, then for the edge $(8, 11)$, if we need $d(10, 12)=d(6, 13)=2$, then it must be $17\sim 12$ and $17\sim 13$, which would result the edge $(13, 17)$ to share three $C_4$s. Then we need $d(10, 13)=d(6, 12)=2$ which only implies $17\sim 12$.  Then $d(17)=4$, let $z$ represent the fourth neighbor of vertex $17$.  Consider the edge $(10, 17)$, as $d(2, 12)=d(2, z)=3$ and $16\not\sim 12$, we need $16\sim z$. However, the edge $(10, 16)$ would share two $C_4:=10-16-9-2-10, 10-16-z-17-10$ but the third pair neighbors has $d(12, 13)=2$, a contradiction to ``Type 6b".

\item Next we consider the case when $9\sim 15$. 
Since we need $d(9, 13)=2$, there are two cases: Case 1 $15\sim 13$, we ignore this case since it can be compared when $5\sim 14\sim 12\sim 11,  9\sim 14$.

Case 2: $9\sim 16\sim 13$. Consider the edge $(2, 9)$ which should be in a $C_4$ that passes through vertex $10$, there are two cases: Case a: $10\sim 15$. Then we need $d(4, 16)=2$. 
Obviously, $d(15)\neq 3$. Let $d(15)=4$. For the edge $(1, 5)$  which should be in a $C_4$ that passes through vertex $6$, assume $6\sim 15$, then consider the edge $(6, 15)$, we need $d(5, z)=3$ where $z$ is the fourth neighbor of $6$. Note $6\not\sim 12, 13$ since the edge $(8, 11)$ should be in exactly one $C_4$, $6\not\sim 16$. Then we need a new vertex $17$  as the fourth neighbor of $6$. Consider the edge $(8, 10)$, we need the fourth neighbor of $10$ has distance $2$ to vertex $7$, since $10\not\sim 12, 13$, then $10$ must be adjacent to vertex $14$. However, this is not good for the edge $(5, 15)$. Thus under the condition $10\sim 15$, we have $6\not\sim 15$. Then for the edge $(1, 5)$  to be in a $C_4$, we need $6\sim 14$. Consider the edge $(8, 6)$, note $10\sim 14$ is not good for the edge $(5, 14)$, then we need a new vertex for the edge $(8, 6)$ to be in a $C_4$, let it be $10\sim 17\sim 6$. Consider the edge $(5, 15)$, note $d(15)\neq 3$, since both $d(10, 1)=d(9, 1)=3$. However, $15\not\sim 12, 13, 16, 17$. Let $15\sim 18$, then we need $14\sim 18$, however, this is not good for the edge $(5, 14)$. 
Now we consider the last one. 
Case b: $10\sim 16$ . For the edge $(1, 5)$, assume $6\sim 15$. Note $d(15)\neq 3$, since both $d(5, 8)=d(9, 8)=3$. Then $d(15)=4$, note $15\not\sim 16$ for the edge $(9, 15)$. Assume $15\sim 13$, then consider the edge $(2, 9)$, we need $d(4, 15)=2$, then $4\sim 13$. Then for the edge $(11, 13)$ to be in a $C_4$, we must have $4\sim 12$. Consider the edge $(4, 12)$, since both vertex $14$ and the fourth neighbor of $12$ have distance $3$ to vertex $0$, and $d(1, 14)=2$, a contradiction. Thus $15\sim 17$ a new vertex and $4\sim 17$. Consider the edge $(6, 15)$, we need $d(8, 17)=2$, then $17\sim 10$, the edge $(2, 4)$ is in two $C_4$, which is rejected before. Thus we have arrived another contradiction. 
\end{itemize}
\end{itemize}
\end{itemize}
\end{itemize}
\end{enumerate}

\end{proof}

Next we will reject the other case ``Type 5b" for the edge with endpoint degrees $\{3, 4\}$. 
\begin{theorem}\label{34withonec42}
Let G be a  graph  with maximum degree at most $4$, let $(x, y)$ be an edge in G with $d(x)=3,  d(y)=4$ that satisfies ``Type 5b". Then $G$ cannot be a Ricci-flat graph. 

\end{theorem}
\begin{proof}
 Assume $G$  contains an edge $(0, 1)$ with $d(0)=3,  d(1)=4$, $(0, 1)$ is contained in exactly one $C_4:0-2-4-1-10$, and $d(3, 5)=2, d(3, 6)=3, d(2, 6)=2$. 
By contradiction, we assume $G$ is Ricci-falt. As this is the last case of ``Type 5", then for any other edge $(x', y')$ in $G$ with endpoint degrees $\{3, 4\}$, there must be an isomorphism from the neighborhood $N(0)\cup N(1)$ to neighborhood $N(x')\cup N(y')$. 
Let $5\sim 7\sim 3$, $2\sim 8\sim 6$. Since $d(3, 6)=3$,  we need a new vertex as the third neighbor of $3$. Note $2\not\sim 7$, otherwise $(0, 2)$ would be in two $C_4$s. Thus we need a new vertex as the fourth neighbor of $2$. Let $2\sim 9$. $3\sim 10$.
Note for the edge $(0, 2)$ since $(d(0), d(2))=(3, 4)$, we also have either $d(3, 9)=2, d(3, 8)=3, d(1, 8)=2, d(1, 9)=2$ or $3$ or $d(3, 8)=2, d(3, 9)=3, d(1, 9)=2, d(1, 8)=2$ or $3$. The latter case is not preserved under isomorphism then we consider the former case. For $d(3, 9)=2$, either $9\sim 10$ or $9\sim 7$, since the edge $(0, 3)$ share two separate $C_5$s, we have $9\sim 10$. 
\begin{figure}[H]
\centering
\includegraphics[scale=0.4]{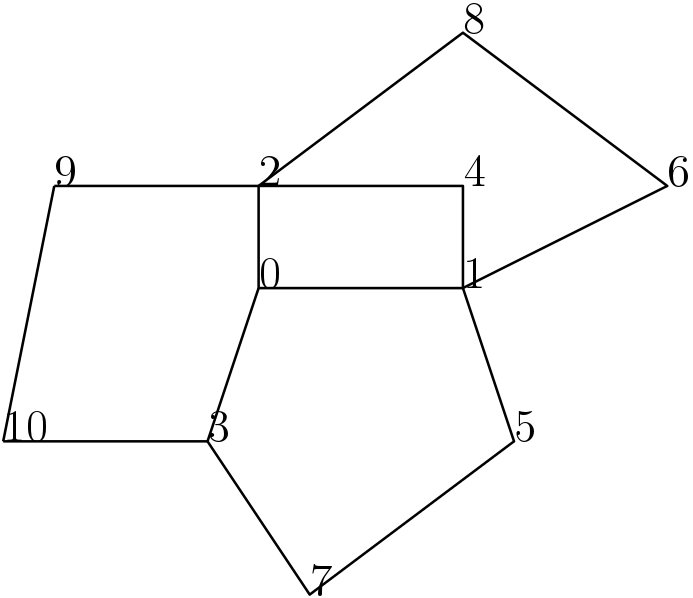} 
\label{type5b-iso}
\end{figure} 

Note $7\not\sim 6, 8$, $10\not\sim 6, 8$. If $d(7)=3$, then the edge $(3, 7)$ shares two $C_5$s, which implies $d(10)=3$. Similarly, if $d(10)=3$, then $d(7)=3$. We will first exclude the following case: 
\begin{itemize}
\item 
If $d(10)=4$, then the edge $(3, 10)$ must be in a $C_4$ which  must pass through vertex $7$. Then $d(7)=4$. 
Assume $10\sim 5$ and $10\sim 11$, a new vertex. Then we need $d(7, 11)=2$. Let $7\sim 12\sim 11$. Consider the fourth neighbor of $7$, note $7\not\sim 9$. Otherwise  the edge $(3,10)$ would be contained in two  $C_4$, a contradiction.  Then $7\sim 13$, a new vertex. Then consider the edge $(3, 7)$, we have both $d(0, 12), d(0, 13)$ are $3$, a contradiction. Thus  for the $(3, 10)$ to be in a $C_4$ that passes through vertex $10$, let $10\sim 11\sim 7$, $7\sim 12$ and $10\sim 13$ and we still need $d(7, 13)=2$, for the edge $(3, 7)$, we need $d(10, 12)=2$. Let $12\sim 13$. 
\begin{center}
\includegraphics[scale=0.4]{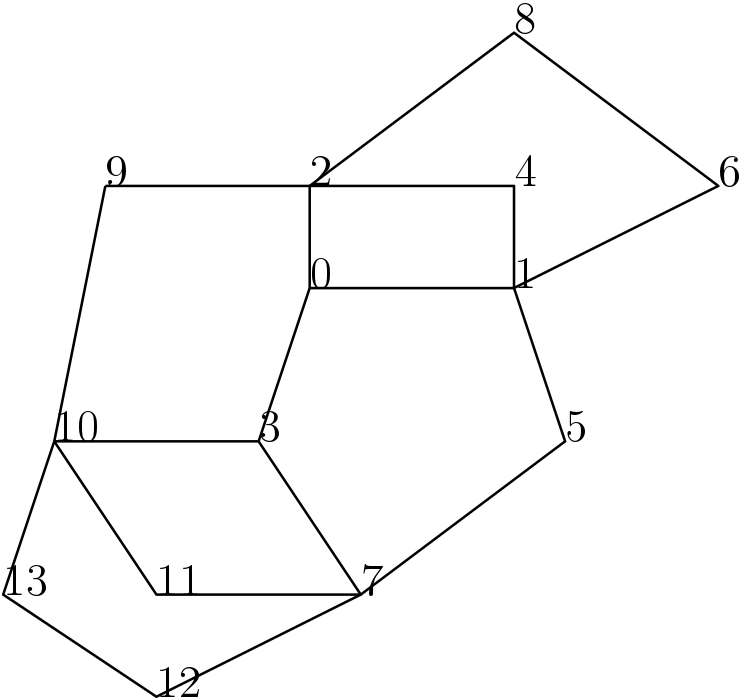}
\end{center}
Consider $d(5)$, if $d(5)=2$, we need one of $d(1, 12)=2$, $d(1, 11)=2$ for the edge $(5, 7)$, note we also need one of $d(7, 6)=2$, $d(7, 4)=2$ for the edge $(1, 5)$.  By symmetry of the current graph, there are three different cases for the two conditions. Case 1: Assume $12\sim 4$ and $12\sim 6$, then $d(12)=4$, the edge $(7, 12)$ must be in a $C_4$ that passes through vertex $11$, then either $11\sim 6$ which is not good for the edge $(6, 12)$, or $11\sim 4$ which is not good for the edge $(1, 4)$. Case 2: Assume $12\sim 4$, $11\sim 6$. Then still the edge $(7, 12)$ must be in a $C_4$ that passes through vertex $11$, then  $11\sim 4$ 
which is not good for the edge $(1, 4)$. Case 3: $11\sim 4$. If $d(4)=3$, then $d(11)=3$, the current graph is all good except the edge $(6, 8), (12, 13)$, however, there is no way to generate a Ricci-flat graph. Thus $d(4)=4$. Since $4\not\sim 6, 8, 12, 13$, we need a new vertex as its fourth neighbor, let it be $14$. Note $11\not\sim 14$, let $11\sim 15$. Since we need the edge $(4, 11)$ to be in a $C_4$, it must be $14\sim 15$. Then we need $d(6, 14)=d(8, 14)=d(12, 15)=d(13, 15)=2$ for edges $(1, 4), (2, 4), (7, 11), (10, 11)$ respectively. Then $d(6)\neq 2$, however the edge $(1, 6)$ cannot be in any $C_4$, a contradiction. Thus $d(5)\neq 2$. 
If $d(5)=3$, let $z$ represent the third neighbor, then the edge $(5, 7)$ must be in exactly one $C_4$ passing through edge $(5, z)$,  the edge $(1, 5)$ should also be in exactly one $C_4$ passing through ege $(5, z)$,  then the edge $(5, z)$ is in two $C_4$s, a contradiction to Lemma \ref{34notwoc4}. Thus $d(5)=4$. Similarly $d(9)=4$. 

Consider the edge $(1, 5)$, if it is contained in a $C_3$, then $5\sim 6$, then we need the fourth neighbors of vertex $5$ and vertex $6$ to have distance $3$ from vertex $4$ . Assume $5\sim 12$, then $d(12)=4$. 
For the edge $(5, 6)$, assume $d(8, 7)=2$, however neither $8\sim 12$, otherwise the $(5, 6)$ would share a $C_3$ and $C_4$, nor $8\sim 11$ considering the edge $(5, 7)$. Similarly, $13\not\sim 6, 4$. Assume $d(8, 12)=2$, consider the case $8\sim 13$. Note if the edge $(6, 8)$ is in a $C_3$, it must be $C_3:=6\sim z\sim 8$, where $z$ is the a new vertex, however, we have $d(1, 2)=d(5, 13)=2$, a contradiction. Thus edge $(6, 8)$ must be in a $C_4$. Let $6\sim 14$, assume $14\sim 13$. Then the edge $(12, 13)$ must be in a $C_4$ that passes through the fourth neighbor of $12$,let it be $15$. Note $15\not \sim 8$ for the edge $(8, 13)$, then it must be $15\sim 14$. Note we need the edge $(10, 13)$ to be in a $C_4$, however, there is no $C_4$ satisfying the current structure. Thus for the edge $(6, 8)$ to be in a $C_4$, it must be $14\sim 15\sim 8$. Similarly we need $12\sim 16\sim 17\sim 13$ for edge $(12, 13)$ to  be in a $C_4$. Then consider the edge $(8, 13)$, it must be in the  $C_4:=8-13-17-15-8$. Then there is no $C_4$ for the edge $(10, 13)$. Thus $8\not\sim 13$.
Hence, for $d(8, 12)=2$, let $8\sim 14\sim 12$. For the edge $(6, 8)$ in a $C_4$, let $6\sim 15\sim 16\sim 8$. 
Then consider the edge $(2, 8)$ which must be in a $C_4$, assume the $C_4$ passes through vertex $4$, then $d(4)=4$. Consider the edge $(1, 4)$ since each of the third and fourth neighbors of $4$ has distance $3$ to vertex $5$, then one of them must have distance $1$ to 
vertex $6$, then the edge $(1, 6)$ shares a $C_3$ and $C_4$, a contradiction. Then the $C_4$ for edge $(2, 8)$ must pass through vertex  $9$. Assume $9\sim 14$. Observe that since $4\not\sim 14, 15$, then any neighbor of vertex $4$ has distance $3$ to vertex $5$, in addition, we have $d(2, 5)=3$, then $d(4)\neq 3$. Let $d(4)=4$ with the third and fourth neighbors $w, l$, then we need $d(8, w)=d(9, l)=2$. However,  $d(8, w)=2$ cannot be true. 

Thus the edge $(1, 5), (5, 7)$ must be in $C_4$s. Assume the $C_4$ passes through vertex $4$, then $d(4)=4$, let $5\sim 14\sim 4$. Then we need the fourth neighbor of vertex $5$ to have distance $2$ from vertex $6$. 
Assume $5\sim 13$,  then $d(13, 6)=2$. Consider the edge $(5, 13)$, since both $d(10, 1)=3, d(12, 1)=3$, then $d(13)\neq 3$. Let $d(13)=4$, we need the fourth neighbor of $13$ to have distance $1$ to vertex $14$. Note if $13\sim 4$, we need $d(13, 6)=3$ for the edge $(1, 4)$, a contradiction. Assume $13\sim 8$, then we need $d(14, 10)=2$ for the edge $(5, 13)$ and $d(14, 11)=2$ for the edge $(5, 7)$. Thus $14\sim 9$ for $d(14, 10)=2$, let $14\sim 15\sim 11$ for $d(14, 11)=2$. Then vertex $4$ must be adjacent to a new vertex $16$. Note $d(6, 16)=d(8, 16)=3$ for the edges $(1, 4), (2, 4)$ respectively, then for the edge $(16, 4)$ to be in any $C_4$, it must pass through vertex $9$, that is $9\sim 16$. However, it is not good for the edge $(2, 9)$.  
Thus let $13\sim 15\sim 14$. For the edge $(5, 7)$, we need $d(14, 11)=2$.  
 Assume $11\sim 15$, then we need $d(5, 9)=2$ for the edge $(10, 13)$ which implies $14\sim 9$. Consider the edge $(9, 10)$ which must be in a $C_4$ that passes through the fourth neighbor of $9$. Note $9\not\sim 15$ for the edge $(9, 15)$, if $9\sim 12$, then $6\not\sim 12$ or $6\not\sim 13$ for the edge $(10, 13)$, a contradiction to $d(13, 6)=2$. Thus for the edge $(9, 10)$ to  be in a $C_4$, we need $9\sim 16\sim 11$. However, the edge $(10, 11)$  is contained in three separated $C_4$s, a contradiction. Thus  there is no way to form a $C_4$ for the edge $(9, 10)$. Then for $d(14, 11)=2$, we assume $11\sim 4$. Then the edge $(4, 11)$ must be in a $C_4$ that passes through vertex $14$, let $14\sim 16\sim 11$, however, we cannot guarantee the edge $(4, 14)$. Then for $d(14, 11)=2$, assume $14\sim 9$. Note $11\not\sim 4$ or $11\not\sim 15$. Then still for the edge $(9, 10)$ to  be in a $C_4$, we need $9\sim 16\sim 11$. Then  $d(11)=4$, let $z$ be the fourth neighbor of $11$, we need $d(11, 13)$ for the edge $(10, 11)$, However, there is no way to form a $C_4$ for the edge $(10, 13)$. A contradiction. Thus   for $d(14, 11)=2$, we need a new vertex $16$ with $14\sim 16\sim 11$. However, now the edge $(5, 14)$ is in two $C_4$s, but the third pair of neighbor $(16, 7)$ which are not on the $C_4$s, has distance $d(16, 7)=2$, a contradiction. 
 
 Thus the $C_4$ for the edge $(1, 5)$ does not passes through vertex $4$, then it must pass through vertex $6$. Then the edges $(1, 6)$ must be in a $C_4$, so do edges $(7, 12), (2, 8), (10, 13)$.  Let $5\sim 14\sim 6$. Then we need the fourth neighbor of vertex $5$ to have distance $2$ from vertex $4$. Similarly, we need the fourth neighbor of vertex $5$ to have distance $2$ from vertex $11$.  Note $14\sim 8$ otherwise, the edge $(6, 14)$ would share a $C_3$ and $C_4$. Thus for the edge $(2, 9)$ to be in a $C_4$, let $9\sim 15\sim 8$. Then we need the fourth neighbor of vertex $9$ to have distance $2$ from vertex $4$. Thus $5\not\sim 9$. 
 Assume $5\sim 8$, then for the edge $(5, 7)$ to be in a $C_4$, it must be $C_4:=5-7-12-14$. Note we need $d(8, 11)=2$, since $11\not\sim 15$, then $11\sim 6$. Consider the edge $(6, 8)$,  we need $d(11, 15)=2$. For the edge $(5, 8)$, we need $d(7, 15)=2$, thus $15$ must be adjacent to vertex $12$ and $11\sim 16\sim 15$. Since for the edge $(1, 6)$, we need $d(4, 11)=2$, then $4\sim 16$. However, we cannot guarantee the edge $(6, 11)$. 
 
 Thus $5\not\sim 8$, similarly, $5\not\sim 13$, $9\not\sim 6$ and $9\not\sim 12$. If $5\sim 15$, we need a new vertex $16$ for $d(4, 15)=2$, let $15\sim 16\sim 4$, since $d(15)$ is at most $4$, then for $d(11, 15)=2$, we need $11\sim 16$. Now observe that the $C_4$  for the edge $(5, 15)$ must pass through vertex $14$. Note $14\sim 9$ is not good for the edge $(9, 15)$, then it must be $14\sim 16$, also we need $d(7, 9)=2$, which implies $9\sim 11$ or $9\sim 12$, contradiction. Thus $5\not \sim 15$. Similarly, $9\not\sim 14$. 
 
 Let $5\sim 16$. For $d(16, 4)=2$, assume $4\sim 12\sim 16$,  since we also need $d(z, 4)=d(w, 4)=2$ where $z, w$ represent the neighbors of vertices $6, 8$ respectively, then the fourth neighbor of vertex $4$ should be adjacent to vertices $4, z, w$, then the fourth neighbor of vertex $4$ should have degree at most $1$ in the current graph, which implies it must be a new vertex. Let $4\sim 17$. Note $d(17)\geq 3$, then the edge $(4, 17)$ must be in a $C_4$ that can only pass through vertex $12$, then either $17\sim 16$ or $17\sim 13$. Assume $17\sim 16$, then we need $d(2, 13)=2$ for the edge $(4, 12)$, which implies $8\sim 13$. Then $17\sim 13$ for $d(13, 4)=2$, however the edge $(8, 13)$ cannot be in any $C_4$, a contradiction. Assume $17\sim 13$, then we need $d(2, 16)=2$ for the edge $(4, 12)$, note $16\sim 8$ implies $17\sim 8$, thus we need $16\sim 9$. Note we need $d(4, 11)=2$ for the edge $(7, 12)$, then $11\sim 17$. However, we cannot guarantee the edge $(10, 11)$. 
Thus for $d(16, 4)=2$, we need a new vertex $17$ such that $4\sim 17\sim 16$. 

Now consider the fourth neighbor of vertex $6$. Assume $6\sim 12$, then $12\sim 17$, however the edge $(5, 7)$  cannot be in any $C_4$ through vertex $12$. Assume $6\sim 13$, then $13\sim 17$, however, we cannot guarantee the edge $(6, 13)$. Note $6\not\sim 15$, otherwise the edge $(8, 15)$ would share a $C_3$ and $C_4$. Assume $6\sim 16$, then the edge $(6, 16)$ is in two $C_4:=6-16-5-1-6, 6-16-5-14-6$, thus $d(16)=4$. However, $16\not\sim 12$ for the edge $(5, 16)$, thus we need $14\sim 12$ for the edge $(5, 7)$ in a $C_4$. Note for the edge $(6, 16)$, we need $14$ to have  distance $2$to the fourth neighbor of $16$, but we also need the fourth neighbor of vertex $14$ has distance $3$ to vertex $16$ considering edge $(5, 14)$, a contradiction. Thus the fourth neighbor of vertex $6$ should be a new vertex $18$, and $18\sim 17$. 

Assume $8\sim 19\sim 17$ for the edge $(1, 8)$. Then $d(17)=4$, consider the edge $(4, 17)$ 
which must be in a $C_4$ that passes through vertex $16$. Now consider the fourth neighbor of vertex $4$. Assume $4\sim 12$, then $16\sim 12$, then we need $d(2, 13)=2$ for the edge $(4, 12)$ which implies $13\sim 9$, a contradiction. Assume $4\sim 13$, then $16\sim 13$, then we need $d(1, 12)=2$ for the edge $(4, 12)$ which is not possible. Note $4\not\sim 14$ and $4\not\sim 15$. Thus we need a new vertex $4\sim 20\sim 16$. Then consider the edge $(5, 16)$ which must be in a $C_4$ that passes through vertex $14$. If $14\sim 20$, then we need the fourth neighbor of vertex $16$ to have distance $2$ to vertex $7$, then $16\not\sim 12$, we must have $14\sim 12$ for the edge $(5, 7)$  to be in a $C_4$. However, this is not good for the edge $(5, 14)$. Assume $17\sim 14$, then we need $16\sim 12$  for the edge $(5, 7)$  to be in a $C_4$. However, this is not good for the edge $(5, 16)$. For all the other possible cases such that the edge $(5, 16)$ is in a $C_4$(that is $d(14, 16)=2$), we must have $4\sim 12$ consider the edge $(5, 7)$ in a $C_4$. Then the edge $(5, 14)$ shares two separated $C_4$s, we need the fourth neighbor of $14$ to have distance $3$ from vertex $16$, a contradiction. 

Assume $8\sim 18$ for the edge $(1, 8)$. The the edge $(6, 8)$ shares a $C_3:=6-8-18-6$. Thus we need $d(14, 15)=3$.  Consider the edge $(6, 14)$, note the new neighbor of vertex $14$ should have distance $2$ from vertex $8$, then it must be adjacent to vertex $15$ which would make $d(14, 15)=2$, a contradiction. 

\item We consider the case $d(7)=d(10)=3$. Let $7\sim 11, 10\sim 12$, then $11\sim 12$ for the edge $(3, 7)$. Consider the edge $(3, 7)$ and edge $(3, 10)$, we have $d(5)=d(11)$ and $d(9)=d(12)$ by Lemma \ref{twoc5property}. 
There are two cases for $d(4)$. 

Case 1:  assume $d(4)=3$. Note $4\not\sim 11$, otherwise there is an isomorphism  $\phi$ that maps the current subgraph to subgraph \ref{type5b-iso}. In the following, we call the latter graph as `` base graph" We must have $\phi(4)=0, \phi(0)=4$. There are two cases for $\phi$ depending on $\phi(1)$.  If $\phi(1)=1$, then follow the isomorphism we must have $\phi(10)=12, \phi(9)=9$. Since $10\sim 9$ in base graph, we must have $12\sim 9$, a contradiction. 
Similar analysis for case $\phi(1)=2$. 
Thus the neighborhood of edge $(1, 4)$ is not preserved under isomorphism. Similarly, $4\not\sim 12$. Note vertex $4$ cannot be adjacent to any other existing vertices. 
Thus let $4\sim 13$,  a new vertex. Consider the isomorphism $\phi$ from current graph to base graph, we have $d(13)=d(3)=3$. 
Let $13\sim 14, 15$,  then $d(14)=d(7)=3, d(15)=d(10)=3$. If $\phi(1)=1$, then $\phi(5)=6$ and $\phi(6)=6$. Since $\phi(14)=7$, then we must have $14\sim 5$ as $7\sim 5$ in the base graph. Similar analysis for $15\sim 9$.  

Consider the neighbors for vertices $14$ and $15$. Note $14\not\sim 6$ or $14\not \sim 8$ or or $14\not \sim 9$ for edges $(1, 4), (2, 4), (13, 14)$ respectively. Assume $14\sim 11$. Then $\phi(11)=11$. Consider the edge $(13, 14)$, we must have $15\sim 12$ and $\phi(12)=12$. 
Then there is a $C_4:=5-14-11-7-5$, thus $d(5)=d(11)=4$. 

Now we consider a new isomorphism $\phi$ between the neighborhood of edges $(7, 5)$ and $(0, 1)$. Let $\phi(5)=1, \phi(7)=0$, then $\phi$ is determined.  We need the fourth neighbors of $5$ to be mapped to $6$ and the fourth neighbor of $11$ to be mapped to $8$. 
Consider vertices $6$ and $8$, note $6\not\sim 5$ for edge $(1, 5)$, similarly, $8\not\sim 9$. If $5\sim 8$, then $11\sim 6$. However, we cannot guarantee the edge $(1, 6)$. Thus $5\not\sim 8, 11\not\sim 6$. Similarly, $9\not\sim 6, 12\not\sim 8$, then all vertices $5, 6, 8, 9, 11, 12$ need new vertices, however, the edges $(1, 6), (2, 9)$ are not good under this situation as they cannot be in any $C_3$ or $C_4$ and  $d(6)=d(8)=2$ is not true.  A similar analysis for case $\phi(5)=2$.
Thus  $11\not\sim 14$.    Similarly, $12\not\sim 15$. 

Assume $14\sim 12$. Still consider the isomorphism $\phi$ between the current graph to the base graph that satisfies $\phi(1)=1, \phi(4)=0$.
Since $14\sim 12$, then $\phi(12)=11, \phi(11)=12$, we must have $phi(10)\sim \phi(12)$, i.e $15\sim 11$.
Then we consider a new isomorphism $\phi$ between the neighborhood of edges $(7, 5)$ and $(0, 1)$.
Similarly we cannot guarantee the edge $(1, 6)$.

Consider the last case: let $14\sim 16$ and $15\sim 17$. Under the isomorphism $\phi$ between the current graph to the base graph that satisfies $\phi(1)=1, \phi(4)=0$, we must have $\phi(16)=11, \phi(17)=12$, then $16\sim 17$. If $d(5)=3$, then the edge $(5,7)$ must share two $C_5$s which needs $11\sim 16$ and $12\sim 17$. However, the edge $(1, 5)$ cannot be in any $C_4$, a contradiction. 
Thus we conclude that $d(4)=4$.


Case 2: $d(4)=4$. Note $4\not\sim 11$, otherwise $4\sim 12$ by symmetry, then $d(11)=d(12)=4$ and they are adjacent to two new vertices by ``Type 6a". Let $11\sim 13, 12\sim 14$. We need $d(2, 13)=d(2, 14)=2$. Then the edge $(7, 11)$ must be in a $C_4$ passing through vertices $5$ and $13$. Let $5\sim 13$, similarly, let $9\sim 14$. Then we consider the isomorphism $\phi$ between the current graph to the base graph with $\phi(7)=0$ and $\phi(11)=1$. Then we have $\phi(3)=3, \phi(12)=5, \phi(4)=6, \phi(2)=12$, since $2\sim 4$, then we have $12\sim 6$, a contradiction. A similarly analysis for the isomorphism $\phi$ with $\phi(7)=0$ and $\phi(11)=2$.

 As $4$ cannot be adjacent to any other existing vertices, then let $4\sim 13, 14$. 
 Consider the edge $(1, 4)$, if it satisfies ``Type 6b", then wlog, either $5\sim 13$ or $6\sim 13$. 
Case 2a: Assume $5\sim 13$, then $d(6, 14)=3$. 
 
 If further $d(5)=3$, then $6\not\sim 13$. Observe the edge $(1, 6)$, it cannot be in any $C_3$ or $C_4$, thus $d(6)=2$. Then we have $d(8)=4$ considering the edge $(6, 8)$. Then under the isomorphism between neighborhoods of edges $(1, 5)$ and $(0, 1)$, we must have $13\sim 8$ and need a new vertex $15$ such that $11\sim 15\sim 13$. For the edge $(2, 8)$, we need the fourth neighbor of $8$ to have distance $2$ from vertex $9$ by ``Type 6b".  Note  $8\not\sim 14$ or $9\not\sim 14$ considering edge $(1, 4)$ and $(2, 4)$ respectively. 
Then the edge $(4, 14)$ can only be in a $C_5$ which implies $d(14)=2$ and the fourth neighbor of $14$ has distance $2$ from vertices $2$ and $13$. Note $8\not\sim 11, 12$ considering edge $(4, 14)$. Then we need a new vertex $16$ as the fourth neighbor of vertex $8$ such that $d(16, 9)=2$. Then the edge $(2, 9)$ which cannot be in any $C_3$ or $C_4$, then $d(9)=2$, however, we then cannot guarantee $d(16, 9)=2$. A contradiction. 
Thus $d(5)=4$.  The edge $(5, 7)$ must be in a $C_4$ passing through vertex $11$, implying $d(11)=4$. Note if $13\sim 11$, one cannot find isomorphism from the current graph to the base graph. Thus we need a new common vertex of $5$ and $11$. 
 Note $5\not\sim 12$ as under any isomorphism between neighborhood of $(0, 1)$ and $(7, 5)$, vertex $12$ is mapped to vertex $4$, vertex $3$ is mapped to vertex $4$, but the $d(3, 12)\neq d(3, 4)$. 
One can check vertex $9$ or $8$ or $14$ is not the common for vertices $5$ and $11$.
%
Thus we need a new vertex $15$ as the common for vertices $5$ and $11$. Let $5\sim 15\sim 11$. 

 Consider two different maps between neighborhood of $(0, 1)$ and $(7, 5)$ specified by whether vertex $5$ is mapped to vertex $1$ or $2$. For the latter isomorphism $\phi$. See Figure \ref{illustrationforType5b} for illustration. 
 Observe that $4$ is not mapped to vertex $6$, otherwise $4\sim 11$ as $11$ is mapped to vertex $1$ and $1\sim 6$, then $d(4)\geq 5$, a contradiction. Thus we need a new vertex as the image of vertex $6$, let it be $16$, then $11\sim 16$, we also have $13\sim 16$ as $13$ is mapped to vertex $8$ and $8\sim 6$. Note vertices $12, 9$ should share one more vertex by the isomorphism. Assume $14$ is this common vertex. As vertices $9$ and $12$ cannot have other common vertex except $10, 14$, then vertex $14$ must be mapped to vertex $15$. Then vertex $4$ is mapped to vertex $14$. Since $4\sim 14$, then $14\sim 15$ by isomorphism, again $15\sim 4$ as $14, 15$ are mapped to $15, 4$. Then $d(4)\neq 5$, a contradiction.  It is easy to check no existing vertices can be the common for vertices $9$ and $12$. Thus we need a new vertex $17$ and $9\sim 17\sim 12$. Then $\phi(4)=17$. Since $4\sim 13$, then $17\sim 8$ by $\phi$ and $\phi(17)=15$.
 Consider vertex $8$, note $d(8)\neq 3$ as there is no isomorphism between neighborhood of edges $(2, 8)$ and $(0, 1)$. Then $d(8)=4$. Consider the edge $(8, 17)$, which is in the $C_4=8-17-9-2-8$. Note $12\not\sim 8$, otherwise vertex $8$ is the image of vertex $13$, since which implies $12\sim 13$ as $12, 13$ are the images of vertex $2, 8$ and $2\sim 8$. This would result $12\neq 5$, a contradiction. Similar analysis for $12\not\sim 13$ or $12\not\sim 14$ or $12\not\sim 6$. Thus we need a new vertex as the fourth neighbor of vertex $12$, let $12\sim 18$. Then $\phi(18)=3$ and $\phi(8)=18$. Since  $\phi(4)=15$, then $15\sim 18$. 
 Now consider the edge $(11, 12)$ which satisfies ``Type 6c". Then we need $d(16, 17)=2$, then $d(\phi(16), \phi(17))=d(15, 6)=2$. Then consider the edge $(1, 5)$, we must have $6\sim 13$ and $d(4, 15)=2$ which implies $14\sim 15$. Thus $\phi(14)=14$, and $17\sim 14, 18\sim 16$. 
 
Then consider the edge $(1, 6)$ which satisfies ``Type 6c". We need a common vertex for $14$ and the fourth neighbor of vertex $6$. 
Similarly, for the edge $(2, 4)$ which satisfies ``Type 6c", we need a common vertex for $14$ and the fourth neighbor of vertex $9$.
Thus vertices $6$ and $9$ should share the fourth neighbor, otherwise $d(14)\geq 5$. Let $9\sim 19\sim 6$ and $14\sim 19$. We also have $\phi(19)=16$. Since $19\sim 14$, we have $\phi(19)\sim \phi(14)$, i.e. $16\sim 14$, then $d(14)\geq 5$. A contradiction. 
\begin{figure}[H]
\centering
\includegraphics[scale=0.4]{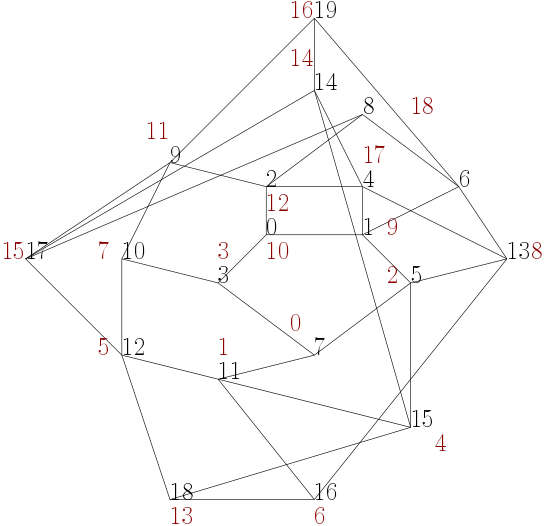}
\caption{The red labels represent the image of black labels}
\label{illustrationforType5b}
\end{figure}

Case 2b: The case  $6\sim 13$ is included in above case. We omit the proof here.

\end{itemize}
We have rejected all cases based on the structure stated in the theorem, thus there is no Ricci-flat graph of ``Type 5b".  

\end{proof}

\begin{theorem}\label{33twoC5}
Let G be a Ricci-flat graph  with maximum degree at most $4$, if there exists an edge $(x, y)$ with endpoint degree $(d(x), d(y))=(3, 3)$, and $(x, y)$ is contained in exactly two separate $C_5$s. Then $G$ is one of the following:
the dodecahedral graph,
the Petersen graph,
the half-dodecahedral graph and
the Triplex graph. 
\end{theorem}

\begin{proof}

See the following subgraph of $G$, let $x=0, y=1$. 

\begin{center}
\includegraphics[scale=0.3]{33C5C5.png}
 \end{center}

By Theorem \ref{23determined}, if $d(2)=d(3)=2$ then $d(4), d(5)$ must be $3$, and the graph is Half-dodecahedra, refer to Theorem 2. If  $d(2)=d(3)=4$ or  $d(4)=d(5)=4$, refer  to Theorem \ref{34withonec41}. 
Thus we only need to consider the following case: 
$d(2)=d(3)=3$, $d(4)=d(5)=3$.
Since $(0,3)$ cannot be in any $C_4$, then $3$ is not adjacent to $6$ or $4$, similarly,  $2$ is not adjacent to $7$ or $5$, $4$ is not adjacent to $7$ and $5$ is not adjacent to $6$.  Thus vertices $2, 3, 4, 5$ need new neighbors, let $8, 9$ be  the neighbors of vertices $2, 3$ respectively. 
If $d(6)=2$ or $d(7)=2$, then the edge $(2, 6)$ or $(3, 7)$ satisfy the condition in  Theorem \ref{23determined}. Thus $G$ is the half-dodecahreal graph. 
If $d(6)=4$ or $d(7)=4$, then the edge $(2, 6)$ or $(3, 7)$ satisfy the condition in  Theorem \ref{34withonec41}. Thus $G$ is determined. 

Then we consider the case when $d(6)=3, d(7)=3$.  Note all edges in the current structure are not in any $C_4$. Consider the edge $(2, 8)$ if it is in a $C_4$, then the $C_4$ must  pass through the edge $(2, 6)$, a contradiction. Thus $d(8)=3$. The edge $(2, 8)$ cannot be in any $C_4$, then the third neighbor of vertex $8$ has degree $3$. Similarly,  $d(9)=3$ and its thrid neighbor also has degree $3$. Following this process, there is no vertex with degree $4$, then no edge in any $C_4$. Refer to \cite{LLY},  then $G$ must be one of the dodecahedral graph,
the Petersen graph, the Triplex graph.

\end{proof}

\section{Ricci-flat graphs with vertex degrees $\{2, 4\}$}

In previous sections, we have finished all cases when an edge has endpoint degrees $(2, 2), (2, 3), (3, 3), (3, 4)$. 
In this section, we consider the Ricci-flat graphs  that contain edges with endpoint degrees $(2, 4), (4,  4)$,  we first classify these that contain a copy of $C_3$.
\subsection{Ricci-flat graphs that contain $C_3$}
\begin{theorem}\label{thm:44C3}
Let $G$ be a Ricci-flat graph with maximum degree at most $4$ and there exist an edge $e=(x, y)$ with $d(x)=d(y)=4$ such that $e$ is contained in a $C_3$. Then the graph is isomorphic to graphs $G_5, G_6$ and $G_7$. 
\begin{figure}[H]
\centering
\includegraphics[scale=0.4]{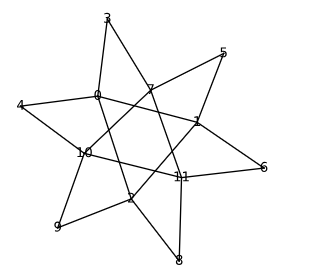}
\hfil
\includegraphics[scale=0.4]{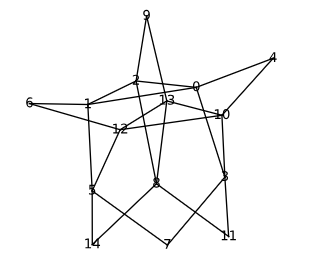}
\caption{Ricci-flat graph $G_5, G_6$}
\label{44Ricci-flat1}
\end{figure} 

\begin{figure}[H]
\centering
\includegraphics[scale=0.5]{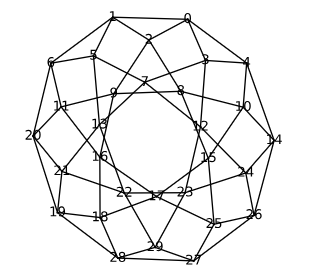}
\hfil
\includegraphics[scale=0.5]{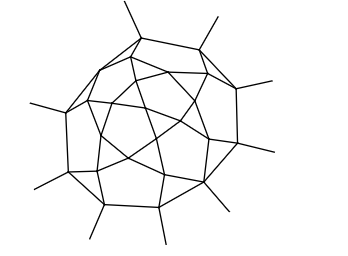}
\caption{Ricci-flat graph $G_7$ and its local structure}
\label{44Ricci-flat2}
\end{figure} 
\end{theorem}

\begin{remark}
In Geometry, $G_7$ is  a polyhedron called icosidodecahedron with $20$ triangular faces, $12$ pentagonal faces,  and $30$ vertices, $60$ edges. 
\end{remark}

\begin{proof}
 Look at the follow subgraph with $x=0, y=1$. Since $(0, 2)$ is in $C_3$, then $d(2)=4$. Note vertex $2$ cannot be adjacent to any existing vertices, let $8, 9$ be  the its new neighbors. Note we will exclude the case when there is an edge with endpoints degree $3, 4$. 
 \begin{center}
\includegraphics[scale=0.4]{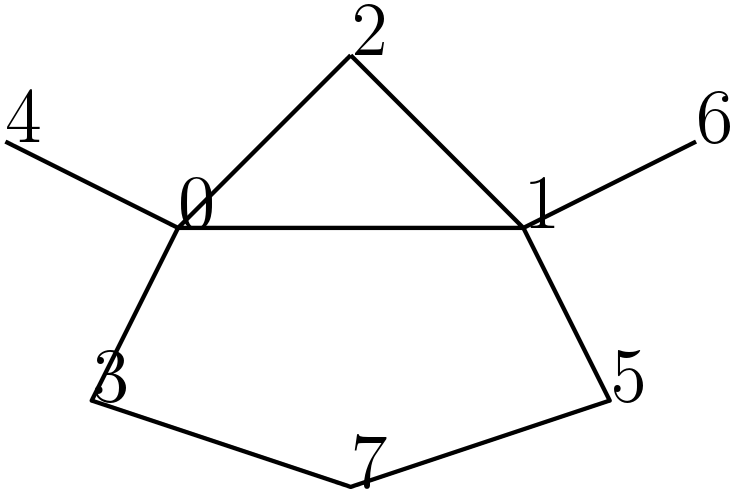}
\end{center}
 
 \begin{itemize}
 
\item  When $d(3)=2$, then $d(7)=4$. Then at least one of $d(7, 2)=2, d(7, 4)=2$ is true consider edge $(3, 7)$. 
\begin{itemize}
\item

     Assume $d(7, 2)=2$, wlog, let $7\sim 8$. Then  $d(4, 9)=3$ for the edge $(0, 2)$ and  $d(6, 9)=3$ for the edge $(1, 2)$. In this case, $7\not\sim 9$.  Otherwise, we have $d(5, 8)=3$ for the edge $(1, 2)$ and $d(8)=4$,  let $8\sim 10, 11$, since $d(3, 10)=d(3, 11)=3$, we need $d(5, 10)=1$ or $d(5, 11)=1$ which would result $d(5, 8)=2$,  a contradiction. 
Note the edge $(0, 4)$ cannot be in any $C_3$ or $C_4$, thus $d(4)=2$, and the edge $(0, 4)$ is in a $C_5$. 
     Assume $d(5)=2$, then similarly, $d(6)=2$. 
Let $4\sim 10$, note for the edge $(0, 4)$ to be in a $C_5$, if $d(10, 1)=2$ then  $10\sim 6$ which would result $d(4, 6)=2$, a contradiction. Thus we need $d(10,3)=d(10, 2)=2$, then $10\sim 7$ and $10\sim 8$. Since vertices $8, 10$ are in a $C_3$, then $(8)=d(10)=4$. Then the edge $(2, 8)$ must be in a $C_4$ that passes through vertices: $C_4:=2-8-9-11-2$.
Now consider the edge $(1, 6)$ since $d(6)=2$, and $6$ cannot be adjacent to any existing vertices, then let $6\sim 12$, however, none of $d(12, 5), d(12, 2), d(12, 0)$ can be $2$. A contradiction. 
    Let $d(5)=4$. For the edge $(0, 4)$ to be in a $C_5$, if we still assume $10\sim 7$ and $10\sim 8$. Then  $(8)=d(10)=4$.
By same reasons as above, let $8\sim 11\sim 9$, $10\sim 12$. 
    Observe that the edge $(5, 7)$ cannot be in any $C_3$ or $C_4$, a contradiction. Thus for the edge $(0, 4)$ to be in a $C_5$, we must have $10\sim 5$ and exactly one of $10\sim 7$ and $10\sim 8$. Assume $10\sim 7$. Note $5\not\sim 6$. Let $5\sim 11$, $10\sim 12$. Then the edge $(1, 5)$ must be in the  $C_4:=11-5-1-6-11$. Consider the vertex $8$, we claim 
$d(8)=2$, otherwise $d(8)=4$. then the edge $(7, 8)$ must be in a $C_4$ which cannot happen. Now consider the edge $(2, 9)$ 
which cannot be in any $C_3, C_4$ or $C_5$. Thus $d(5)\neq 4$. A contradiction. 
Thus $7\not\sim 8$. 

\item So we need $d(7, 4)=2$ for the edge $(0, 3)$. Let $4\sim 10\sim 7$. Note $7\not\sim 8, 9$.  Observe the edge $(0, 4)$ cannot be in any $C_3$ or $C_4$, thus $d(4)=2$. 
Assume $d(5)=4$. For the edge $(0, 4)$, we assume $10\sim 5$. For the edge $(1, 5)$, if it is contained in a $C_3$, then $5\sim 6$. For vertices $7, 10$, let $7\sim 11, 10\sim 12$. Consider vertex $6$, we have $d(6)=4$,  since  $6$ cannot be adjacent to any existing vertices, then let $6\sim 13, 14$. Then for the edge $(5, 6)$ we need either $d(7, 13)=2$ then $d(6, 11)=2$ which is not good for the edge $(5, 7)$;   or $d(10, 13)=2$  then  $d(6, 12)=2$ which is not good for the edge $(5, 10)$. 
Thus vertex $5$ needs a new vertex as its fourth neighbor. Let  $5\sim 11$.  Observe that the edge $(1, 5)$ must be in a $C_4$, then $6\sim 11$. 
Consider the edge $(1, 5)$, either $d(7, 2)=2$ or $d(10, 2)=2$, since $7\not\sim 8$, then we need $10\sim 8$. For the edge on $C_3:=5-10-7-5$, we need $d(11, 12)=d(11, 8)=d(12,8)=3$. Then $d(11)=2$. Then the edge $(8, 10)$ cannot be in ant $C_3$ or $C_4$, we have $d(8)=2$. Consider the edge $(2, 9)$, it cannot be in any $C_3, C_4$, then $d(9)=2$, however, the second neighbor of $9$ cannot have distance $2$ to vertices $8$ and $0$, which is a contradiction for the edge $(2, 9)$ to be in two $C_5$s. Thus $d(5)\neq 4$.
 
Let $d(5)=2$, then similarly, $d(6)=2$. Let $6\sim 11$. Thus since at least two of $d(11, 5), d(11, 5), d(11, 0)$ are $2$. Then it must be $11\sim 7$. Now we have $d(10)=d(11)=4$. Then for the edge $(0, 4)$, we must have $d(10, 2)=2$ which implies $10\sim 9$. 
Note the edges $(7, 10), (7, 11)$ cannot be in any $C_4$, they must be in $C_3$, thus $10\sim 11$. Then for the edge $(1, 6)$, we must have $11\sim 8$. Observe that $d(8), d(9)$  must be $2$. The resulting graph is $G_5$.
\end{itemize}

\end{itemize}

\item  Now we let $d(3)=d(5)=4$. 
 For the edge $(0, 3)$,  if it is contained in a $C_4$ then must pass through vertex $4$. There are two cases: 
 \begin{itemize}
 \item[Case 1:] $C_4:=3-0-4-7-3$. Let $3\sim 10, 11$, then for the edge $(0, 3)$, we need $d(10, 2)=d(11, 1)=2$. For $d(10, 2)=2$, wlog, let $10\sim 8$.  
    For $d(11, 1)=2$, note for the edge $(0, 1)$, $11\not\sim 6$. Thus let $11\sim 5$. Note for the edge $(0, 1)$, $7\not\sim 6$, for the edge $(3, 7)$, $7\not\sim 8$, or the edge $(0, 2)$, $7\not\sim 9$. Let  $7\sim 12$. Consider the edge $(0, 4)$, it is in the $C_4:=0-4-7-3-0$ and a $C_5:=0-4-7-5-1-0$, thus $d(4)=4$. For the edge $(0, 4)$, we need that the new neighbors of $4$ have distance $2$ to vertices $1, 2$ respectively. Since $d(4, 6)=d(4, 9)=3$, then we need that the new neighbors of $4$ have distance $1$ to vertices $8, 5$ respectively.
    
    Assume $4\sim 10$ and $4\sim 11$. Then vertices $10, 11$ are not adjacent to any existing vertices. Let $10\sim 13, 11\sim 14$. For the edge $(3, 10)$, we need $d(7, 13)=2$. Note $13\not\sim 12$ for the edge $(4, 7)$, thus $13\sim 5$. Then consider the edge $(1, 5)$, it must be in a $C_4$ which passes through vertex $6$. Then it must be $6\sim 13$. Then for the edge $(5, 13)$, we need the fourth neighbor of vertex $13$ has distance $2$ to either vertex $7$ or $11$. By symmetry of $7, 11$, wlog, let $13\sim 12$. 
 Then consider the edge $(1, 2)$, we have both $d(5, 8)=3, d(5, 9)=3$, a contradiction. 
 
   Now assume $4\sim 13\sim 5$ where $13$ is a new vertex. Consider the edge $(1, 5)$ which must be in a $C_4$ that passes through vertex $6$. Then either $6\sim 11$ or $6\sim 13$, however, both cases are not good for the edge $(0, 1)$. 
   Thus $4\not\sim 10$.
   
   Let $4\sim 13\sim 8$. Then assume $4\sim 11$. By symmetry of $7, 11$,  vertex $11$ is not adjacent to $6, 8, 9, 10, 13$,
   if $11\sim 12$, then consider the vertex $5$, $5\not\sim 13$ for the edge $(4, 7)$, $5\not\sim 10$ for the edge $(3, 7)$. Let $5\sim 14$. For the edge $(5, 7)$, we need either $d(4, 14)=2$, then $13\sim 14$ which is not good for the edge $(4, 11)$, or $d(1, 12)=2$, then $6\sim 12$. Then for the edge $(5, 11)$, we need $d(14, 12)=2$, then $12\sim 13$, however, this situation is not good for the edges $(4, 7)$ and $(4, 11)$. 
   Thus we need new vertex as the fourth neighbor of $11$, let $11\sim 14$. Then consider the edge $(1, 5)$, it must be in a $C_4$ which passes through vertex $6$ and new neighbor of $5$, since $6\not\sim 13$ for the edge $(0, 1)$,  let $5\sim 15\sim 6$. Consider the edge $(1, 2)$, assume $d(5, 8)=2$, then $8\sim 15$. Consider the edge $(2, 8)$ which must be in a $C_4$ that passes through vertex $9$. However the vertex $9\not\sim 10, 15, 13$ for the edge $(0, 2)$ or $(1, 2)$, a contradiction. Thus for the edge $(1, 2)$, we need $d(5, 9)=2$, then $9\sim 15$. Then consider the edge $(1, 5)$, we need either  $d(2, 7)=2$ or $d(2, 11)=2$, both cannot happen. A contradiction. 
   For the requirement that new neighbors of $4$ has distance $1$ to vertices $5$. We let $4\sim 14\sim 5$.  Then consider the edge $(1, 2)$, assume $d(5, 8)=2$, then either $8\sim 11$ or $8\sim 14$.   Note $9\not\sim 10, 13, 14, 11$ for the edge $(0, 2)$, then edge $(2, 8)$ cannot be in any $C_4$. Thus for the edge $(1, 2)$, we need $d(5, 9)=2$, then either $9\sim 11$ or $9\sim 14$, a contradiction for the edge $(0, 2)$.

\item[Case 2:]
$C_4:=3-0-4-10-3$, where $10$ is a new vertex. Let $11$ be the fourth neighbor of $3$. 
A similar analysis for the edge $(1, 5)$, it must be in a $C_4:=1-5-12-6-1$. 
We need $d(2, 11)=2$ for the edge $(0, 3)$ which lead to $d(8, 11)=1$. For the edge $(0, 2)$, we have $d(3, 8)=2$ thus $d(4, 9)=3$. Then we consider a similar situation for the edge $(0, 2)$, we have $d(8)=4$ and $(2, 8)$ is in a $C_4$ that pass through vertex $9$. Let $8\sim 13\sim 9$. By symmetry, let $5\sim 14\sim 8$. 

Consider the edge $(0, 4)$, assume $d(4)=4$, let $4\sim z, w$. Then we need $d(1, z)=d(2, w)=2$, which means $4$ should be adjacent to these vertices: $7, 11, 13, 12, 14$. If $4\sim 7$, by symmetry, $6\sim 7$, then $d(4, 6)=2$, a contradiction. Similar analysis, we would get contradictions for all other cases. 

Thus $d(4)=2$, then $d(10)=2$. By symmetry, $d(6)=d(9)=2, d(12)=d(13)=4$. Consider the edge $(3, 10)$, let $10\sim z, w$, we need $d(7, z)=d(11, w)=2$. This fact also implies that the edges $(3, 7). (3, 11)$ cannot be in any $C_3$ or $C_4$, then $d(7)\neq 4, d(11)\neq 4$. Note for $d(7, z)=d(11, w)=2$, $10\not\sim 14$, otherwise by symmetry $7\sim 13, 11\sim 12$, a contradiction.  Thus it must be $10\sim 12, 13$, then consider the edge $(5, 12)$, we need the fourth neighbor of $12$ has distance $2$ to vertex $14$, similarly, for the edge $(8, 13)$, we need the fourth neighbor of $13$ has distance $2$ to vertex $14$. By symmetry $12\sim 13$. The resulting graph is $G_6$ and it is Ricci-flat. 

\item[Case 3:]
Now we consider the case when the edge $(0, 3)$ is in a $C_3$,  which must be $C_3:=0-3-4-0$. Similarly, 
the edge $(1, 5)$ must be in a $C_3$ which must be $C_3:=1-5-6-0$ and the edge $(2, 8)$ must be in a $C_3$ which must be $C_3:=2-8-9-2$. Let $3\sim 12, 5\sim 13$. We must have $12\sim 7$ and $13\sim 7$ for the edge $(0, 3)$ and $(1, 5)$, etc. 
\begin{figure}[H]
\centering
\includegraphics[scale=0.4]{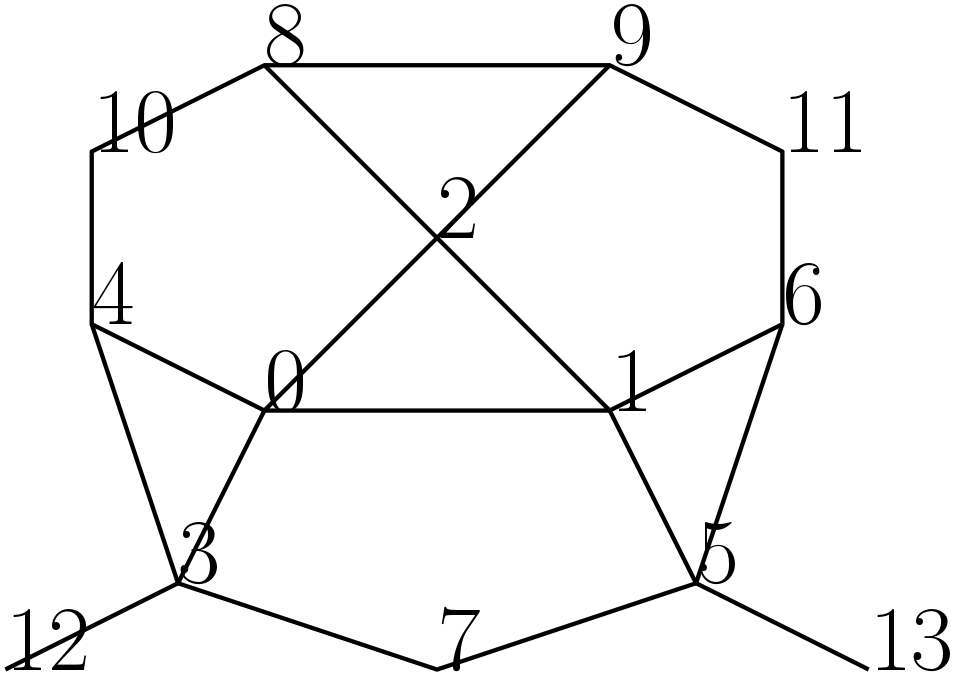}
\hfil
\includegraphics[scale=0.4]{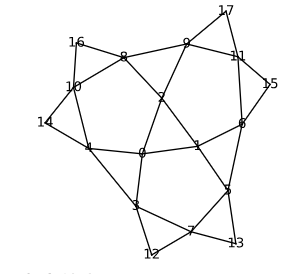}
\end{figure}

Following this structure, we are actually constructing a simple graphs in which every vertex has degree $4$ and every edge shares one $C_3$ and one $C_5$, then eventually the graph is isomorphic to $G_7$. 

\end{itemize}
\end{proof}

\subsection{Ricci-flat graphs that contain $C_4$}

Now all Ricci-flat graphs containing $C_3$ are determined. In the following section we will only consider the Ricci-flat graphs containing $C_4$s.  Basically these can be classified into three cases: 
\begin{itemize}
\item all $C_4$s are vertex-disjoint; 
\item exist two $C_4$s sharing an edge;  
\item no two $C_4$s sharing an edge, but exist two $C_4$s sharing a vertex. 
\end{itemize}

We consider each case. 
\begin{theorem}\label{vertex-disjoint}
Let $G$ be a Ricci-flat graph with maximum vertex degree $4$, assume $G$ contains no edges with endpoint degree $(3, 4)$, and all $C_4$s of G are vertex-disjoint. Then $G$ is isomorphic to $G_8$. 
\begin{figure}[H]
\centering
\includegraphics[scale=0.4]{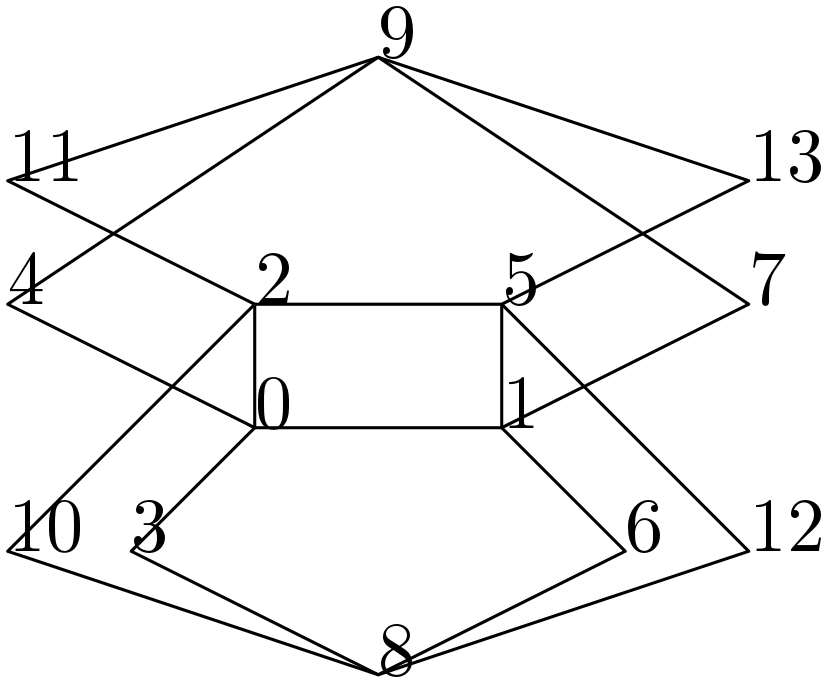}
\caption{Ricci-flat graph $G_8$}
\label{44Ricci-flat3}
\end{figure}
\end{theorem}

\begin{proof}
Let $G$ be a Ricci-flat graph in class $\mathcal{G}$ with the following sub-structure,  where $d(0)=d(1)=4$ and the edge $(0, 1)$ is in the $C_4:=0-1-5-2-0$. 
\begin{center}
\includegraphics[scale=0.3]{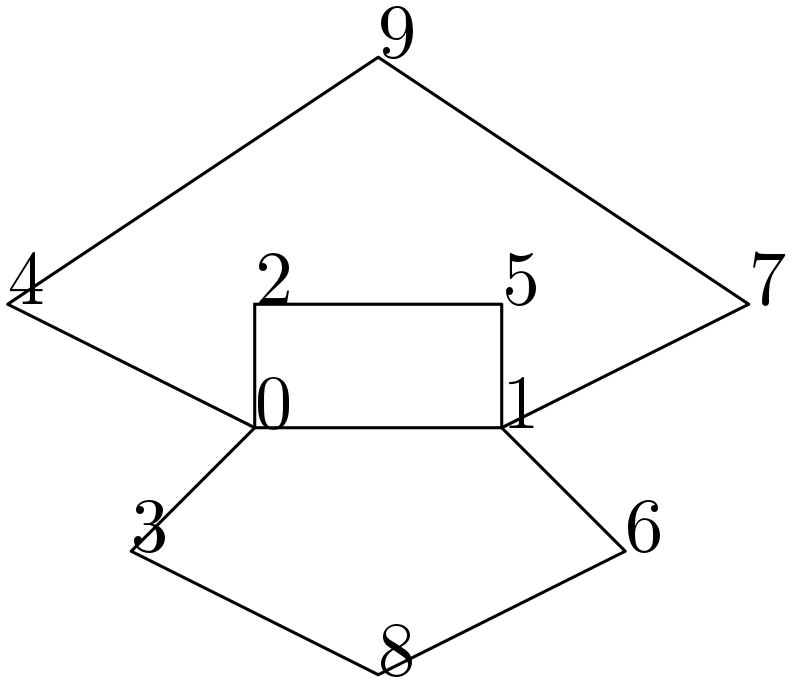}
\end{center}
Since $C_4$s in $G$ are vertex-disjoint, then $8\neq 9$, and the edges $(0, 2), (0, 3), (0, 4)$ are not in any other $C_4$s, which implies $d(3), d(4)$ must be $2$. Similarly, $d(5)=d(6)=2$. Then $d(8)=d(9)=4$. Consider vertex $2$, it is not adjacent to $8$ or $9$.  Let $2\sim 10, 11$. Since the edge $(2, 5)$ is not in any $C_4$ or other $C_4$, let $5\sim 12, 13$. Then wlog, we must have $d(10, 12)=d(11, 13)=2$ and we must make connection with vertices $8, 9$. Then $10, 12$ must be adjacent to $8$ or $9$, wlog, let $10\sim 8\sim 12$. Then $11\sim 9\sim 13$. Now the graph is $G_8$, this  is the unique Ricci-flat graph such that $C_4$s are vertex-disjoint. Note $G_8$ is isomorphic to the second graph found in \cite{HLYY}. 

\end{proof}

In the following parts, we will see many infinite Ricci-flat graphs, we call these as {\it primitive graphs}, as from each, we can obtain finite Ricci-flat graphs which are obtained by a vertex-preserve projection from the primitive graph. 

For the second and third cases, we first consider the case when a vertex in some $C_4$ has degree $2$, see the following results:

\begin{theorem}\label{thm:C424}
Let $G$ be a Ricci-flat graph with maximum degree at most $4$.  Let edge $(x, y)\in E(G)$ be in a $C_4$  with $d(x)=4$ and $d(y)=2$. Then the $G$ is isomorphic to graphs with the primitive graphs showing in the Figures \ref{24Ricci-flat1},  \ref{24Ricci-flata}, \ref{24Ricci-flat8}, \ref{24Ricci-flat9}. 
\end{theorem} 
\begin{proof}
See the following graph, let $x=0, y=1$, $C_4:=0-1-5-2-0$, then $d(3, 5)= d(4, 5)=3$ by ``Type 3".   Since edge $(1, 5)$ is in $C_4$, then $d(5)$ must be $4$. 
\begin{center}
\includegraphics[scale=0.5]{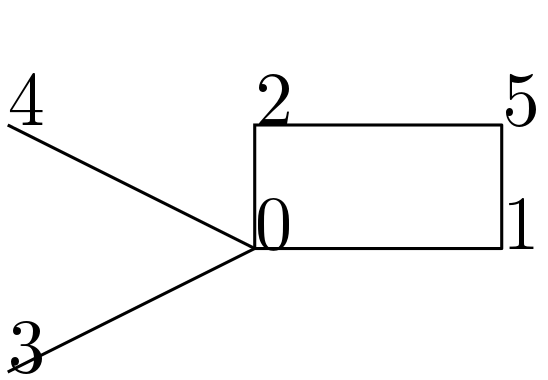} 
\end{center}

We focus on vertex $2$,  note we will exclude the case when  an edge has endpoints degree $3, 4$.  Thus there are two cases:
 \begin{itemize}

\item[Case 1] Assume $d(2)=2$. Let $5\sim 6, 7$. Note $d(0, 6),  d(0, 7)$ must be $3$. 
\begin{center}
\includegraphics[scale=0.5]{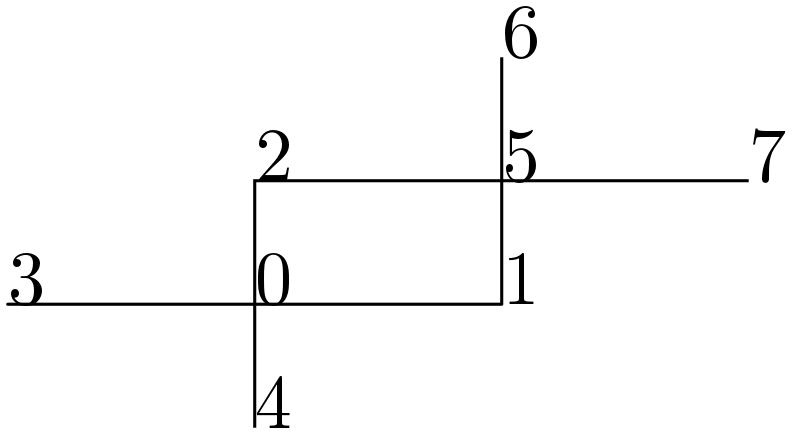} 
\end{center}

Assume  $d(6)=3$, then the edge $(5, 6)$ must be contained in a $C_4$  passing through vertex $7$
let it be $C_4:=6-5-7-8-6$, and $9$ be the third neighbor of $6$, so we need one of $d(9, 1), d(9, 2)$ to be $2$ which cannot happen. 

Assume  $d(6)=4$, then the edge $(5, 6)$ is  in a $C_4$ which must pass through vertex $7$, still we need the neighbors of $6$ have distance $2$ from vertex $1$ and vertex $2$, which cannot happen. 

Thus  $d(6)=2$. Similarly,  $d(7)=2$.  We assume the edge $(5, 6)$ must be in a $C_5$ passing through vertex $7$ and new vertices, let this $C_5:=6-5-7-9-8-6$. However we need  then $d(8, 2)$ or $d(8, 1)$ to be $2$ which cannot happen.
Thus the edge $(5, 6)$ must be in a $C_4$ passing through vertex $8$, $C_4:=6-5-7-8-6$. Now compare the $C_4:=0-1-5-2-0$ and $C_4:=5-7-8-6-5$, the situation is same as the initial stage, we have that 
$d(8)$ must be $4$ and its neighbors must have degree $2$.  Then  we can extend this process infinitely times to get an infinite graph  that consists of a sequence of $C_4$s  where the vertices degrees are $2, 4, 2, 4$ in the cycle order and  each two consecutive $C_4$s sharing one vertex of degree $4$.  To get a finite graph, the only way is to combine the neighbors of two degree $4$ vertices only if the distance of these two vertices is at least  $3$.  For example, we could have $8\sim 3, 4$, then get the smallest Ricci-flat graph of this type, we could also let $8\sim 9, 10, 11\sim 9, 10, 3, 4$.  Thus there are infinitely many such Ricci-flat graphs.  

\begin{figure}[H]
\centering
\includegraphics[scale=0.4]{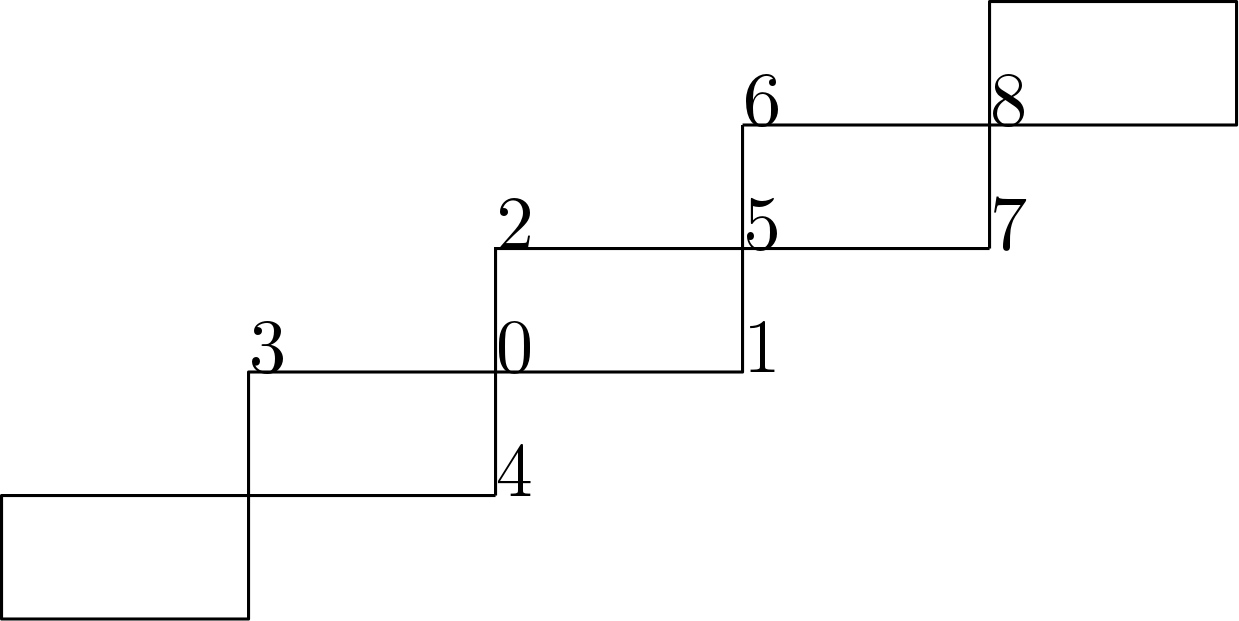} 
\hfil
\includegraphics[scale=0.4]{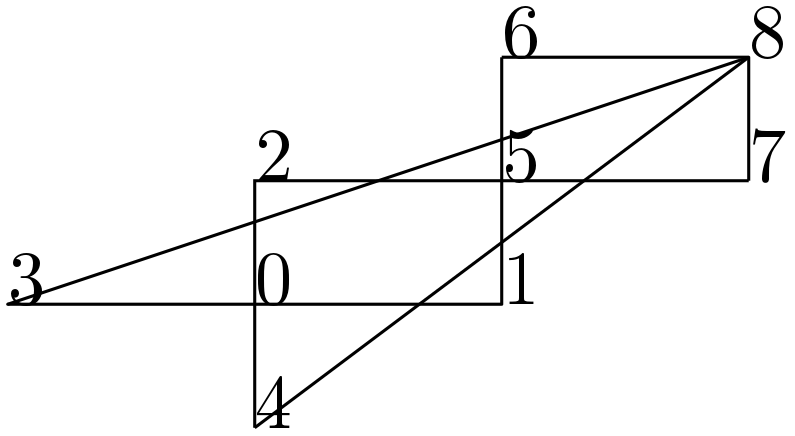} 
\caption{Ricci-flat graphs: the primitive graph and  an example of finite type}
\label{24Ricci-flat1}
\end{figure}

\item[Case 2]  Assume $d(2)=4$. Let vertices $6, 7$ be neighbors of $2$, vertices $8, 9$ be neighbors of $5$. Note $d(0, 8)=d(0, 9)=3$. 
We claim the edge $(0, 3), (0, 4)$ must be in $C_4$. Suppose $(0, 3)$ is in a $C_5$, then $d(3)=2$, let $3\sim 10\sim 6$. We have $d(3, 6)=2$, thus $4\not\sim 7$ for the edge $(0, 2)$. Assume $d(4)=2$, then we must have $d(4, 6)=2$, note $4\not\sim 10$, otherwise the edge $(0, 3)$ would be in a $C_4$. Let $4\sim 11\sim 6$. Consider the edge $(0, 4)$, by ``Type 3" and the fact $d(11, 1)=3$, we also need $d(11, 3)=2$. Then $11\sim 10$, a $C_3$ appears, a contradiction. Thus $d(4)=4$ and the edge $(0, 4)$ is in a $C_4$ that must be $C_4:=4-0-2-6-4$. Consider the edge $(0, 4)$, let $z, w$ represent the third and fourth neighbor of vertex $4$, then we need $d(z, 1)=2$, which implies $4\sim 5$, then $d(4, 5)=2$, a contradiction. Thus the edge $(0, 3), (0, 4)$ must be in $C_4$. 

Note $d(3), d(4)$ cannot be both $2$.  Wlog, let $d(3)=4$.
Now we assume neither the $C_4$ for edge $(0, 3)$ nor the $C_4$ for the edge $(0, 4)$ pass through edge $(0, 2)$. Let $3\sim 10 11, 12$ with  $C_4:=0-3-10-4-0$.  Consider the edge $(0, 3)$, we need at least one of $d(1, 11), d(1, 12)$ is $2$, however, both would result $d(3, 5)=2$, a contradiction.
Thus one of $C_4$s for edges  $(0, 3), (0, 4)$  must pass through edge $(0, 2)$. Wlog, let $3\sim 6$. Then we need $d(4, 7)=3$ for the edge $(0, 2)$. In the following, there are two cases for the $C_4$ that passes through edge $(0, 4)$.
 \begin{itemize}
\item 
If the edges  $(0, 4)$  also  passes through edge $(0, 2)$, then it must be $4\sim 6$. Then $d(3, 7)=3$. 
In this case,  $d(3)=d(4)=4$ considering the edge $(0, 3)$. Let $3\sim 10, 11$, since $d(3, 5)=d(4, 5)=3$, then $d(1, 10)=d(1, 11)=3$ which implies $d(4, 10)=1$. 
 \begin{center}
 \includegraphics[scale=0.5]{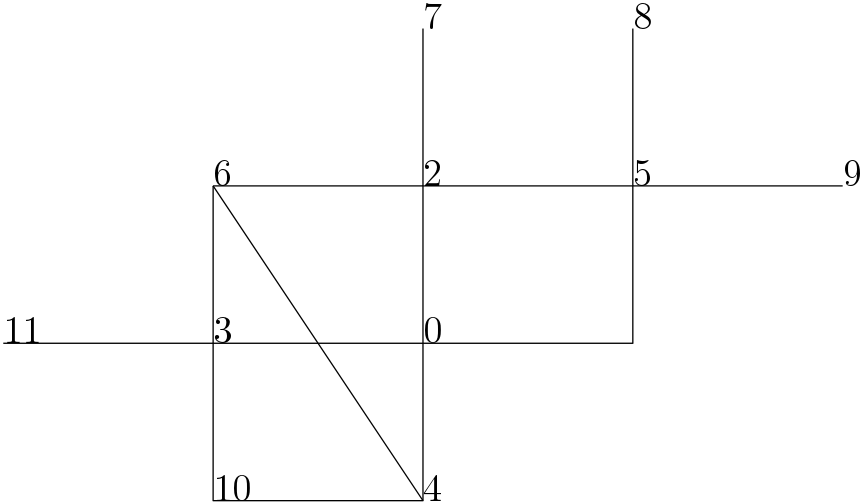}
 \end{center}

             Consider vertex $4$, note either $4\sim 11$ or  $4\sim 12$ where $12$ is a new vertex. Assume the second case, then consider the edge $(2, 6)$, let $6\sim t$, then we need $d(5, t)=1$ or $d(7, t)=1$. For $d(5, t)=1$, wlog, let $6\sim 8$, still consider the edge $(5, 9)$ which must be in a $C_4$ that pass through vertex $8$, let $9\sim w\sim 8$, where $t$ could be $11$ or $12$ or a new vertex.  Then consider the edge $(2, 7)$ which must be in a $C_4$ that passes through vertex $8$, let $7\sim 8$. Observe that we have $d(7, 9)=d(7, 3)=d(7, 4)=3$, then vertex $7$ must be adjacent to new vertices which have distances $3$ from vertices $0$ and $5$, then it is not good for the edge $(2, 7)$.

             Thus $4\sim 11$.  Note  one of the edge $(5,  8), (5, 9)$ must be in a $C_4$ that pass through the edge $(2, 5)$.  Assume $6\sim 8$, then consider the edge $(5, 9)$ which must be in a $C_4$ that pass through vertex $8$, since $d(6)=4$ in the current graph, we need a new vertex as common for $8, 9$. Observe that $d(8, 10)=d(8, 11)=3$ for the edge $(3, 6)$, then we need a new vertex $12$ such that $9\sim 12\sim 8$. Then consider the edge $(6, 8)$, we have both $d(3, 12), d(4, 12)$ are not $1$, then at least one of them must be $2$ which need $12\sim 10$ or $12\sim 11$, however, we  would have $d(8, 10)=2$ or $d(8, 11)=2$, a contradiction. 

             Thus when  $4\sim 11$, we need $7\sim 8$. Let a new vertex $12$ as the fourth neighbor of $6$. Consider the edge $(2, 6)$, since $d(4, 5)=d(3, 5)=3$,  we need $d(7, 12)=1$. 
              Consider the edge $(5, 8)$, $d(8)\neq 2$, otherwise the edge $(5, 9)$ cannot be in any $C_4$. Let $d(8)=4$, consider the edge $(5, 8)$, since vertex $1$ has distance $3$ to any new neighbor of vertex $8$, then we need vertex $9$ to have distance $1$ from neighbor of vertex $8$, that is, a common vertex for $8, 9$. Note $d(10, 12)=d(11, 12)=3$ for the edge $(3, 6)$. Then $d(12)\neq 4$. Otherwise the new neighbor of vertex $12$ must have distance $2$ from $3$ or $4$, then it must be adjacent to vertex $10$ or $11$ which would make $d(10, 12)=2$ or $d(1, 12)=2$. 
              Consider the vertices $10, 11$, we have $d(10)=4$ and the edge $(3, 10)$ satisfies ``Type 6b", such that the second $C_4:=3-10-w-11-3$ where $w$ is a new neighbor of vertex $10$. A similar situation for the vertex $11$ and the edge $(3, 11)$. Thus there is at least new common vertex for vertices $10$ and $11$. Let $z, t$ be the fourth neighbors of vertex $10$ and $11$ respectively.  Now the $C_4$ for edge $(z, 10)$ must pass through edge $(10, c)$ and the $C_4$ for edge $(11, t)$ must pass through edge $(11, c)$. Then either $z=t$ or they are two different vertices.  
%
              
               Consider the vertex $7$, there are two cases: either $7\not\sim 9$ or $7\sim 9$. Next, we consider the different situations under each case:
              
              \textbf{Cases for the fourth neighbor of vertex $7$. }
              \begin{itemize}

              \item Assume $7\sim 9$. Then $d(9)=4$. By above analysis, there are four different situations:
               \begin{figure}[H]
              \centering
              \includegraphics[scale=0.3]{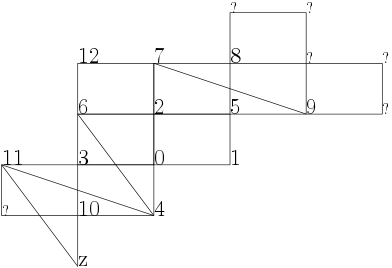}
              \hfil
              \includegraphics[scale=0.3]{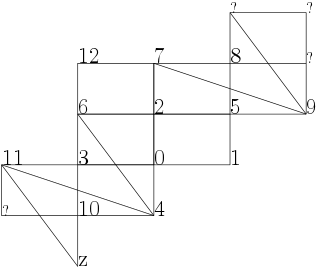}
           
              \includegraphics[scale=0.3]{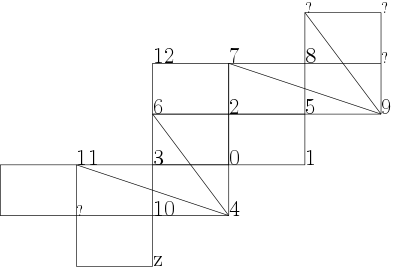}
              \hfil
              \includegraphics[scale=0.3]{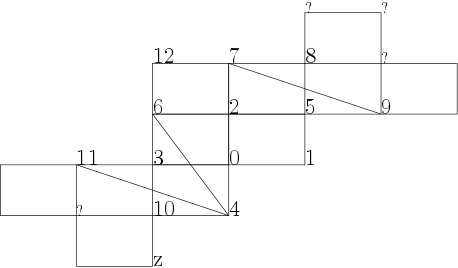}
                    \caption{Structures of the Ricci-flat graphs when $7\sim 9$}
              \end{figure}

              \item Assume $7\not\sim 9$. Let $7\sim 13$. Then we have $d(13)=2$ considering the edge $(7, 13)$, and $d(9)=2$ considering the edge $(5, 9)$. Still we need a common for vertices $8, 13$ and a common for vertices $8, 9$. A similar analysis for the vertices $10$ and $11$, the difference lies in whether the vertices $10, 11$ share two more common neighbors or just one more common neighbor, thus we have two different infinite structures. See Figure \ref{24Ricci-flat4}:
               \begin{figure}[H]
                 \includegraphics[scale=0.3]{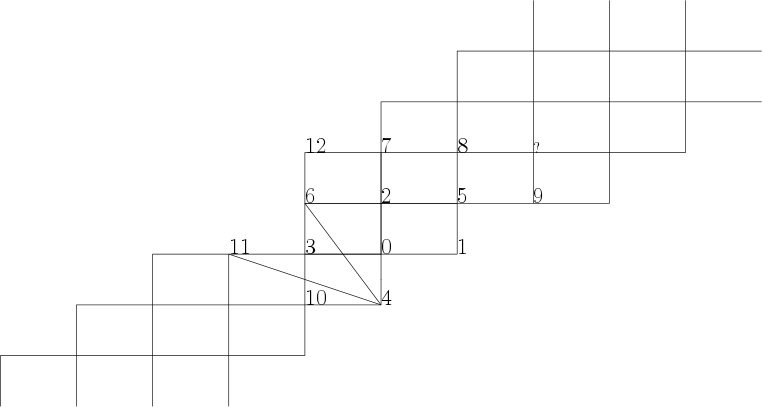}
                 \hfil
                  \includegraphics[scale=0.3]{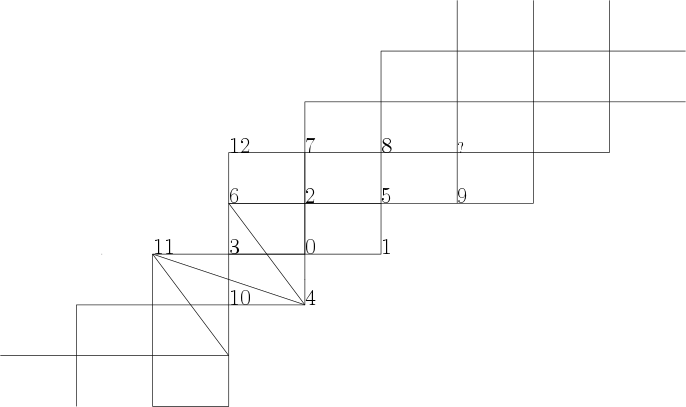}
                   \caption{Structure of Ricci-flat graph(when $7\not\sim 9$)}
                   \label{24Ricci-flat4}
              \end{figure}

%
%
              
             \end{itemize}
As we can see, we can also combine above different situations ($7\sim 9$ and $7\not\sim 9$) into one graph. The construction is: we can add the ``slash edges"(except $(4, 6), (4, 11)$) at any stage. To get a finite one, always merge the vertices on the left-upper corner region and vertices on the right-down conner region.  See the following illustration:
              \begin{figure}[h!]
              \centering
              \includegraphics[scale=0.4]{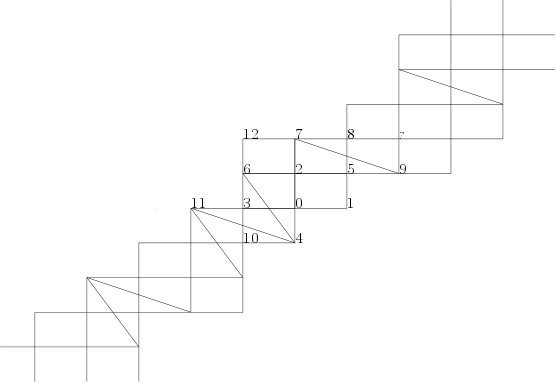}
              \caption{Ricci-flat graphs: the primitive graph}
              \label{24Ricci-flata}
              \end{figure}
              
For examples:
\begin{figure}[H]
 \centering
              \includegraphics[scale=0.4]{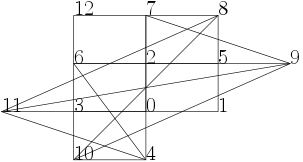}
              \hfil
              \includegraphics[scale=0.4]{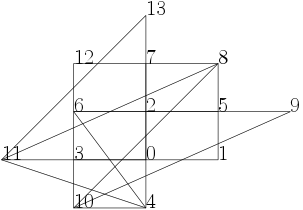}
\caption{Examples of finite Ricci-flat graphs}
\end{figure}

                 \begin{figure}[H]
                 \centering
                 \includegraphics[scale=0.4]{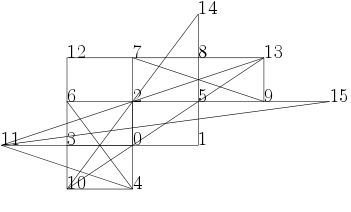}
                 \hfil
                 \includegraphics[scale=0.4]{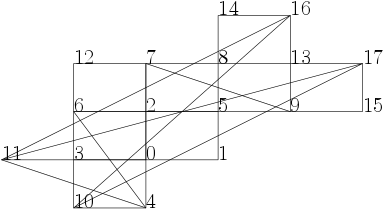}
                    \caption{Examples of finite Ricci-flat graphs}
                 \end{figure}

 \item      We consider the other case for the $C_4$ that passes through edge $(0, 4)$, that is, this $C_4$ passes through edge $(0, 3)$ but not vertex $6$. Let $3\sim 10\sim 4$ and $3\sim 11$.  Observe that $d(4)=2$ considering the edge $(0, 4)$. By symmetry, we consider the same situation for edges $(5, 8), (5, 9)$, wlog, let $d(8)=2, d(9)=4$. 
Then the edge $(5, 8)$ is in a $C_4$ that passes through $9$.  Consider the edge $(5, 9)$ , since any new neighbor of vertex $9$ has distance $3$ from vertex $1$, we need one of them to be adjacent to vertex $2$ such that  $(5, 9)$ shares two $C_4$. Then either $9\sim 6$ or $9\sim 7$. Note we must have $9\sim 7$ considering the edge $(2, 7)$.  Then $d(6, 8)=3$ for the edge $(2, 5)$. Now consider the edge $(2, 6)$ according to the following cases:
\begin{itemize}
\item Assume $d(6)=2$. Then $d(7)=2$, $d(11)=2$ considering the edge $(2, 7), (3, 11)$ respectively. Similarly, the fourth neighbor of vertex $9$ has degree $2$, and the fourth neighbor of vertex $10$ has degree $2$ as $d(4)=2$. 
 \begin{center}
\includegraphics[scale=0.3]{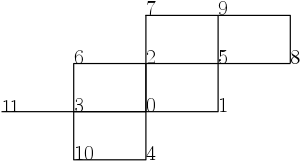}
\end{center}

Continue with the similar process, we will get the following structure which is an infinite Ricci-flat graph: 

\begin{figure}[H] 
\centering
\includegraphics[scale=0.4]{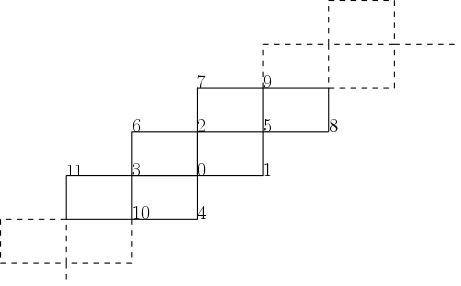}
\caption{The primitive Ricci-flat graph}
\label{24Ricci-flat5}
\end{figure}

To get a finite one we need to merge the vertices on the ends. For example, let vertex $10$ be the common neighbor of vertices $8$ and $9$. 
Then we need $9\sim 11$ for the edge $(3, 10)$ as both $d(6, 8), d(6, 9)$ are $3$. We obtain a Ricci-flat graph:
\begin{figure}[H]
\centering
\includegraphics[scale=0.4]{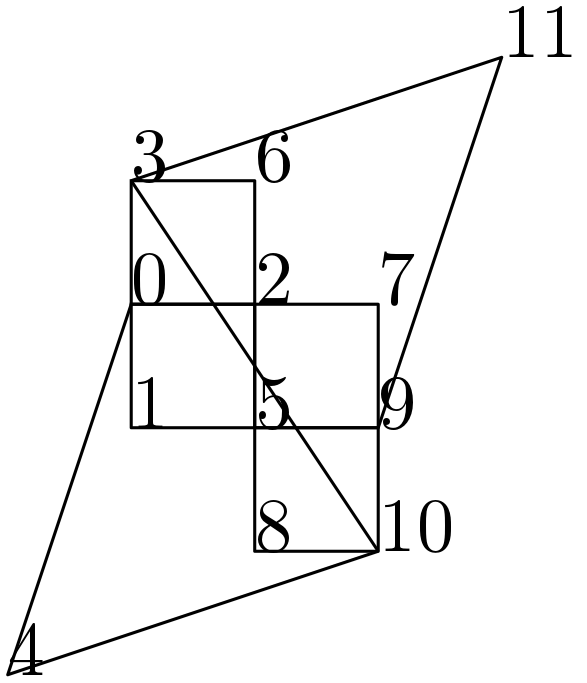}
\caption{Ricci-flat graph}
\end{figure}
Assume we need a new vertex $12$ as the common of $8, 9$, then $d(12)=4$, we need the edge $(9, 12)$ to share two $C_4$s, note $9\not\sim 11$ as $d(11)=2$. Let $9\sim 13$, thus we need a common vertex for $12, 13$. Note if $10\sim 13$, then $10\sim 12$, consider the edge $(3, 10)$, we have both $d(6, 12), d(6, 13)$ are $3$, then we need $11\sim 12$. 
\begin{figure}[H]
\centering
\includegraphics[scale=0.4]{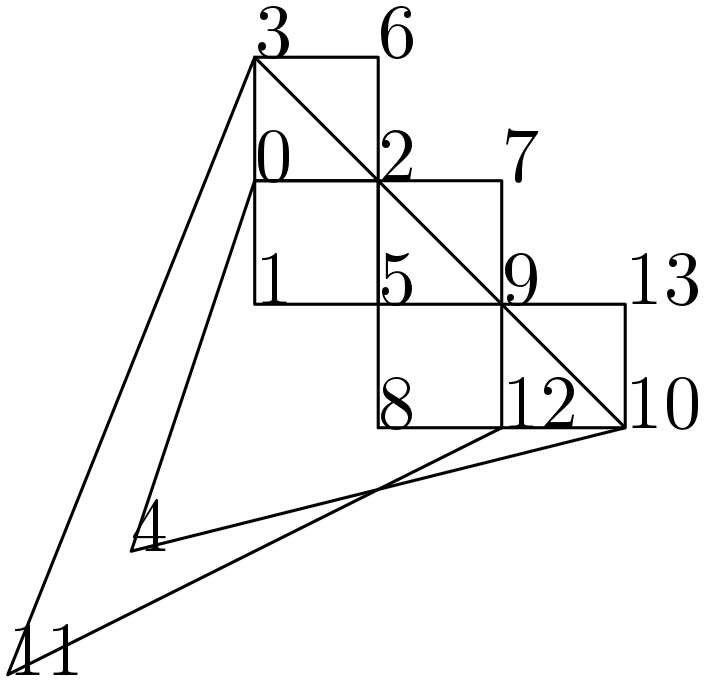}
\caption{Ricci-flat graph}
\end{figure}
Assume $10$ is not the common vertex for $12, 13$. Then let $12\sim 14\sim 13$, then $d(14)=4$, for the fourth neighbor of vertex $12$, assume $12\sim 11$, then we have $12\sim 10$, or assume $12\sim 10$, then $12\sim 11$ which would make $d(12)=5$. Thus let $12\sim 15$, we need a  common vertex for $14, 15$. Still assume $10\sim 14, 14$. Then we need $14\sim 11$ for the $C_4$ that passes through edge $(3, 10)$. 
\begin{figure}[H]
\centering
\includegraphics[scale=0.4]{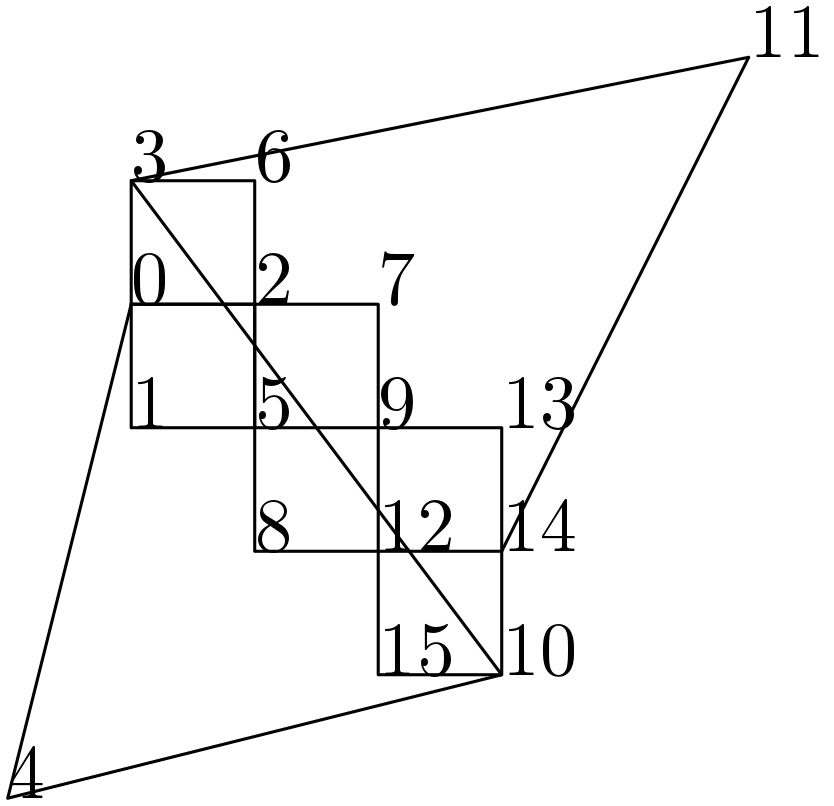}
\caption{Ricci-flat graph}
\end{figure}
Continue with the similar process, we will obtain the following Ricci-flat graphs. 
\begin{figure}[H]
\centering
\includegraphics[scale=0.4]{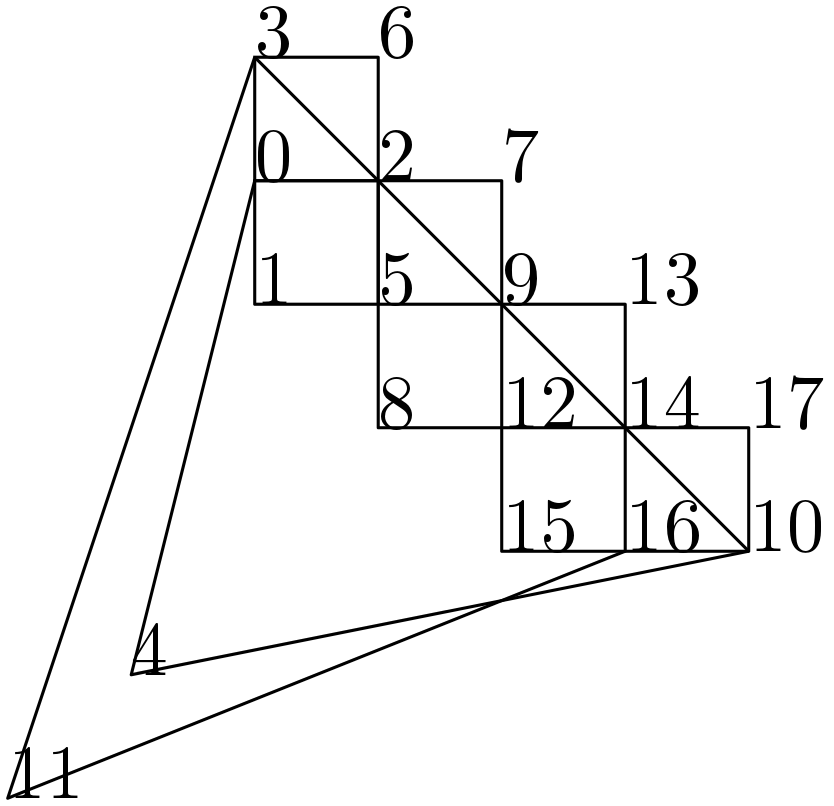}
\hfil
\includegraphics[scale=0.4]{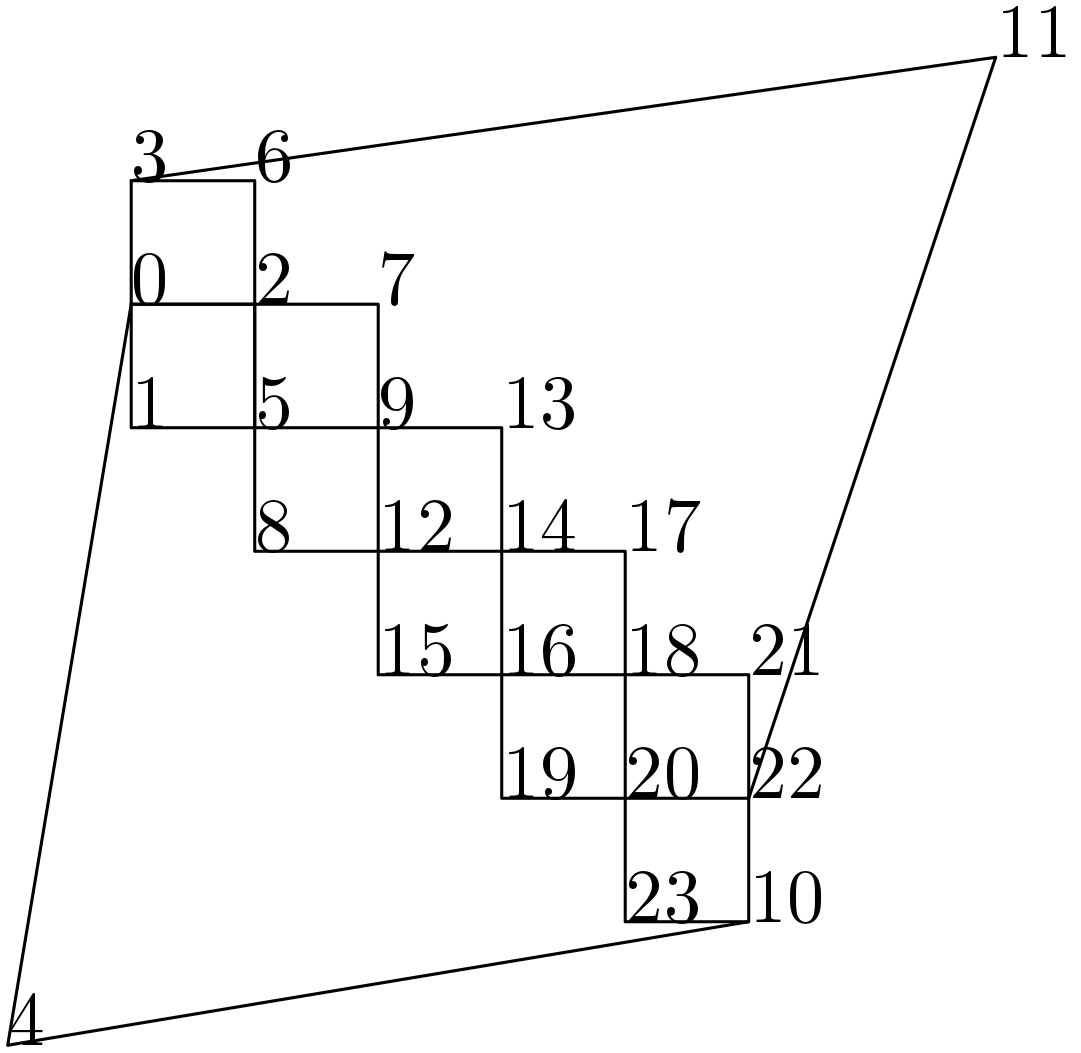}
\hfil
\includegraphics[scale=0.4]{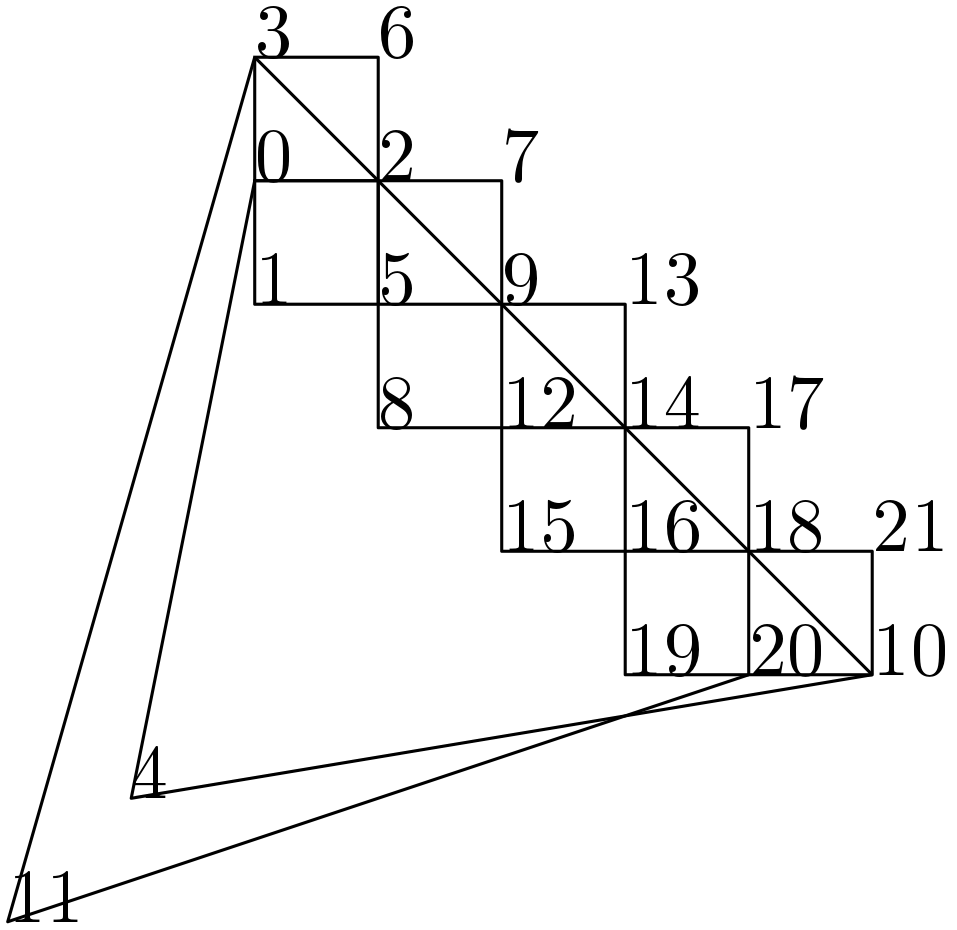}
\caption{Ricci-flat graphs}
\end{figure}

The idea behind these Ricci-flat graphs is simple, we can first have the infinite Ricci-flat graph(see Figure \ref{24Ricci-flat5}),
then  two ways to obtain a finite Ricci-flat:

\begin{figure}[H] 
\centering
\includegraphics[scale=0.4]{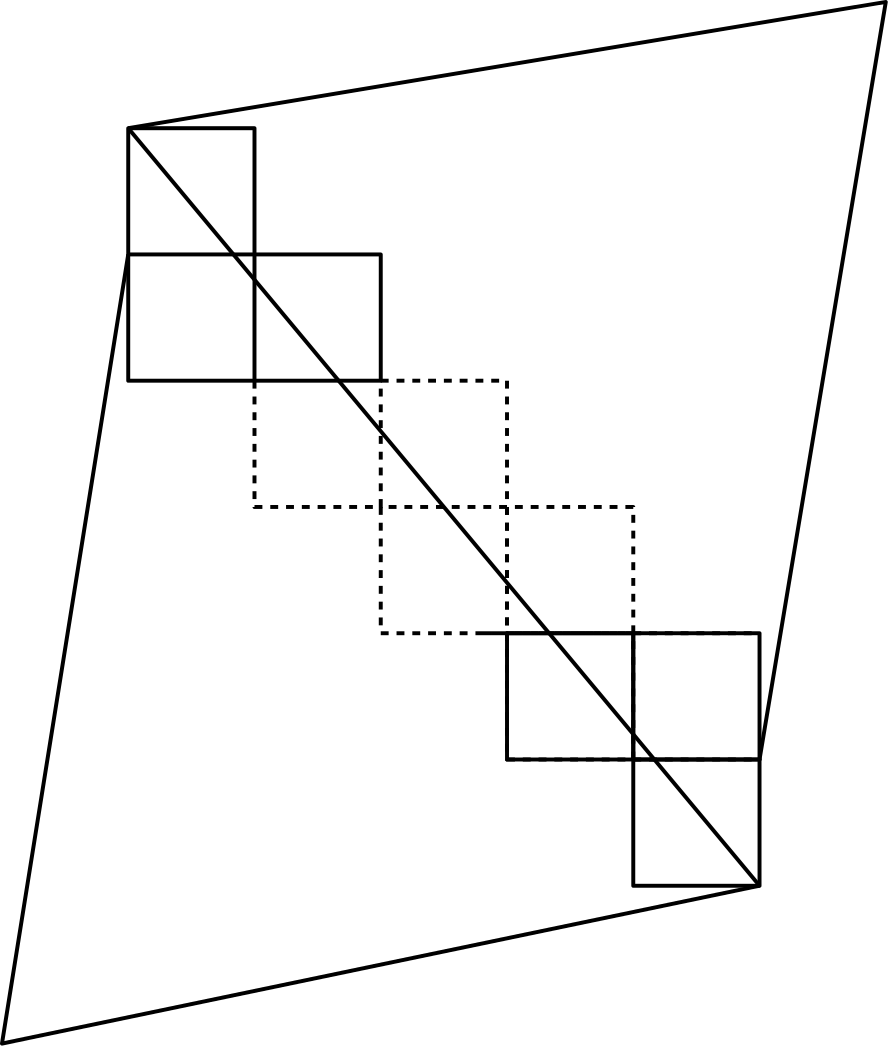}
\hfil
\includegraphics[scale=0.4]{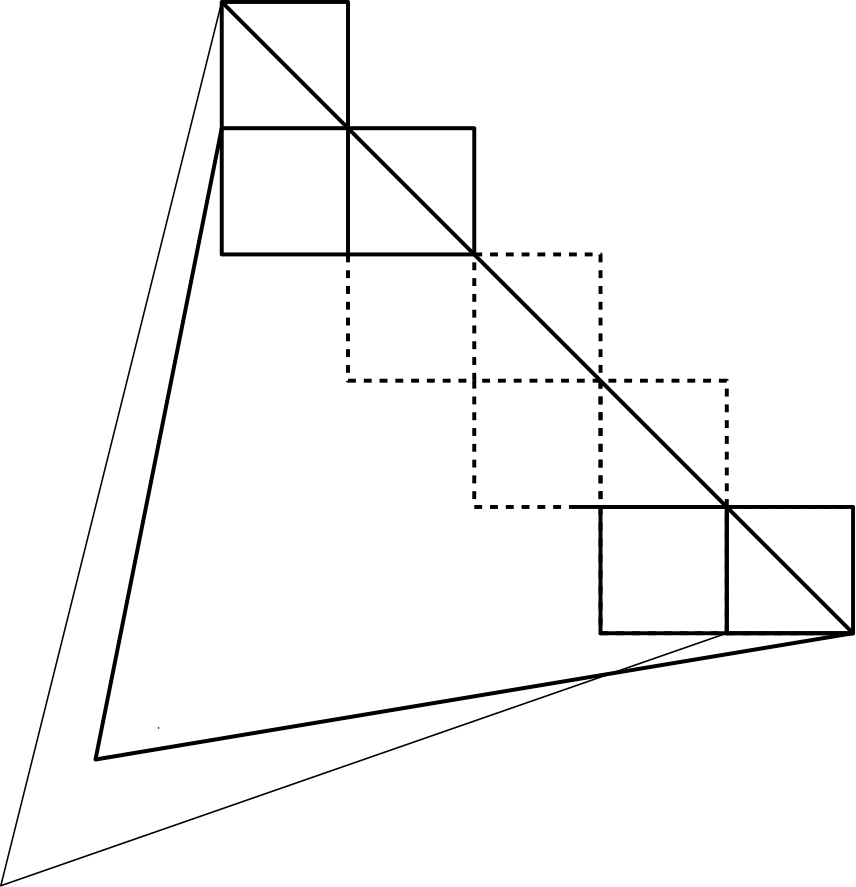}
\caption{Two ways to generate finite types}
\end{figure}

\item Assume $d(6)=4$. Then $d(7)=4$.  For the edge $(2, 6)$ in a $C_4$, there are two cases: either  $6\sim 12, 13$ and $7\sim 12$, or $6\sim 9, 12$. Consider the latter case:  $6\sim 9, 12$, then we need $d(7, 12)=d(5, 12)=3$ for the edge $(2, 6)$. For the edge $(2, 7)$, note $7\not\sim 11$, as both $d(5, w), d(6, w)$ cannot be $2$ where $w$ is the fourth neighbor of vertex $7$. Thus $7$ is adjacent to two new vertices, let $7\sim 13, 14$, since
the  vertex $0$ has distance $3$ from both vertices $13, 14$, we need $5\sim 13$ or $6\sim 14$, a contradiction. Thus $6\not\sim 9$. 

Let $6\sim 12, 13$ with two new vertices. Note we still need the ``Type 6b"  for the edge $(2, 6)$, otherwise suppose $d(5, 13)=2$, then either $9\sim 12$ or $8\sim 12$.  The latter case implies the former case, since $d(8)=2$ and we need a common vertex for $8, 9$. While the former case is not good for the edge $(5, 9)$ as the vertex $1$ always has distance $3$ from the fourth neighbor of vertex $9$, thus we cannot have ``Type 6c".  Using ``Type 6b",  we need $12\sim 7$. 
\begin{figure}[H]
\centering 
\includegraphics[scale=0.4]{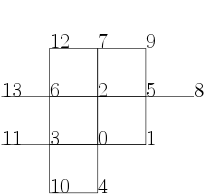}
\end{figure}
\begin{itemize}
\item Assume $d(12)=2$. Then both vertex $13$ and the fourth neighbor of vertex $7$ have degree $2$. For the edge $(6, 13)$, we need $13\sim 10$ or $13\sim 11$. Note if $3\sim 10$, then we need $d(11, 12)=3$ for the edge $(3, 6)$. Then the $C_4$ for the edge $(3, 11)$ must pass through edge $3, 10$, then we need $d(13, 6)=2$ for the edge $(3,10)$, contradict to the fact $d(13, 6)=1$. Let $13\sim 11$. Continue with this process, we will get the following infinite Ricci-flat graph:

\begin{figure}[H]
\centering 
\includegraphics[scale=0.42]{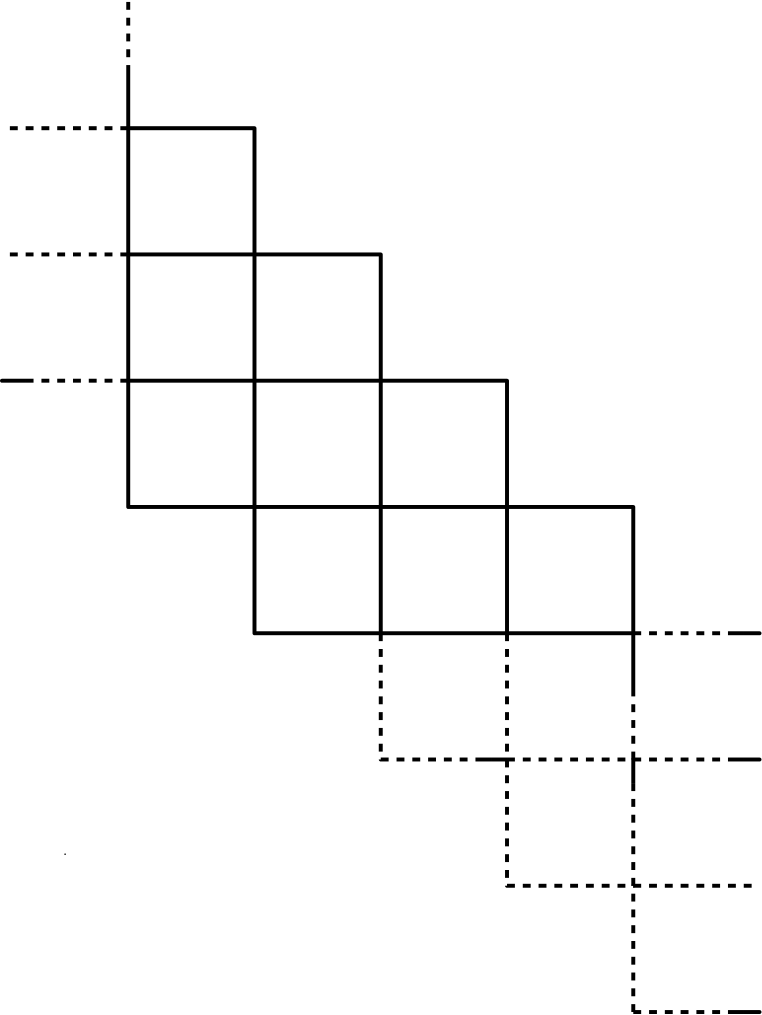}
\caption{The primitive Ricci-flat graph}
 \label{24Ricci-flat7}
\end{figure}

To get a finite one,  need to merge the vertices on the ends.
Now consider the common vertex for $8, 9$, from the current vertices, it can be vertex $10$ or $11$. Assume $11\sim 8, 9$, then we need $9\sim 10$ however, we cannot guarantee the edge $(9, 10)$. Assume $10\sim 8, 9$, then we need $9\sim 11$ for the edge $(9, 10)$. Then for the edge $(3, 11)$, we need  $11\sim 13$, then we need $7\sim 14 \sim 11$ for the edge $(9, 11)$.  Another Ricci-flat graph:
\begin{figure}[H]
\centering 
\includegraphics[scale=0.4]{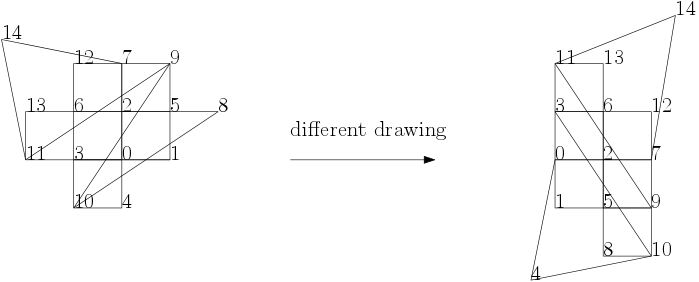}
\caption{Ricci-flat graph}
\end{figure}

If we need a new vertex $14$ as common for vertex $8, 9$.  We can have the following Ricci-flat graphs by similarly arguments:
\begin{figure}[H]
\centering 
\includegraphics[scale=0.4]{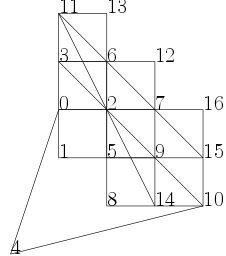}
\hfil
\includegraphics[scale=0.4]{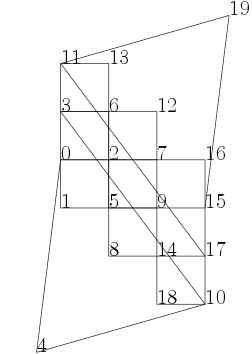}
\caption{Ricci-flat graphs}
\end{figure}

Now we can see, there are also two ways to generate finite Ricci-flat graphs:
\begin{figure}[H]
\centering 
\includegraphics[scale=0.4]{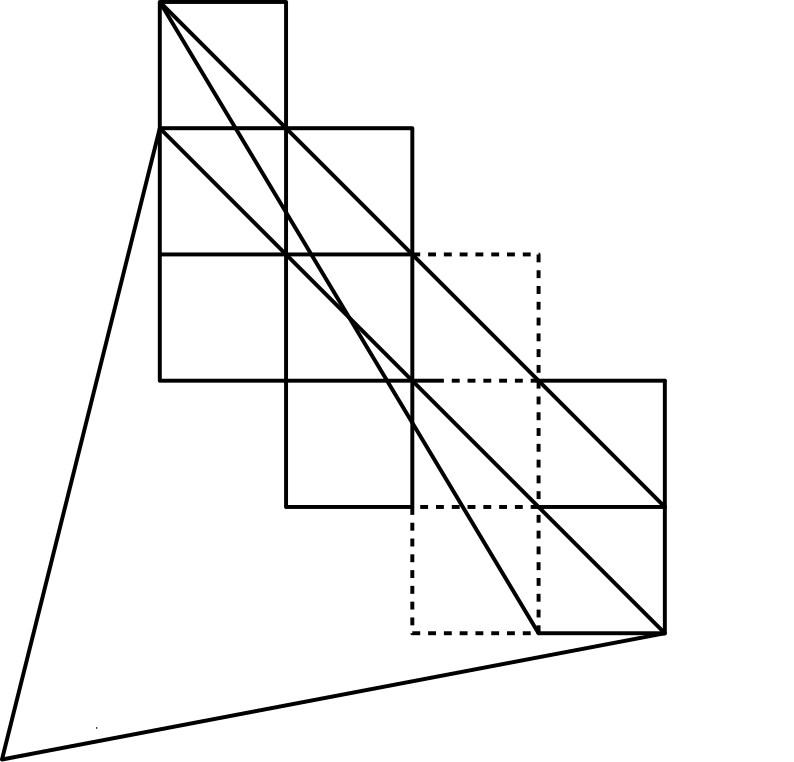}
\hfil
\includegraphics[scale=0.4]{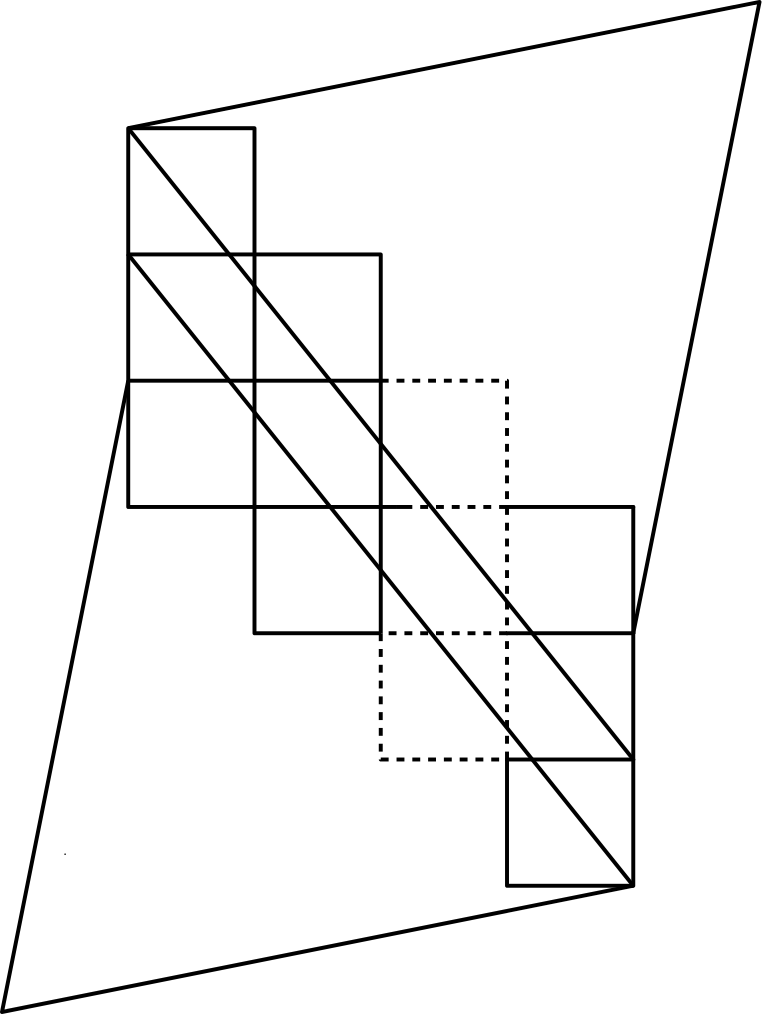}
\caption{The primitive Ricci-flat graph}
\end{figure}

\item Assume $d(12)=4$. Let $12\sim 14$. Then assume $d(14)=2$ or $4$. For each case, the similar arguments are used, here we omit the details. 
We will have the following infinite and finite Ricci-flat graphs: 
\begin{figure}[H]
\centering 
\includegraphics[scale=0.6]{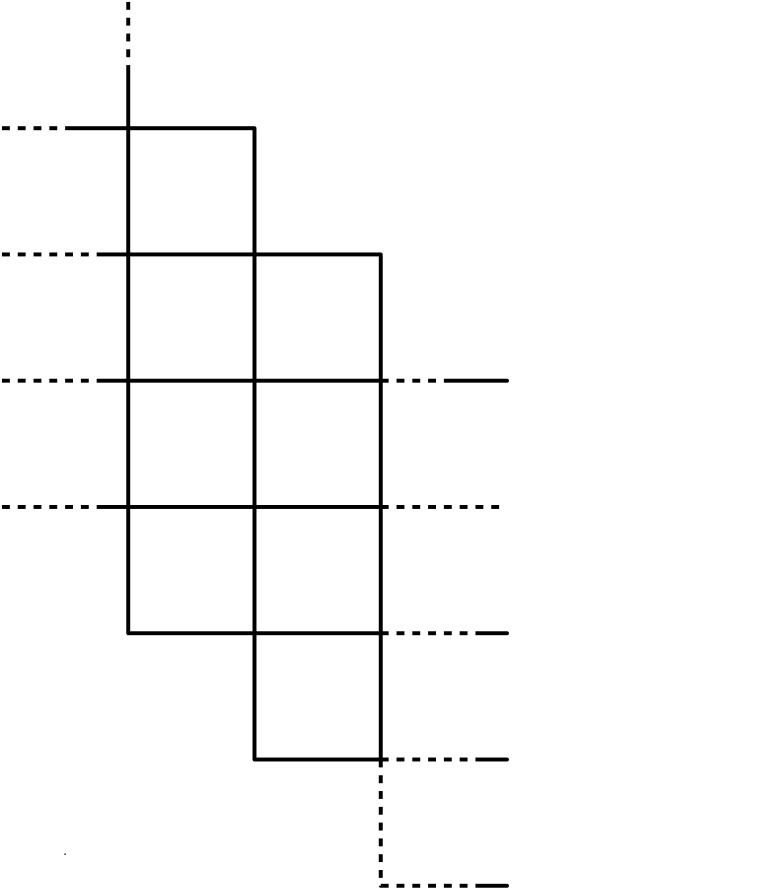}
\hfil
\includegraphics[scale=0.4]{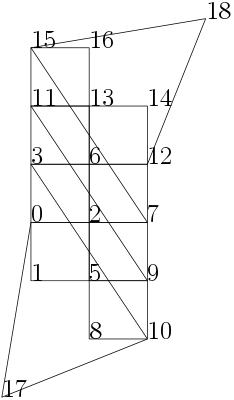}
\hfil
\includegraphics[scale=0.4]{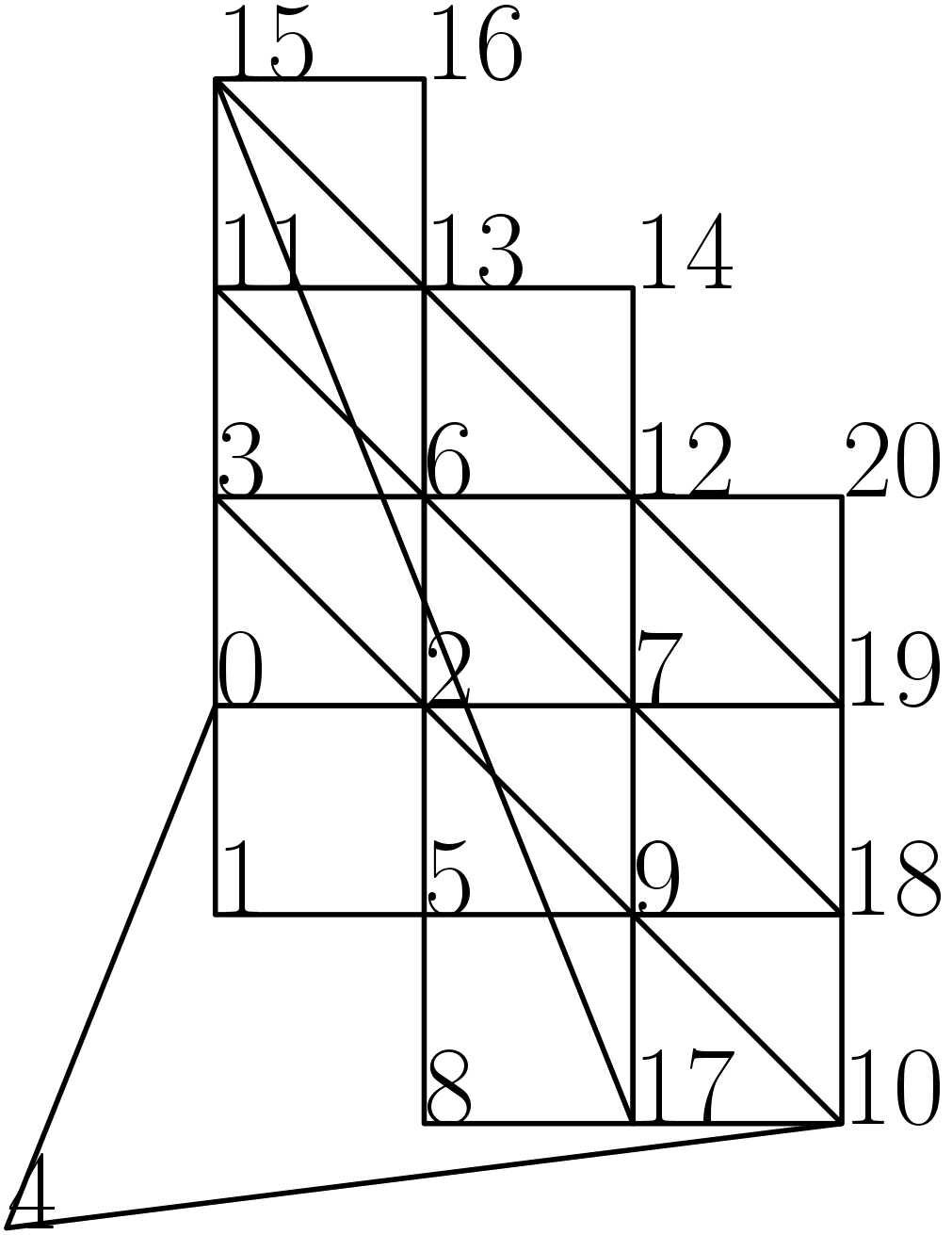}
\caption{Ricci-flat graphs}

\end{figure}

\end{itemize}

Thus there are infinitely many finite Ricci-flat graphs. To describe them,  we first specify the drawing as showing above, define the number of $C_4$ in each column as the length of the graph, define the width as the number of columns in the drawing.  For example above finite graphs have length  $4$ (if $d(14)=4$, we could have a larger length),  width $2$ and width $3$ respectively.   We use an example with length $5$ to illustrate the infinite Ricci-flat graphs and the first way to generate finite types:
\begin{figure}[H]
\centering
\includegraphics[scale=0.35]{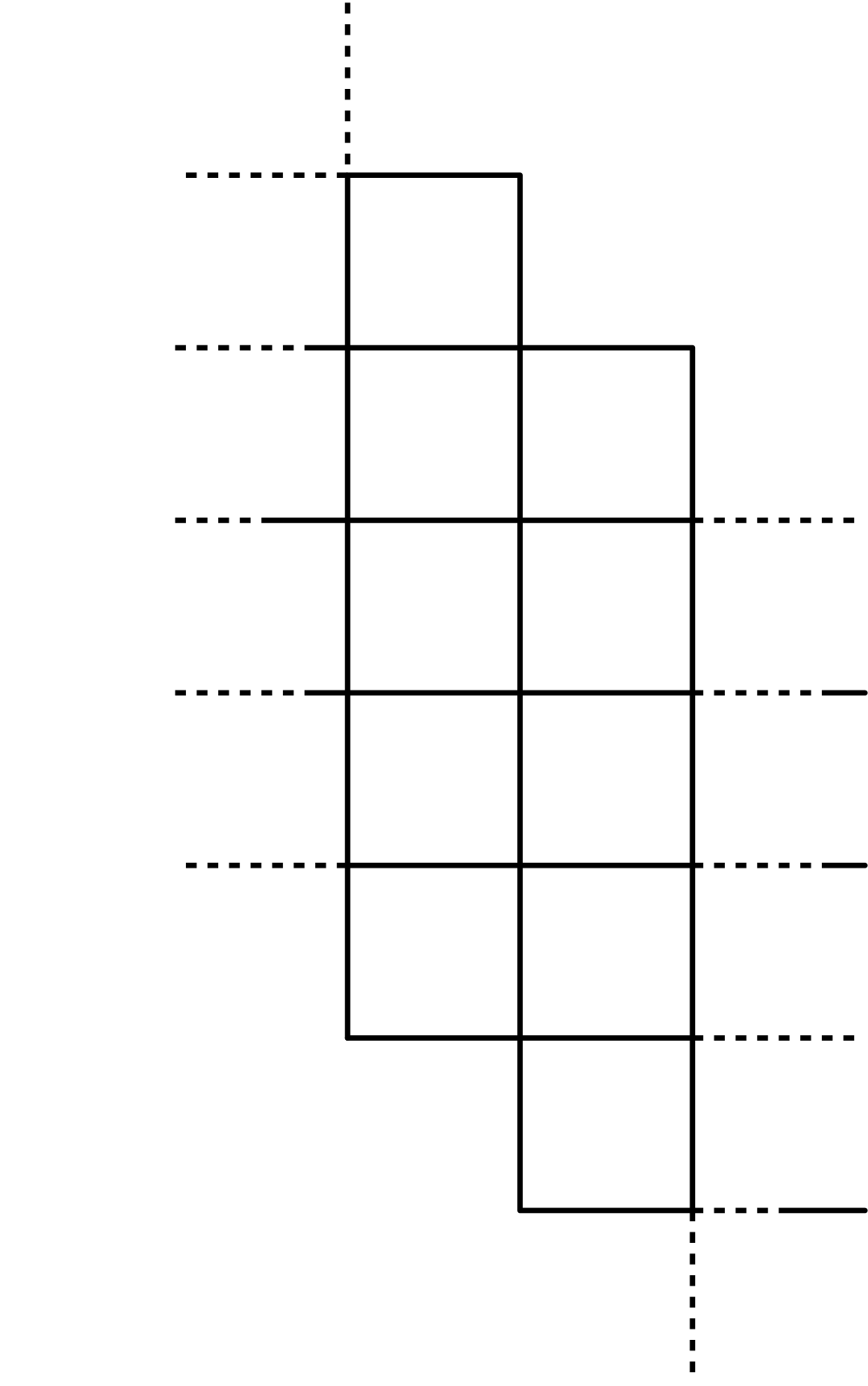}
\hfil
\includegraphics[scale=0.4]{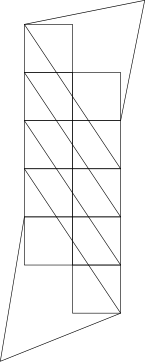}
\caption{Infinite Ricci-flat graph of length $5$,  the first way to generate a finite type with width $2$.}
 \label{24Ricci-flat8}
The second way to generate finite Ricci-flat graph starts with width at least $3$, see the following example. 
\end{figure}
\begin{figure}[H]
\centering
\includegraphics[scale=0.4]{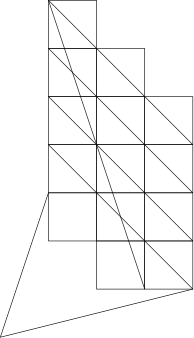}
\caption{The  second way to generate a finite type with width $3$.}
 \label{24Ricci-flat9}
\end{figure}

\end{itemize}

\end{itemize}
\end{itemize}
\end{proof}

The other situation for the second case is that  all four vertices in any $C_4$ have degree $4$, then $G$ must be $4$-regular, we have the following results. 
\begin{theorem}\label{thm:C444edge}
Let $G$ be a $4$-regular Ricci-flat graph  that contains two $C_4$s sharing one edge. Then G is isomorphic to to graphs with the primitive graphs showing in the Figures \ref{44C4Ricci-flat1}, \ref{44C4Ricci-flat2}. 
\end{theorem}

\begin{proof}
Start with the following structure where the edge $(0, 1)$ shares two $C_4$s:
\begin{center}

\includegraphics[scale=0.4]{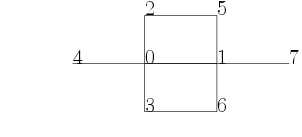}

\end{center}
There are two cases for the edge $(0, 2)$:

Case 1:  $(0, 2)$ shares two $C_4$, let the second one be $C_4=0-2-8-4-0$. 
\begin{center}

\includegraphics[scale=0.4]{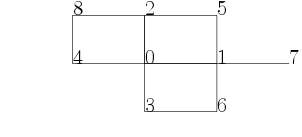}

\end{center}
Note for the edge $(0, 2)$, we need the fourth neighbor of $2$ has distance $3$ from vertex $3$, thus $2\not\sim 6$. 
 there are two cases, either $2\sim 7$ or $2\sim 9$ where $9$ is a new vertex.  
 \begin{itemize}

 \item 
 Assume $2\sim 7$.
Then $d(3, 7)=d(3, 5)=d(4, 7)=d(4, 5)=3$.  We will show that the edge $(1, 6)$ does not satisfy ``Type 6c".
We first assure the $C_4$ for edge $(1, 6)$ is shown in above subgraph. Suppose we also have $C_4:=1-6-0-4-1$. Then according to `Type 6c", we need at least one of  $d(3, 5), d(3, 7)$ to be $2$, a contradiction. Thus $6\not\sim 4$. Then vertex $6$ must be adjacent to at least one new vertex, let it be $6\sim 9$. Wlog, we need $d(5, 9)=2$ for the edge $(1, 6)$. Then we need a new vertex as common for vertices $5, 9$, let $5\sim 10\sim 9$. Consider the edge $(1, 5)$, since $d(4, 5)=d(3, 5)=3$, as the vertex $0$ has distance $3$ from any new neighbor of vertex $5$, then we need vertex $6$ to be adjacent to the fourth neighbor of vertex $5$, a contradiction. Thus the neighborhood for the edge $(1, 6)$ must be ``Type 6b".  
That is,  the edge $(1, 6)$ must share a $C_4$ that passes through either vertex $5$ or vertex $7$. Assume it passes through vertex $5$. Then we need a new vertex $10$ as the common for $5, 6$. Let $5\sim 10$, we need $d(0, 10)=d(7, 10)=3$ for the edge $(1, 5)$. Then let $7\sim 11, 12$, we have the $C_4$s for the edge $(7, 11)$ or edge $(7, 12)$ must pass through edge $(7, 2)$ then $(2, 8)$, which implies $8\sim 11, 12$. However, for the edge $(2, 5)$, since both vertices $0$ and $7$ have distance $3$ from the third and fourth neighbors of vertex $5$, then we need vertex $8$ to have distance $1$ to them, that is,  $5\sim 11$ or $5\sim 12$, however, both are not good for the edge $(2, 7)$. A contradiction. Thus  let $6\sim 9\sim 7$. 
\begin{center}

\includegraphics[scale=0.4]{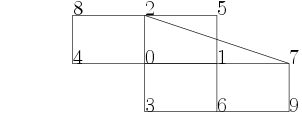}

\end{center}
Let $7\sim 10$. We need $d(5, 10)=3$ and $d(0, 10)=3$ for the edge $(1, 7)$. Consider the edge $(2, 7)$ which is in the $C_4:=2-7-1-5-2$ and $2-7-1-0-2$, since $d(5, 9)=d(5, 10)=d(0, 9)=d(0, 10)=3$, we need either $8\sim 9$ or $8\sim 10$ or both. For any of these three cases, we need vertex $5$ to have distance $3$ to new neighbors of vertex $8$. Let $s$ represent the new neighbor of vertex $5$, we have $d(8, s)=3$. Then the situation is not good for the edge $(2, 5)$. A contradiction again. Thus $2\not\sim 7$.

\item Assume $2$ is adjacent to a new vertex $9$. Then we need $d(3, 9)=d(4, 7)=3$. We will show that the edge $(0, 6)$ does not satisfy ``Type 6c".
Otherwise, consider the $C_4$ which would be used for ``Type 6c" for edge $(1, 6)$, there are three cases:
\begin{itemize}
\item Assume $C_4:=1-6-4-0-1$. Then $d(3, 7)=3$ for the edge $(0, 1)$. Let $6\sim 10$, we need $d(3, 5)=d(7, 10)=2$. Then the $C_4$ the edge $(1, 7)$ must pass through edge $(1, 5)$. Then we need a new vertex as the common for vertices $5, 7$, let it be $5\sim 11\sim 7$. Note for $d(3, 5)=2$, $3\not\sim 11$ as $d(3, 7)=3$, let $3\sim 12\sim 5$, then we have $d(6, 12)=2$, however, this will contradict to the edge $(1, 5)$ which requires $d(6, 12)=3$. 
\item Assume $C_4:=1-6-3-5-1$.  Then we need the fourth neighbor of vertex $5$ to have distance $3$ from vertex $7$, then the $C_4$ the edge $(1, 7)$  cannot pass through edge $(1, 5)$, it must be pass through edge $(0, 1)$ which implies $3\sim 7$. Then consider the edge $(3, 5)$, which now shares two $C_4s:=3-5-1-7-3, 3-5-2-0-3$, then we need the fourth neighbor of vertex $5$ to have distance $3$ from vertex $6$. Since $d(3, 9)=2$, then $6\not\sim 8, 9$. Let $6\sim 10, 11$,  then we have $d(5, 10)=d(5, 11)=3$. Then we need $d(0, 10)=2$ for the edge $(1, 6)$ which implies $4\sim 10$, which is not good for the edge $(0, 1)$. A contradiction.

\item Assume $C_4:=1-6-3-7-1$. Then $d(4, 6)=3$ for the edge $(0, 1)$. Consider the edge $(0, 3)$, since $d(4, 7)=d(4, 6)=3$, then we either need $d(2, 6)=2$( or $d(2, 7)=2$ by symmetry of vertices $6$ and $7$) or $d(2, z)=1$ where $z$ is the fourth neighbor of vertex $3$. Consider  $d(2, 6)=2$, then either $6\sim 8$ or $6\sim 9$. However $6\sim 8$ would contradict $d(4, 6)=3$,  $6\sim 9$ would contradict $d(3, 9)=2$. Thus let $d(2, z)=1$. Note $3\sim 5$ would give us a same situation as previous item. Thus we need $3\sim 8$. Then $d(4, 9)=3$ for the edge $(0, 2)$. 
 Consider the fourth neighbor of vertex $8$, let $8\sim 10$, since $d(2, 6), d(2, 7), d(4, 6), d(4, 7)$ are all $3$, we need either $6\sim 10$ or $7\sim 10$, wlog, let $6\sim 10$. Then consider the edge $(3, 6)$ which shares two $C_4:=3-6-1-0-3, 3-6-10-8-3$, then we need $d(7, 11)=3$ where $11$ is the fourth neighbor of vertex $6$, which would contradict to ``Type 6c" for the edge $(1, 6)$. 
\end{itemize}
Since all possible cases lead to contradictions, we conclude that  the edge $(1, 6)$ must satisfy ``Type 6b", that is,  the edge $(1, 6)$ must share a $C_4$ that pass through either vertex $5$ or vertex $7$. Assume it pass through vertex $5$. Then we need a new vertex $10$ as the common for $5, 6$. Then we need vertex $7$ to have distance $3$ from the fourth neighbors of vertices $5$ and $6$ respectively, which would result that the edge $(1, 7)$ cannot be in any $C_4$, a contradiction. Thus the $C_4$ for the edge $(1, 6)$  must pass through vertex $7$, let $6\sim 10\sim 7$. Similarly,  let $6\sim 11$.  We obtain the following structure. 
\begin{figure}[H]
\centering
\includegraphics[scale=0.4]{C4C4.png}
\hfil
\includegraphics[scale=0.4]{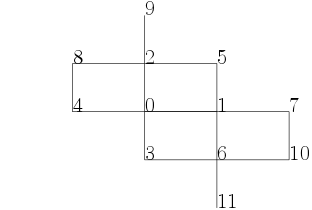}
\caption{The left structure generate the unique right structure.}
\end{figure}

Continue with the similar arguments, we will have the following structure: 
\begin{figure}[H]
\centering
\includegraphics[scale=0.4]{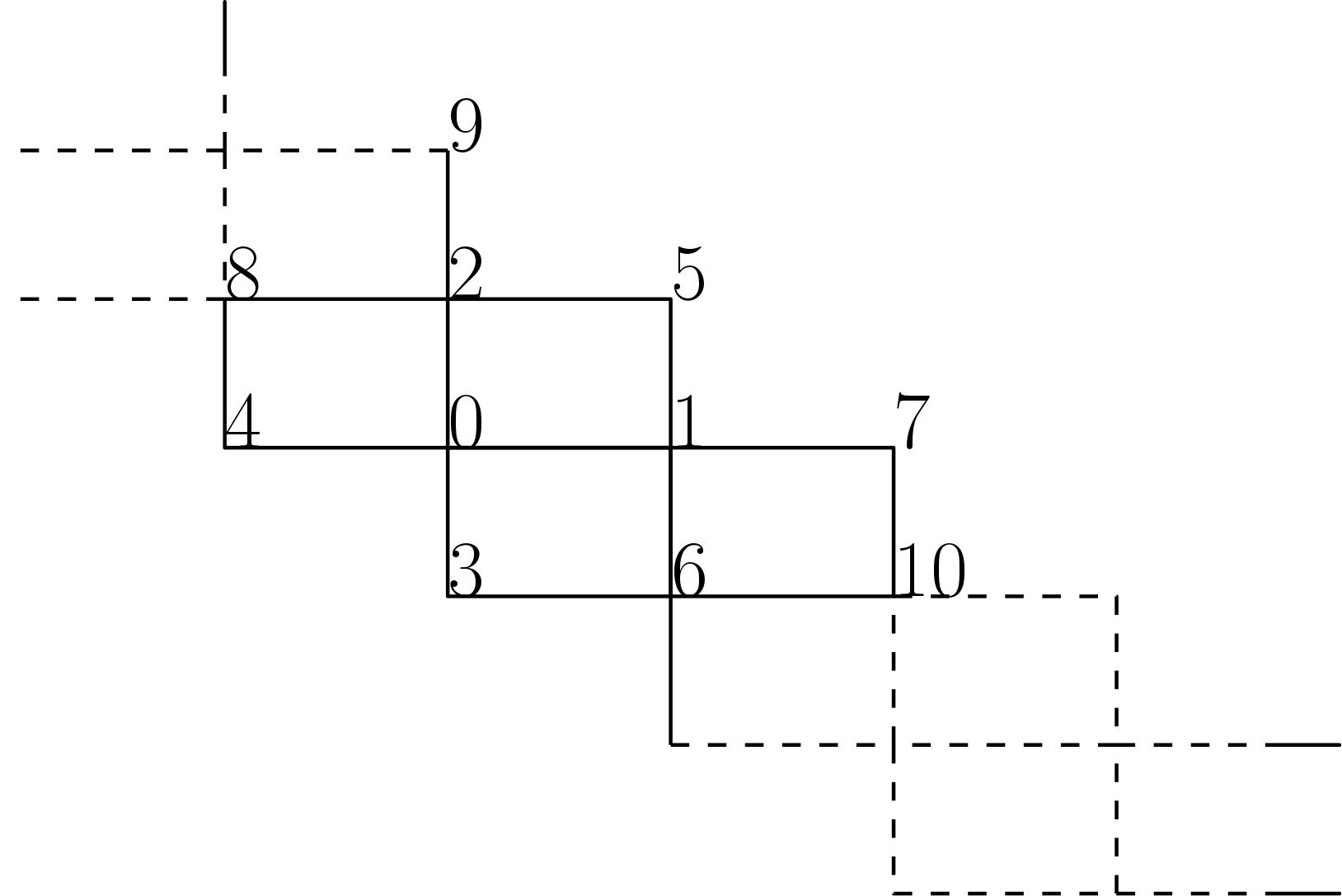}
\end{figure}

The above is a Ricci-flat primitive graph, as we require that all vertices have degree $4$, we consider two different situations according to the neighbors of vertex $4$. 
\begin{itemize}
\item Assume $4\sim 5$.  Consider the fourth neighbor of vertex $5$, note $5\not\sim 10$ as there is no way to generate a $C_4$ passing through edge $(5, 10)$ in above structure. The other possible case is $5\sim 3$ or $5\sim 12$, a new vertex. Consider $5\sim 12$, then $d(9, 12)=3$ for the edge $(2, 5)$. Note for the edge $(2, 8)$, which should share two $C_4$s, we need the fourth neighbor of vertex $8$ to have distance $3$ from vertex $5$, thus  the only possible way to form a $C_4$ passing through the edge $(5, 12)$ is through edge $(5, 1)$ then edge $(1, 7)$. Let $7\sim 12$. Then consider the edge $(5, 12)$, since we need vertex $12$ to have distance $3$ to the fourth neighbor of vertex $4$, then vertex $4$ has distance $3$ to the third and fourth neighbor of vertex $12$, then for the edge $(5, 12)$, we need $12\sim 8$. However, there is no way to generate a $C_4$ passing through vertices $2, 9, 8$. A contradiction. 
 Thus, when $4\sim 5$, we need $5\sim 3$. Similarly, $3\sim 7, 4\sim 9$. 
 Following the similar arguments, we have the following infinite Ricci-flat graph:
 \begin{figure}[H]
\centering
\includegraphics[scale=0.4]{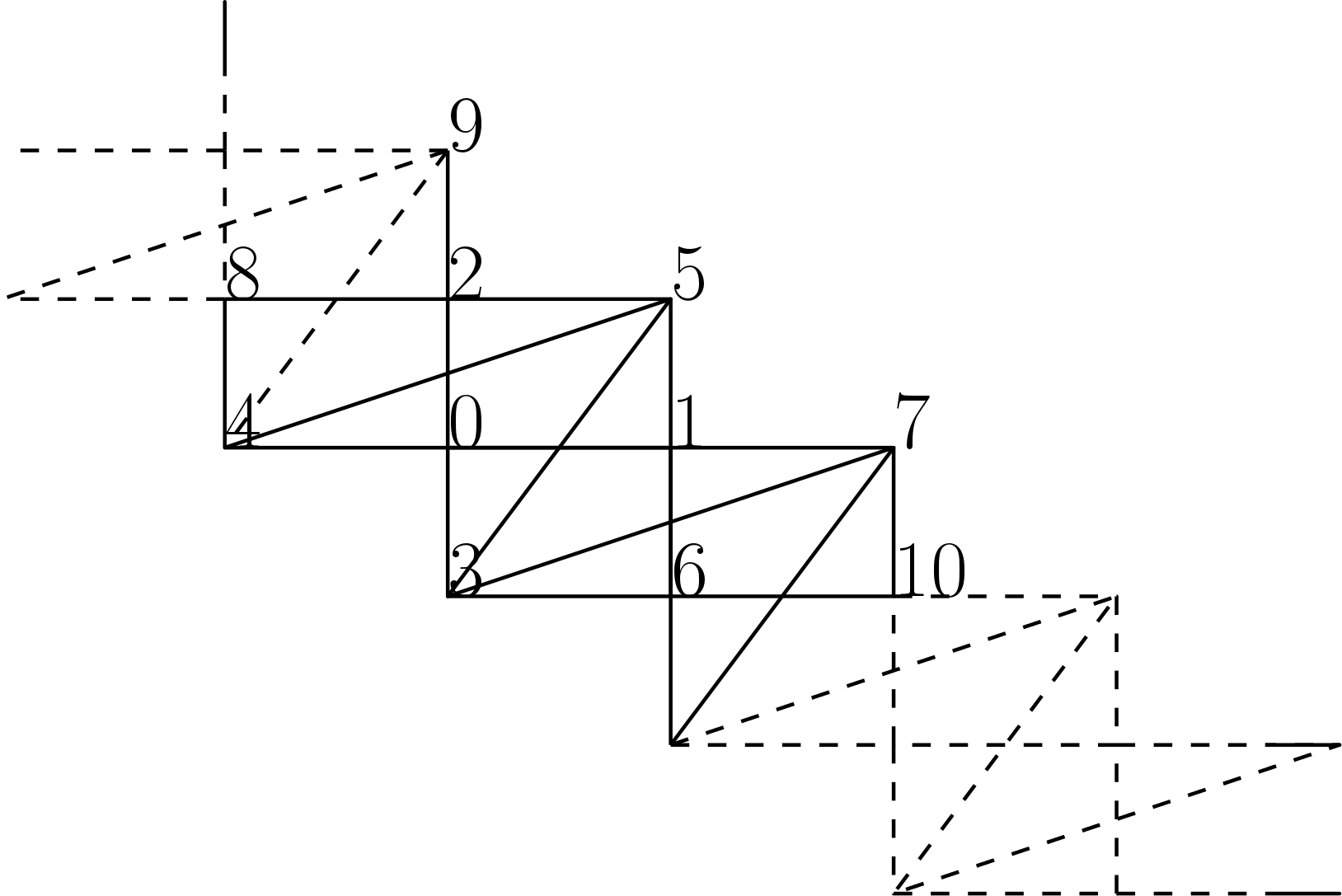}
\caption{A primitive  Ricci-flat graph}
\label{44C4Ricci-flat1}
\end{figure}

 To get a finite one, always merge the vertices on the left-upper corner region and vertices on the right-down corner region. For example, let $10, 11\sim 8, 9$ or new vertex as common for vertices $8$ and $9$. See the following finite Ricci-flat graphs based on above structure: 

\begin{figure}[H]
\centering
\includegraphics[scale=0.4]{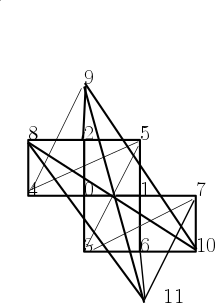}
\includegraphics[scale=0.4]{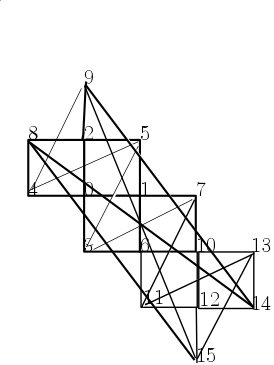}
\caption{Two  finite $4$-regular  Ricci-flat graphs}
\end{figure}

\item Now we consider the case when $4\not\sim 5$, by symmetry $3\not\sim 5$. Assume $3\sim 8$, then $4\sim 6$. Then $3\not\sim 7$ and $4\not\sim 9$. Let $3\sim 10$, note $4\not\sim 10$ for the edge $(0, 3)$. For the edge $(1, 7)$ to be in  a $C_4$, which must pass through edge $(1, 6)$, let $6\sim 12\sim 7$. Similarly, let $9\sim 13\sim 8$. Note $12\neq 13$, since then the edge $(3, 10)$ cannot be in any $C_4$.  Consider the edge $(0, 3)$, $(3, 6), (3, 8)$, we need $d(2, 10)=3, d(10, 12)=3$, $d(10, 13)=3$ respectively, then the edge $(3, 10)$ cannot be in any $C_4$. A contradiction. 

Now consider the case when $3\sim 10, 11$, then the edge $(3, 10)$  must be in a $C_4$ that pass through edge $(0, 4)$. Let $4\sim 10, 4\sim 12$. Continue with this process we can get an infinite lattice. 

\begin{figure}[H]
\centering
\includegraphics[scale=0.3]{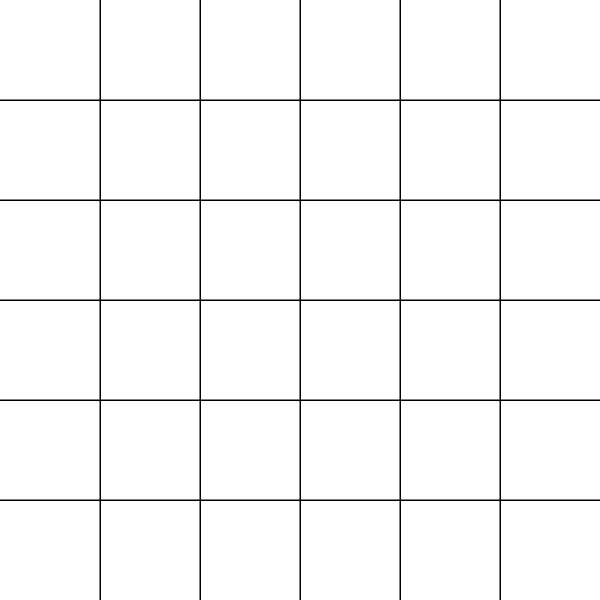}
\caption{A primitive $4$-regular Ricci-flat graph of ``lattice type"}
\label{44C4Ricci-flat2}
\end{figure}

We define the length as the number of $C_4$s in each row and  the length as the number of $C_4$s in each column, 
then width and length of the lattice could be extended. Note the length and width should be at least $6$ as for each edge there is a pair of neighbors with distance $3$. To get a finite one,  the vertices on parallel boundaries can be connected in the following ways. 
\begin{figure}[H]
\centering
\includegraphics[scale=0.3]{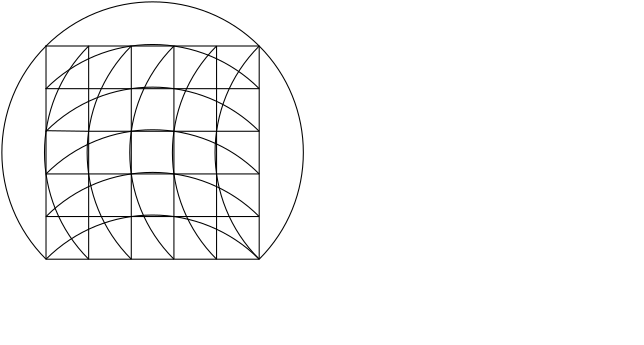}
\includegraphics[scale=0.3]{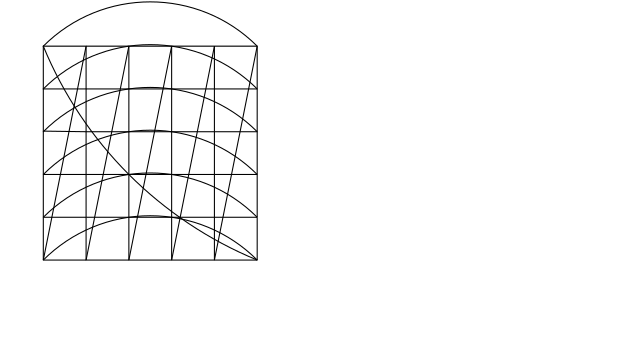}
\includegraphics[scale=0.55]{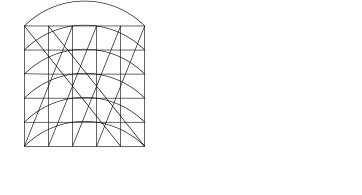}
\includegraphics[scale=0.55]{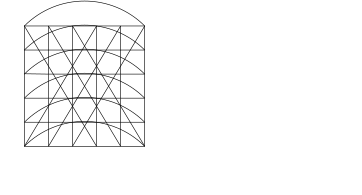}
\caption{ Finite $4$-regular Ricci-flat graphs of ``lattice type"}
\end{figure}
 \end{itemize}

\end{itemize}

Case 2: All edges $(0, 2), (0, 3), (1, 5), (1, 6)$ are ``Type 6c".  Then the $C_4$ for edge $(0, 4)$ cannot pass through vertex $2$ and its neighbor, then $4\sim 5$ or $4\sim 6$. Wlog, let $4\sim 5$. Then $d(2, 7)=3$ for the edge $(0, 1)$. Consider  the $C_4$ for edge $(0,7)$, we have $7\sim 3$. See the following graph.

\begin{figure}[H]
\centering
\includegraphics[scale=0.4]{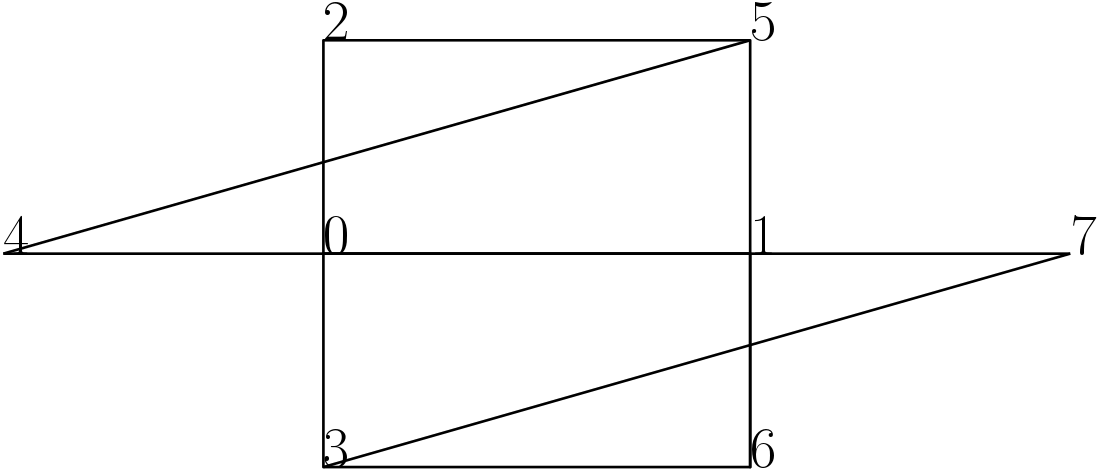}
\end{figure}
 Then $d(2, 6)=3, d(4, 6)=3$. Then consider the edge $(1, 5)$, since we have $d(2, 7)=d(2, 6)=d(4, 7)=d(4, 6)=3$,   then we need the fourth neighbor of vertex $5$ to have distance $1$ to either vertex $6$ or vertex $7$. Let $5\sim 9\sim 6$.  Then edge $(1, 6)$ shares two $C_4:=1-6-3-0-1, 1-6-9-5-1$, a contradiction.

\end{proof}

\begin{theorem}\label{miss24}
Let $G$ be a  Ricci-flat graph  that contains two $C_4$s sharing one edge and all vertices on any $C_4$ have degree $4$.  Assume G is not 4-regular, G is isomorphic to Figure 32. 
\end{theorem}

\begin{proof}
Start with the following structure where the edge $(0, 1)$ shares two $C_4$s and  $d(2)=d(5)=d(3)=d(6)=4$. Note $d(4),  d(7)\neq 4$, otherwise, we would get a $4 $-regular Ricci flat graph, thus $d(4)=d(7)=2$. For edge $(2, 5)$ to be in two $C_4$'s, let $2\sim 8\sim 9\sim 5$. As vertices $2, 5$ cannot be adjacent to any existing vertices, let $2\sim 9, 5\sim 11$, similarly, $d(10)=d(11)=2$. Observe that $4\not\sim 10$, otherwise $d(4)=d(10)=4$. However vertex $4$ cannot not be adjacent to new vertex, thus it has be adjacent to vertex $9$ considering edge $(0, 2)$.  Similarly, $7\sim  8, 10\sim 6, 11\sim 3$. Continue with this process, if vertex $8$ is adjacent to a new vertex called $a$, then $a$ must be adjacent to vertex $11$, a contradiction to $d(11)=2$, thus vertex $8$ is adjacent to existing one and it can only  vertex $3$. Similarly, vertex $9\sim 6$. Now the resulting graph is done and it is a Ricci-flat graph. 
\begin{figure}[H]
\centering
\includegraphics[scale=0.4]{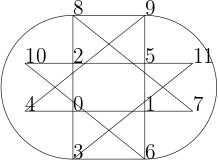}
\caption{A Ricci-flat graph}
\label{24C4C4miss1}
\end{figure}

\end{proof}

For the third case, let $G$ be a $4$-regular Ricci-flat graph in class $\mathcal{G}$ such that there exist two $C_4$s sharing one vertex and any two $C_4$s don't share an edge.
Then $G$ contains a $C_5$, since every edge on the $C_5$ must be in a $C_4$, there are three distinct substructures:
\begin{figure}[H]
\centering
\includegraphics[scale=0.3]{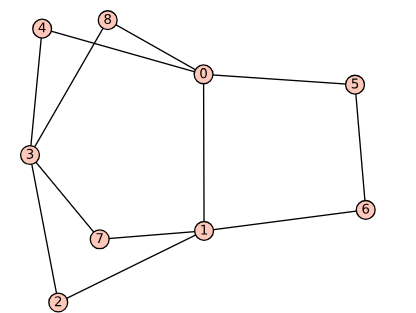}
\hfil
\includegraphics[scale=0.3]{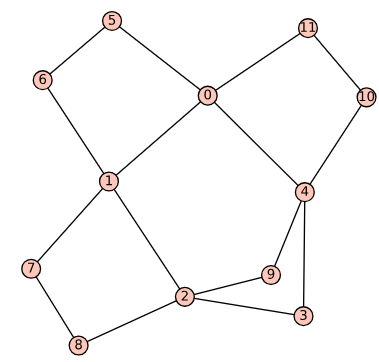}
\hfil
\includegraphics[scale=0.4]{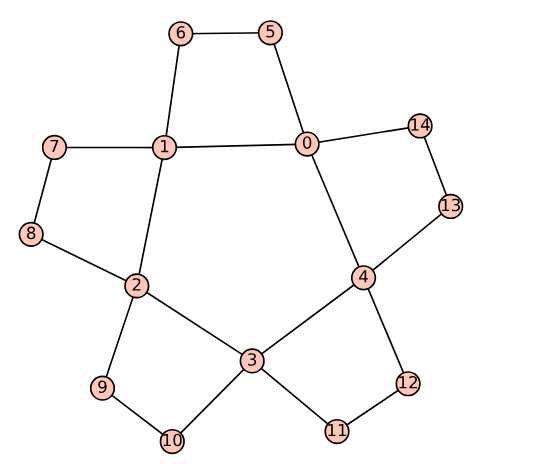}
\caption{Type A, Type B, Type C}
\end{figure}
Next, we consider each case. 
\begin{lemma}\label{thm:44noedgesharec4Typea}
Let $G$ be a $4$-regular Ricci-flat graph, then $G$ does not contain subgraph of ``Type A." 
\end{lemma}
\begin{proof}
We study the graph based on this structure. Note all vertices cannot be adjacent to any other vertices in the current subgraph in order  to avoid $C_4$s sharing on edge. Consider the edge $(0, 4)$, let $4\sim 9, 10$, wlog, then we need $d(1, 9)=2$ and $d(5, 10)=2$. For $d(1, 9)=2$, note $2\not\sim 9$, otherwise the edge $(3, 4)$ would share two $C_4$s. Similarly, $7\not\sim 9$ considering the edge $(3, 7)$. Then it must be  $6\sim 9$. 
Consider the vertex $8$, $8\sim 11, 12$ and $6\sim 11$. Then $d(6)=4$, consider the edge $(1, 6)$, wlog, we need $d(2, 9)=d(7, 11)=2$. Since $2, 7$ are not adjacent to any other vertices in the current subgraph, we need a common vertex $13$ such that $2\sim 13\sim 9$. Note $7\not\sim 13$, then let $7\sim 14\sim 11$. Let $2\sim 15, 7\sim 16$. 

Consider the edge $(2, 3)$, wlog, let $d(13, 4)=2$ and $d(15, 8)=2$. 
Assume $15\sim 11$, then consider the edge $(6, 11)$ in a $C_4$, assume $9\sim 15$, then $10\not\sim 15$ considering the edge $(9, 15)$.  It must be $10\sim 13$ for the $C_4$ that passes through vertices $4, 9, 10$. Consider the edge $(2, 13)$, since we have $d(1, 10)=3, d(3, 10)=2$, we also need $d(1, z)=2$ where $z$ is the fourth neighbor of vertex $13$, However, this cannot happen in the current subgraph. Thus When $15\sim 11$, $15\not\sim 9$. Then for the edge $(6, 11)$ in a $C_4$ that passes through vertex $9$, it must be the case $9\sim 14$. Then still we have  $10\sim 13$ for the $C_4$ that passes through vertices $4, 9, 10$. A contradiction again for the edge $(2, 13)$. Thus $5\not\sim 11$ for $d(15, 8)=2$, we must have $5\sim 12$. Similarly, we have $15\not\sim 9, 16\not\sim 9, 11$, and $16\sim 10$ for the edge $(3, 7)$. 

Now consider the edge $(1, 2)$. Since $d(6, 15)\neq 2$, we must have $d(0, 15)=2$ which implies  $15\sim 5$. 
Then for the edge $(1, 7)$, we need $d(0, 16)=2$ which implies  $16\sim 5$. 
Consider the $C_4$ that passes through vertices $5, 15, 16$, we need a common vertices for $15, 16$. Note  $15\not\sim 13, 9, 11$ and, $16\not\sim 14, 9, 11$. 
Assume $15\sim 10$, 
then consider the edge $(10, 15)$, we have $10\not\sim 13$. Thus we must have $13\sim 12$ for the 
$C_4$ that passes through vertices $2, 15, 13$. 
Consider the edge $(5, 6)$, either $d(15, 9)=2$ or $d(15, 11)=2$. If it is the former case, then $9\sim 10$ or $9\sim 12$, both would result appearances of $C_3$,  a contradiction. Thus it must be $d(15, 11)=2$ which implies $11\sim 10$. 
Then for the edge $(6, 11)$ in a $C_4$ that passes through vertex $9$, it must be the case $9\sim 14$.  It follows that $12\sim 14$ for 
the $C_4$ that passes through vertex $12, 8, 11$. However, the edge $(9, 14)$ shares two $C_4$s, a contradiction. 

We need a new vertex to be adjacent to $15, 16$. Let it be $15\sim 17\sim 16$. Then consider the edge $(2, 15)$, both  $d(3, 5)=3$ and $d(3, 17)=3$, a contradiction. 

Hence, there is no Ricci-flat graph based on ``Type A"

\end{proof}

\begin{lemma}\label{thm:44noedgesharec4Typeb}
Let $G$ be a $4$-regular Ricci-flat graph, then $G$ does not contain subgraph of ``Type B". 
\end{lemma}
\begin{proof}
First note vertices $3, 9$ are not adjacent to any existing vertices in the current subgraph. 
Consider the edge $(0, 1)$, if $d(4, 7)=2, d(11, 2)=2$, then we must have $11\sim 8, 7\sim 10$, then the edge $(7, 8)$ would share two $C_4$s, a contradiction. 
Thus we need $d(11, 7)=2, d(4, 2)=2$ for the edge $(0, 1)$. 
Let $7\sim 12\sim 11$.  Consider the edge $(0, 4)$, since $d(1, 9)=3$, we need $d(5, 9)=2$, note $5\not\sim 12$ and any other existing vertices, let $5\sim 14\sim 9$. Similarly, we need $d(6, 9)=2$ for the edge $(1, 2)$, since $6\not\sim 14$  and any other existing vertices, let $6\sim 15\sim 9$. Now consider the $C_4$ that passes through vertices $9, 14, 15$, we need a vertex as the common for vertices $14, 15$. Consider the edge $(2, 9)$, note $14\not\sim 7$, otherwise the cycle $C_5:=14-7-1-6-5-14$ would be a ``Type A". 
Thus $d(1, 14)=3$, then we need $d(8, 14)=2$, since $8\not\sim 5$, otherwise, there would be a $C_5:=5-0-1-2-8-5$ of ``Type A', also  $8\not\sim 12$ and 
any existing vertices in the current subgraph, we need a new vertex $16$ as the common of $14, 8$. Similarly, we need $d(10, 15)=2$.

Assume $10\sim 16\sim 15$, then we need a common vertex for $8, 10$ to form a $C_4$ considering the $C_5:=8-2-9-15-16-8$. Let $10\sim 17\sim 8$. Then consider the edge $(9, 14)$, since $d(2, 5)=3$, we need $d(2, z)=2$ where $z$ is the fourth neighbor of vertex $14$, then $z$ should be adjacent to neighbor of $2$. However all neighbors of $2$ have degree $4$, a contradiction. 
Thus for $d(10, 15)=2$, let $10\sim 17\sim 15$. Then for the $C_4$ that passes through vertices $9, 14, 15$, since $15\not\sim 16, 14\not\sim 17$, we need a new vertex, let $14\sim 13\sim 15$. 

For the edge $(5, 14)$ in a $C_4$ that passes through vertex $16$, assume $5\sim 17\sim 16$. Then the edge $(5, 17)$would share two $C_4:=5-17-16-14-5, 5-17--15-6-5$, a contradiction. Thus we need a new vertex as common for $5, 16$, 
let $5\sim 18\sim 16$. Similarly, for $d(6, 17)=2$, note both $6, 17$ cannot be adjacent to $18$ considering the edge $(5, 6)$. Thus let $6\sim 19\sim 17$. 

Consider the $C_4$ that passes through vertices $6, 15, 9, 17$, then  $19\not \sim 13$, otherwise $(6, 15)$ would share two $C_4$s.  Similarly, $18\not\sim 13$.

Now consider the fourth neighbor of vertex $7$: 
\begin{itemize}
\item Assume $7\sim 13$, then consider the edge $(1, 7)$, if $d(6, 13)=2$, then $13\sim 19$, a contradiction. Thus we need $d(0, 13)=2$, then $13\sim 11$. 
Consider the edge $(7, 8)$, note $d(16, 13)\neq 2$, otherwise $16$ must be adjacent to one of neighbors of vertex $13$, however, all neighbors have degree $4$. Similar, let $s$ represent the fourth neighbor of vertex $8$, we have $d(z, 13)=3$. Thus $d(8)\neq 4$. a contradiction. Actually, when $d(8)=3$, 
this Ricci-flat graph  is isomorphic to graph $G_4$ in Theorem \ref{34withonec41}. 
\begin{figure}[H]
\centering
\includegraphics[scale=0.34]{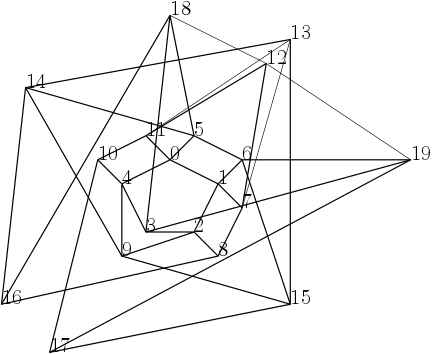} 
\caption{A Ricci-flat graph isomorphic to $G_4$}
\end{figure}

\item Assume $7\sim 17$, consider the $C_4$ that passes through vertices $17, 10, 7, 12$, then $10\sim 12$, a contradiction. 
\item Assume $7\sim 10$, consider the $C_4$ that passes through vertices $7, 10, 17, 12$, then $17\sim 12$. However the edge $(12, 17)$ would share two $C_4$s, with the second one $C_4:=17-12-11-10-17$, a contradiction.

\item Then $7\sim 20$, a new vertex.  For the edge $(1, 7)$, since $d(0, 20)\neq 2$, we need $d(6, 20)=2$ which implies $20\sim 19$. Consider the edge $(6, 19)$, we need $d(5, w)=2$ where $w$ is the fourth neighbor of $19$, then $w\sim 18$. Assume $w=12$, then there is a $C_4$ that passes through vertices $7, 12, 19, 20$. However, consider the edge $(7, 20)$, any new neighbor of vertex $20$ has distance $3$ from vertex $1$, a contradiction.

\end{itemize}
We have considered all possible neighbors for vertex $7$. Thus there is no $4$-regular Ricci-flat graph of ``Type B ".

\end{proof}
Now we consider all graphs in which every $C_5$ must be ``Type C". We have the following result. 
\begin{theorem}\label{thm:44noedgesharec4Typec}
Let $G$ be a $4$-regular Ricci-flat graph that contains a subgraph isomorphic to ``Type C". And $G$ is isomorphic to the graph showing in Figure \ref{44C4Ricci-flat3}.  
\end{theorem}
\begin{proof}
We start from such a subgraph of ``Type C" in $G$ by consider the neighborhood of each edge on the $C_5$, wlog, for the edge $(0, 1)$, there are two cases: 
\begin{itemize}
\item [Case 1:] $d(4, 7)=2, d(14, 2)=2$. For $d(4, 7)=2$, we assume $13\sim 7$. However the $C_5:=7-1-0-14-13-7$ does not satisfy ``Type C" Thus it must be  $12\sim 7$. Then consider the $C_5:=12-4-0-1-7-12$, we need a $C_4$ that passes through edge $(12, 7)$ with new vertices added, let it be  $12\sim 16\sim 15\sim 7$. For $d(14, 2)=2$, similarly, we need $14\sim 9$. 

Now consider the edge $(1, 7)$, there are also two cases: 
\begin{itemize}

\item[Case 1a:] $d(0, 15)=2$, $d(6, 12)=2$. 
   Case 1a: By the ``Type C", it must be $14\sim 15$ and $6\sim 11$. Consider the  $C_5:=1-6-11-12-7-1$, we need two new vertices to form a $C_4$ with $6, 11$, let them be $6\sim 17\sim 18\sim 11$. Consider
   the $C_5:=14-0-1-7-15-14$, to satisfy ``Type C", the edge $(14, 15)$ must be in a $C_4$ that pass through vertex $9$, there are two case: $9\sim 17\sim 15$ or $15, 9$ are adjacent to a new vertex $19$.
    Assume the former.
   Consider the edge $(0, 14)$, there are two cases: Case 1a1: $d(1, 15)=2$, $d(5, 9)=2$; Case 1a2: $d(1, 9)=2, d(5, 15)=2$. For Case 1a1, we need $5\sim 10$, for Case 1a2, we need $5\sim 16$. Assume Case 1b1 is true, then for the edge $(5, 6)$, we need either $5\sim 16$ or $5\sim 19$ where $19$ is a new vertex such that $19\sim 18$. However the latter case would contradict to the ``Type C" for the $C_5:=18-11-6-5-19-18$. Thus for the edge $(0, 14)$, Case 1a2 is always true, that is $5\sim 16$. Assume $5\not\sim 10$. Let $5\sim 19$, a new vertex. Then we need $16\sim 18$ and $19\not\sim 18$. However, then $d(19, 17), d(19, 11)>2$, a contradiction. Thus for the edge $(0, 14)$, we have $5\sim 16$, and $5\sim 10$. 
   
   Note  $16\sim 18$ by ``Type C". 
   By similar reasons, for the edge $(2, 8)$, we need $8\sim 13, 8\sim 18$. The for the edge $(8, 18)$
 to  be in a $C_4$, the only choice is passing through $13$, let $18\sim 19\sim 13$. Then consider the edge $(4, 13)$, we must have $19\sim 16$. Then for the edge $(5, 10)$ to be in a $C_4$, the only choice is $16\sim 19$. Done. 
\begin{figure}[H]
\centering
 \includegraphics[scale=0.35]{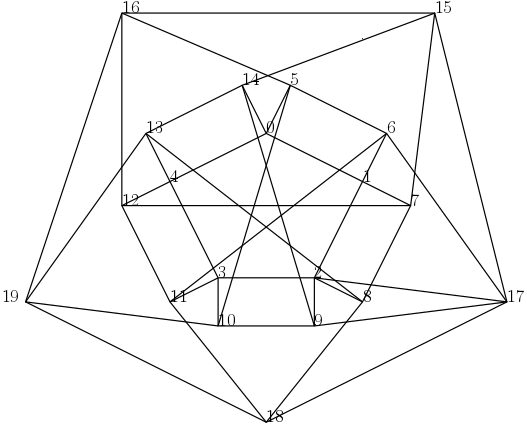}   
 \caption{Ricci-flat graph}   
 \label{44C4Ricci-flat3}  
\end{figure}

\item[Case 1b:] $d(0, 12)=2, d(6, 15)=2$. This case is included by above case.
   \end{itemize}

\item[Case 2:]  $d(4, 2)=2, d(14, 7)=2$. Note we consider the same situation for all the other edges on the $C_5$.  For $d(14, 7)=2$, we need a new vertex $15$. Similarly, we must have $6\sim 16\sim 9, 8\sim 17\sim 11, 10\sim 18\sim 13, 12\sim 19\sim 5$. 
Now, we conclude that except above Ricci-flat graph, any other $4$-regular Ricci-flat graph in class $\mathcal{G}$ such that no edges share two $C_4$s  must contain a subgraph that isomorphic to the following one:
\begin{center}
\includegraphics[scale=0.3]{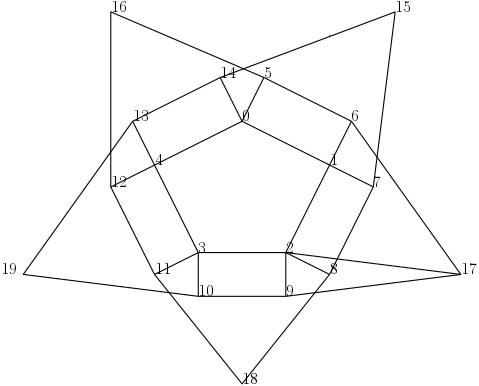}
\end{center}

Consider the $C_5:=5-0-4-12-19-5$, then the edge $(12, 19)$ must be in a $C_4$. Note $12\not\sim 13, 14, 5, 6, 10, 8$. 
Assume $12\sim 9$. Then $19\sim 16$ which lead to the edge $(5, 6)$ share two $C_4$, a contradiction. 
Assume $12\sim 7$. Then $19\sim 15$. Now consider the edge $(1, 7)$, since $d(0, 12)=3, d(0, 15)=2, d(2, 12)\geq 2, d(6, 12)\geq 2$, then we need $d(6, 12)=2$, thus $6\sim 11$. Similarly, $5\sim 10, 19\sim 18, 14\sim 9, 15\sim 16,  16\sim 17,   13\sim 8, 17\sim 18$. Done. The resulted graph is actually isomorphic to previous one. 

Now we consider another case for the circle $C_5:=5-0-4-12-19-5$. That is, we need new vertices to form a $C_4$ for the edge $(5, 19)$ Let $5\sim 20\sim 21\sim 19$.  Note we should not have $15\sim 19$ for this case, otherwise, it would result a same graph as previous one. Consider the edge $(0, 5)$,  assume $d(4, 20)=2, d(14, 19)=2$, then for $d(4, 20)=2$, either $20\sim 13$ or $20\sim 12$.  however both assumptions do not guarantee the ``Type C" for $C_5:=20-13-4-0-5-20$ or $C_5:=20-12-4-0-5-20$. 
Then it muse be the case $(4, 19)=2, d(14, 20)=2$ for the edge $(0, 5)$. For $d(14, 20)=2$, note $14\not\sim 21$ since the $C_5:=14-21-20-5-0-14$ does not satisfy ``Type C". Note $14\not\sim 16$ or $14\not\sim 17$  as the $C_5:=14-0-5-6-16-14$ and $C_5:=14-15-7-8-17-14$ cannot satisfy ``Type C". It is easy to see that $14$ is not adjacent to any other existing vertices in the current subgraph. Thus we need a new vertex for $d(14, 20)=2$, let $14\sim 22\sim 20$. Then consider the edge $(0, 14)$, note $d(1, 22)\neq 2$, otherwise $22\sim 6$ or $22\sim 7$, both cases would generate $C_5$s that do not satisfy ``Type C". Thus we need $d(1, 15)=2$ and $d(5, 22)=2$. The current subgraph meet this requirement. Now consider the edge $(4, 15)$ which must be in a $C_4$ that passes through vertex $22$. Note $15\not\sim 20$ for this requirement as $C_5:=14-15-20-5-0-14$ cannot satisfy ``Type C". Thus let $22\sim 23\sim 15$.

Next, we consider the fourth neighbor of vertex $6$. Consider the edge $(5, 20)$, we need $d(6, w)=2$ where $w$ is the fourth neighbor of vertex $20$, thus $20\not\sim 16$, and  $6\not\sim 22$.  
Note $6\not\sim 23$, otherwise 
the $C_4$ for the edge $(6, 23)$ must pass through the fourth neighbor of vertex $23$ which is adjacent to vertex $16$. Let $23\sim 24\sim 16$. Then consider the edge $(15, 7)$ on the $C_5:=1-6-23-15-7-1$, note $7\not\sim 22, 24$, since these two vertices are on the $C_4$s attached for this cycle. Note $7\not\sim 21$ as we cannot guarantee the edge $(7, 21)$ any more. It is easy to see that $7$ is not adjacent to any other existing vertices in the current subgraph. Thus let $7\sim 25\sim 26\sim 15$. Consider the edge $(5, 19)$, we need the fourth neighbor of vertex $19$ to have distance $2$ form vertex $6$. Consider the edge $(5, 20)$, we also need the fourth neighbor of vertex $20$ to have distance $2$ form vertex $6$, since both the fourth neighbors should be adjacent to vertex $16$ as the there has been $4$ neighbors for vertex $23$, then $d(16)\geq 5$, a contradiction.  It is easy to see that $6$ is not adjacent to any other existing vertices in the current subgraph.

Thus we need a new vertex $24$ as the fourth neighbor of vertex $6$. 
 Note vertex $7$ is not adjacent to any existing vertices, we need a new vertex $26$ as the fourth neighbor of it. Same reason as why $20\sim 22$, we have $24\sim 26$. 
 For the edge $(6, 16)$ to be in a $C_4$ that passes through vertex $24$, assume $24\sim 23\sim 16$. Then the  $C_5:=23-24-26-7-15-23$ cannot satisfy ``Type C". Thus we need a new to form $C_4$ for the edge $(6, 16)$, let $24\sim 25\sim 16$. And we need a new neighbor for vertex $15$ to from a $C_4:=7-26-27-5$. 
 
 We now consider the  fourth neighbor of vertex $8$ and $9$. Assume $8\sim 20$, then for the edge $(8, 17)$ in a $C_4$ that pass through vertex $20$, we must have $17\sim 22$. However, for the edge $(8, 20)$, since $d(2, 5)=3$, we need $d(2, 21)=2$ and $d(5, 9)=2$, which implies $9\sim 21$ and $9\sim 19$, then $d(9)\geq 5$, a contradiction. Similar analysis shows that $8\not\sim 22$, $9\not\sim 20, 9\not\sim 22$. And it is easy to check vertices $8$ and $9$ are  not adjacent to any other  existing vertices in the current subgraph. Thus let $8\sim 28, 9\sim 30$ and we also have $28\sim 30$ by similar reasons for $20\sim 22, 24\sim 26$. Now we analysis the common vertex for $28, 17$. We consider the vertex with degree $2$ in the current subgraph. Assume $28\sim 21\sim 17$. Then we need either $d(22, 17)=2$ or $d(22, 28)=2$ considering the edge $(20, 21)$, and $d(z, 28)=2$ where $z$ is the fourth neighbor of vertex $19$. Since $d(28)=3$ in the current subgraph and the vertices $19, 22$ cannot share a common vertex, we actually need $d(22, 17)=2$. Note $17\not\sim 20$ and $22\not\sim 11$, thus we need a new vertex $29$ such that $22\sim 29\sim 17$, then consider the edge $(8, 17)$, we need $d(7, 29)=2$, then $29\sim 26$. Since we need a $C_4$ for the edge $11, 17$ that passes through vertex $29$ and a $C_4$ for the edge $24, 26$ that pa passes through vertex $29$, then $d(29)\geq 5$, a contradiction. 
 One can also check that $28, 17$ are not both adjacent to vertex $23, 27, 25$. Thus we need new vertex.  Let $8\sim 18\sim 29\sim 17$, similarly, $9\sim 30\sim 31\sim 16$.  

Similar argument will give us $11\sim 32\sim 33\sim 17$, $10\sim 34\sim 35\sim 18$, $32\sim 34$, and $13\sim 36\sim 37\sim 18, 12\sim 28\sim 39\sim 19, 36\sim 38$. See the following structure: 
\begin{center}
\includegraphics[scale=0.3]{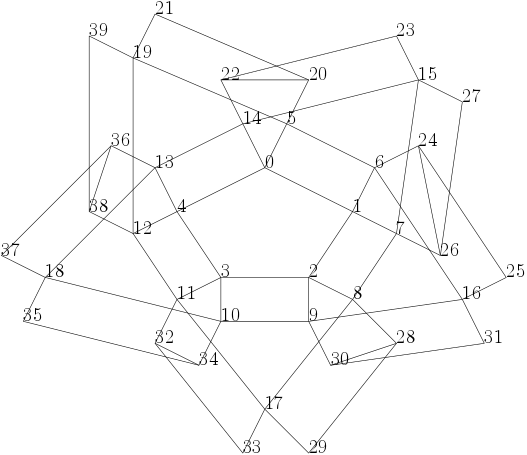}  
\end{center}

Now observe the edge $(5, 6)$, since $d(9, 16)$ cannot be $2$, we must have $d(9, 14)=2$ and $d(20, 6)=2$. Note $24\not\sim 21$, as it would generate the $C_5:=19-21-24-6-5-19$ which does not satisfy ``Type C". Then it must be $24\sim 39$, similarly $20\sim 31, 26\sim 33, 28\sim 23, 30\sim 37, 34\sim 25, 32\sim 21, 38\sim 29, 36\sim 27, 22\sim 35$. 
Consider the $C_4$ for the edge $(20, 22)$, we must have $31\sim 35$, similarly, $33\sim 39$, $23\sim 37, 21\sim 25, 27\sim  29$.  
\begin{center} 
\includegraphics[scale=0.25]{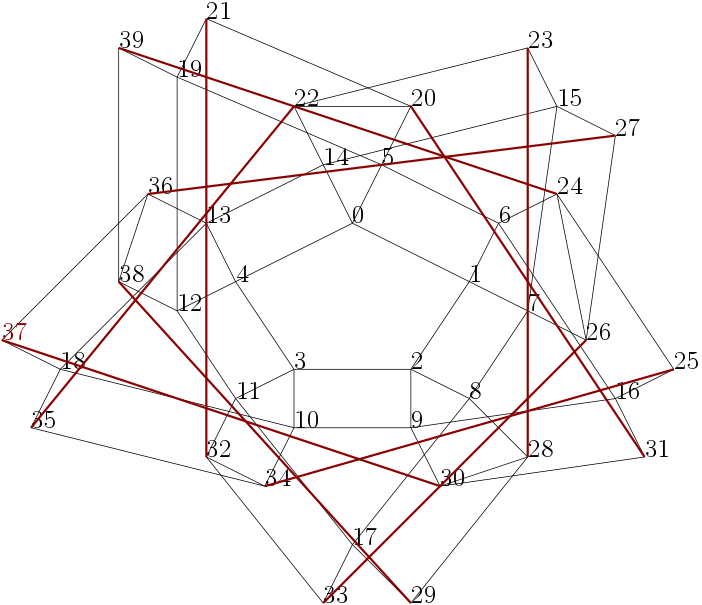}
\includegraphics[scale=0.25]{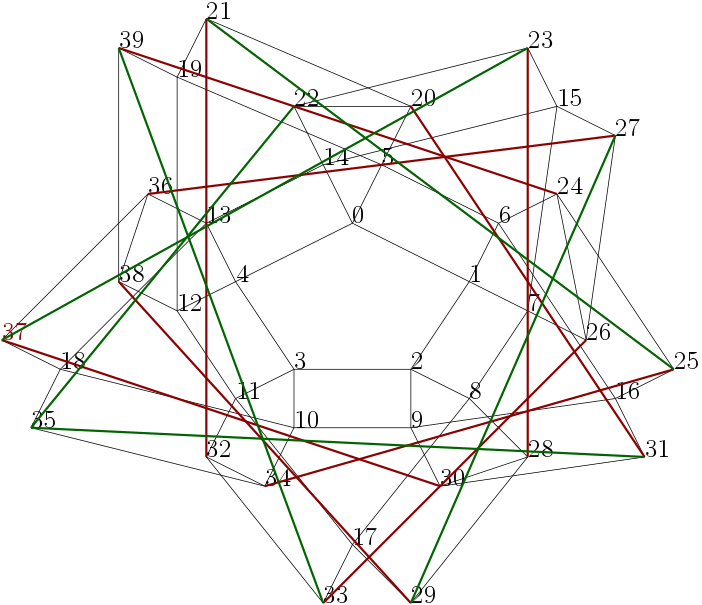}
\end{center}
Now all vertices have degree $4$, however we cannot guarantee the edge $(20, 21)$, since $d(22, 32)$ cannot be $2$. 
Thus there is no Ricci-flat graph for this case. 

\end{itemize}

\end{proof}
\bibliographystyle{plain}
\bibliography{Riccicita}

\end{document}